\documentclass[a4paper,12pt]{amsbook}

\usepackage[latin1]{inputenc}
\usepackage[english]{babel}
\usepackage{latexsym}
\usepackage{amssymb}
\usepackage{amsmath}
\usepackage{amsthm}
\usepackage[all]{xy}
\usepackage[dvips]{graphicx}

\newtheorem*{teos}{Theorem}

\newtheorem{teo}{Theorem}[section]
\newtheorem{lem}[teo]{Lemma}
\newtheorem{prop}[teo]{Proposition}
\newtheorem{cor}[teo]{Corollary}
\theoremstyle{definition}

\newtheorem{dhef}[teo]{Definition}
\newtheorem{exe}[teo]{Example}
\newtheorem{oss}[teo]{Remark}

\newcommand{\epsi}{\varepsilon}
\newcommand{\lrg}{\longrightarrow}

\newcommand{\N}{\mathbb{N}}
\newcommand{\Z}{\mathbb{Z}}
\newcommand{\C}{\mathbb{C}}
\newcommand{\K}{\mathbb{K}}

\newcommand{\kepsi}{\mathbb{K}[\varepsilon]}

\newcommand{\shA}{\mathcal{A}}

\newcommand{\shO}{\mathcal{O}}

\newcommand{\de}{\partial}
\newcommand{\debar}{\overline{\partial}}

\newcommand{\Spec}{\operatorname{Spec}}

\newcommand{\Aut}{\operatorname{Aut}}
\newcommand{\MC}{\operatorname{MC}}
\newcommand{\Def}{\operatorname{Def}}

\newcommand{\Ext}{\operatorname{Ext}}

\newcommand{\Hom}{\operatorname{Hom}}
\newcommand{\Htp}{\operatorname{Htp}}

\newcommand{\image}{\operatorname{Im}}
\newcommand{\coker}{\operatorname{coker}}

\newcommand{\contr}{{\mspace{1mu}\lrcorner\mspace{1.5mu}}}
\newcommand{\bl}{\boldsymbol{l}}
\newcommand{\bi}{\boldsymbol{i}}

\newcommand{\Fun}{\mathbf{Fun}}
\newcommand{\Art}{\mathbf{Art}}
\newcommand{\Set}{\mathbf{Set}}
\newcommand{\Com}{ \mathbf{C}}

\newcommand{\mc}{\operatorname{MC}_{(h,g)}}
\newcommand{\cil}{\operatorname{C}_{(h,g)}}
\newcommand{\hil}{\operatorname{H}_{(h,g)}}

\begin{document}

\begin{titlepage}

\advance\textheight by \headheight

\advance\textheight by \headsep

\advance\textheight by \footskip

\headheight=0pt

\headsep=0pt

\footskip=0pt

\begin{center}
\par\vspace {25mm}
{\Large\sc Universit\`a  \ degli \  Studi \  di \ Roma \ }\\[2ex]
{\Large\sc ``La Sapienza''}\\[4ex]
{\large Facolt\`a di Scienze Matematiche Fisiche e Naturali}\\[1ex]
{\large Dottorato di Ricerca in Matematica }\\[1ex]
 {\large XVIII Ciclo}\\

\end{center}

\begin{center}
\par\vspace {35mm}

{\Huge\bf Differential Graded } \\

\bigskip

{\Huge\bf  Lie Algebras }  \\

\bigskip

{\Huge\bf and }\\

\bigskip

{\Huge\bf  Deformations of}\\

\bigskip

{\Huge\bf Holomorphic Maps}

\par\vspace{45mm}
\end{center}

\par\vspace{30mm}

\noindent

\raggedright {\large \bf Dottoranda} \hfill { \large \bf Relatore}

{\large \it Donatella Iacono} \hfill
\large Prof. {\large \it Marco Manetti}\\

\end{titlepage}

\newpage

{\ }

 \vspace{75mm}

\begin{center}

\Large{\emph{This thesis is devoted}}

\emph{ to the memory }

\emph{of my grandmother Mina.}

\end{center}

\tableofcontents

\chapter*{Introduction}

In the last fifty  years   deformation  theory has played an
important role in algebraic and complex geometry.

The study of small deformations of the complex structures of
complex manifolds started with the works of K.~Kodaira and
D.C.~Spencer \cite{bib Kodaira SpencerII} and   M.~Kuranishi
\cite{bib kuranishi}.

\medskip

Then A.~Grothendieck \cite{bib Grothendieck}, M.~Schlessinger
\cite{bib Schlessinger } and M.~Artin \cite{bib Artin} formalized
this theory translating it into a functorial language. The idea
was that with an infinitesimal deformation of a geometric object
we can associate a deformation functor of Artin rings: that is, a
functor from the category $\Art$ of local Artinian $\C$-algebras
(with residue field $\C$) to the category $\Set$ of sets, that
satisfies some extra conditions (see Definition~\ref{defin
deformation funtore}).

\medskip

A modern approach to   deformation  theory is via differential
graded Lie algebras (DGLA for short).

The guiding principle is the idea  due to P.~Deligne, V.~Drinfeld,
D.~Quillen and M.~Kontsevich (see \cite{bib kontsevich}) that \lq
\lq in characteristic zero  every deformation problem is
controlled by a differential graded Lie algebra".

\medskip

Inspired by this principle, the aim of this thesis is to follow
the modern approach to study the infinitesimal deformations of
holomorphic maps of compact complex manifolds.

\medskip

More precisely,  a DGLA is a differential graded vector space with
a structure of graded Lie algebra, plus some   compatibility
conditions between the differential and the bracket (of the Lie
structure) (see Definition~\ref{definizio DGLA}).

Moreover, using the solutions of the Maurer-Cartan equation
$dx+\displaystyle\frac{1}{2}[x,x]=0$  it is  well known how we can
associate with a DGLA $L$ a deformation functor $\Def_L$. Written
in detail:
$$
\Def_L:\Art \lrg \Set,
$$
$$
\Def_L(A)=\frac{\{ x \in L^1\otimes m_A \ |\ dx+
\displaystyle\frac{1}{2} [x,x]=0 \}}{gauge},
$$
where $m_A$ is the maximal ideal of $A$ and   the gauge
equivalence is  induced by the gauge action $*$ of $\exp(L^0
\otimes m_A)$ on the set of solutions of the Maurer-Cartan
equation (see Definition~\ref{dhef funtore defo DLGA DEF_L}).

Then the idea of the principle is that we can define a DGLA $L$
(up to quasi-isomorphism) from the geometrical data of the
problem, such that the deformation functor $\Def_L$ is isomorphic
to the deformation functor  associated with the geometric problem.

We note that it is easier to study a deformation functor
associated with a DGLA but, in general,  it is not an easy task to
find the right DGLA (up to quasi-isomorphism) associated with the
problem.

\medskip

A first example, in which the associated DGLA is well understood,
is the case of deformations of  complex manifolds.

Let $X$ be a compact complex manifold. Then $X$ is obtained by
gluing together  a finite number of polydisks in $\C^n$. The
fundamental idea of K.~Kodaira and D.C.~Spencer is that \lq\lq a
deformation of $X$ is regarded as the gluing of the same polydisks
via different identifications" (see \cite[pag.~182]{bib kodaira
libro}) .

Translating it into a functorial language we define, for each $A
\in\Art$, an infinitesimal deformation of $X$ over $\Spec(A)$ as a
commutative diagram
\begin{center}
$\xymatrix{X \ar[r]^i\ar[d] & X_A\ar[d]^{\pi} \\
           \Spec(\C) \ar[r]^{a} & \Spec(A),   \\ }$
\end{center}
where $\pi$ is a proper  flat  holomorphic map and $X$ coincides
with the restriction of $X_A$   over the closed point of
$\Spec(A)$ (see Definition~\ref{defin funtore DEF_X}).  Moreover,
we can give the notions of isomorphism  and of trivial deformation
($X_A\cong X\times \Spec(A)$).

Thus we can define the functor associated with the infinitesimal
deformations of $X$:

$$
\Def_X : \mathbf{Art}\to\mathbf{Set},
$$

$$
\Def_X(A)=
\left\{\begin{array}{c}  \mbox{  isomorphism   classes   of}\\
 \mbox{ infinitesimal   deformations}\\
 \mbox{ of  $X$   over   $\Spec(A)$}
\end{array} \right\}.
$$

Therefore, following the idea of the principle,  we are looking
for a DGLA $L$ such that $\Def_L \cong \Def_X$.

Let $\Theta_X$ be the holomorphic tangent bundle of $X$ and
consider the sheaf  $\shA^{0,*}_X(\Theta_X)$ of the differential
forms of $(0,*)$-type, with coefficients in $\Theta_X$.

Then we  define the Kodaira-Spencer algebra $A^{0,*}_X(\Theta_X)$
of $X$ as the graded vector space of global sections of the sheaf
$\shA^{0,*}_X(\Theta_X)$.

Considering the (opposite) Dolbeault differential and the bracket
of vector fields, it is possible to  endow  $A^{0,*}_X(\Theta_X)$
with a natural structure of differential graded Lie algebra (see
Definition~\ref{def kodaira spencer algebra}).

The DGLA $A^{0,*}_X(\Theta_X)$  controls the problem of
infinitesimal deformations of $X$ (see Theorem~\ref{teo def_k
=Def_X}):
\begin{teos}[A]
Let $X$ be a complex compact manifold and $A^{0,*}_X(\Theta_X)$
its Kodaira-Spencer algebra.  Then there exists an isomorphism of
functors
$$
 \Def_{A^{0,*}_X(\Theta_X)} \lrg \Def_X.
$$
\end{teos}

In this case  the correspondence is clear between the solutions of
the Maurer-Cartan equation and the infinitesimal deformations of
$X$, such that the gauge equivalence corresponds to the
isomorphism of deformations. In particular, a solution of the
Maurer-Cartan equation is gauge equivalent to zero if and only if
it induces a trivial deformation of $X$.

\bigskip

The next natural problem is to investigate the embedded
deformations of a submanifold in a fixed manifold.

Very recently, M.~Manetti  in \cite{bib ManettiPREPRINT} studies
this problem using the approach via DGLA.

In his   work, more than proving the existence of a DGLA which
controls this geometric problem, M.~Manetti develops some
algebraic tools related to the DGLA.

More precisely, he describes a general construction to define a
new deformation functor associated with a  morphism of DGLAs.

Given a morphism of differential graded Lie algebras
$$
h:L \lrg M
$$
he defines the functor
$$
\Def_h: \Art \lrg \Set,
$$
$$
\Def_{h}(A)=
$$
$$
\frac{\{(l,m) \in (L^1 \otimes m_ A)  \times
 (M^0 \otimes m_A)  |\  dx +\frac{1}{2}[x,x]=0 \mbox{ and }
e^m*h(l)= 0 \}}{gauge},
$$
where this gauge equivalence is a generalization of the previous
one (see Remark~\ref{oss DEf(hg) reduces to DEF_h e DEF L}). This
new functor is an extension of the functor associated with a
single DGLA: by choosing $h=M=0$, $\Def_h$ reduces to $\Def_L$.

Moreover, using path objects (see Example~\ref{exe definizio
M[t,dt]})  he shows that for every choice of $L,M$ and $h$ there
exists a DGLA $H$ such that $\Def_H \cong \Def_h$.

In particular, this implies that if a deformation functor
associated with a geometric problem is isomorphic to $\Def_h$, for
some $h:L \lrg M$, then automatically we know the existence of a
DGLA that controls the problem and we have an explicit description
of it.

This is what M.~Manetti does in \cite{bib ManettiPREPRINT}: by
suitably choosing    $L,M$ and $h$ he proves that there exists an
isomorphism between the functor $\Def_h$ and the functor
associated with the infinitesimal deformations of a submanifold in
a fixed manifold. Actually, let $X$ be a compact complex manifold
and $Z$ a submanifold.  The infinitesimal emebedded  deformations
of $Z$ can be interpreted as the deformations of the inclusion map
$i:Z \hookrightarrow X$ inducing the trivial deformation on $X$.

Consider the Kodaira-Spencer DGLA $A^{0,*}_X(\Theta_X)$ of $X$ and
the differential graded Lie subalgebra $A_X^{0,*}(\Theta_X(-log\,
Z))$ defined by the following exact sequence
$$
0 \lrg A_X^{0,*}(\Theta_X(-log\, Z)) \lrg A_X^{0,*}(\Theta_X)
  \lrg A_Z^{0,*}(N_{Z|X})\lrg 0,
$$
where $N_{Z|X}$ is the normal bundle of $Z$ in $X$ (see
Section~\ref{sezio submanifold defz L'}).

We have already observed that the DGLA $A^{0,*}_X(\Theta_X)$
controls the infinitesimal deformations of $X$, whereas  the DGLA
$A_X^{0,*}(\Theta_X(-log\, Z))$   controls the infinitesimal
deformations of the pair $Z\subset X$ (each solution of the
Maurer-Cartan equation in $A_X^{0,*}(\Theta_X(-log\, Z))$  defines
a deformation of both $Z$ and  $X$).

Fix $M=A^{0,*}_X(\Theta_X)$, $L=A_X^{0,*}(\Theta_X(-log\, Z))$ and
$h$ the inclusion:
$$
h:A_X^{0,*}(\Theta_X(-log\, Z)) \hookrightarrow
A^{0,*}_X(\Theta_X).
$$

Then it is clear how we can associate with each element $(l,m) \in
\Def_h$ a deformation of $Z$ in  $X$, with $X$ fixed: the
infinitesimal deformation of $Z$ is the one corresponding to the
Maurer-Cartan solution  $l \in A_X^{0,*}(\Theta_X(-log\, Z))$ and
it induces a trivial deformation of $X$, since we are requiring
that $h(l)$ is gauge equivalent to zero in $ A^{0,*}_X(\Theta_X)$.

\medskip

These new ideas developed by M.~Manetti are of fundamental
importance for  this thesis which, in some sense, can be
considered a generalization of them. Actually, we extend these
techniques to study not only the deformations of an inclusion but,
in general, the deformations of holomorphic maps.

\smallskip

These deformations were first studied from the classical point of
view (no DGLA) by E.~Horikawa  \cite{bib HorikawaI} and \cite{bib
Horikawa III}, M.~Namba  \cite{bib Namba}  and by Z.~Ran \cite{bib
RAn mappe}.

\medskip

More precisely, let $f:X \lrg Y$ be a  holomorphic map of compact
complex manifolds.

There are several aspects of deformations of $f$: we can deform
just the map $f$ fixing both $X$ and $Y$, we can allow to deform
$f$ and $X$ or $Y$ or, more in general, we can deform everything:
the map $f$, $X$ and $Y$.

The infinitesimal deformations of $f$, with fixed domain and
target, can be interpreted as infinitesimal deformations of the
graph of $f$ in the product $X \times Y$, with $X\times Y$ fixed
(see Section~\ref{sezio fixed target domain}). Therefore, we are
considering  infinitesimal deformations of a submanifold in a
fixed manifold and so  the DGLA approach to this case  is
implicitly included in M. Manetti's work \cite{bib
ManettiPREPRINT}.

Next, we turn our attention to the general case in which we deform
$f$, $X$ and $Y$.

In a functorial language, for each $ A \in \Art$, we define an
infinitesimal deformation of $f$  over $\Spec(A)$  as a
commutative diagram
\begin{center}
$\xymatrix{X_A\ar[rr]^{F} \ar[dr]_\pi &  & Y_A
\ar[ld]^\mu \\
            & S, &  \\ }$
\end{center}
where $S=\Spec(A)$, $(X_A,\pi,S)$ and $(Y_A,\pi,S)$ are
infinitesimal deformations of $X$ and $Y$ respectively, and $F$ is
a holomorphic  map that restricted to the fibers over the closed
point of $S$ coincides with $f$.

In this case too we can give the notions of isomorphism and of
trivial deformation.

Then we can define the   functor of infinitesimal deformations of
a  holomorphic map $f:X \lrg Y$:
$$
\Def(f): \Art \lrg \Set,
$$
$$
\Def(f)(A)= \left\{
\begin{array}{c} \mbox{ isomorphism  classes  of}\\
\mbox{ infinitesimal   deformations   of }\\
f  \mbox { over }  \Spec(A) \\
\end{array} \right\}.
$$

\medskip

Let $\Gamma$ be the graph of $f$ in $X \times Y$. An infinitesimal
deformation of the map $f$ can be interpreted as an infinitesimal
deformation of $\Gamma $ in $X \times Y$, such that the induced
deformation of the product $ X \times Y$ is a product of
deformations of $  X $ and $ Y$. In general not all the
deformations of a product are products of deformations, as it was
showed by  K. Kodaira and D.C Spencer (see Remark~\ref{oss
DEF(XxY)> DEF(X) x DEF(Y)}).

\medskip

Consider the Kodaira-Spencer algebra $A^{0,*}_{X \times
Y}(\Theta_{X \times Y}) $  of the product $X \times Y$ and the
differential graded Lie subalgebra $A_{X \times Y}^{0,*}(\Theta_{X
\times Y}(-log\, \Gamma))$ defined by the following exact sequence
$$
0 \lrg A_{X \times Y}^{0,*}(\Theta_{X \times Y}(-log\, \Gamma))
\lrg A_{X \times Y}^{0,*}(\Theta_{X \times Y})
  \lrg A_\Gamma^{0,*}(N_{\Gamma|X \times Y})\lrg 0,
$$
where $N_{\Gamma|X \times Y}$ is the normal bundle of the graph
$\Gamma$ in $X\times Y$ (see Section~\ref{sezio submanifold defz
L'}).

As before, we know that $A_{X \times Y}^{0,*}(\Theta_{X \times
Y})$ controls the infinitesimal deformations of $X\times Y$ and
$A_{X \times Y}^{0,*}(\Theta_{X \times Y}(-log\, \Gamma))$
controls the infinitesimal deformations of the pair $\Gamma\subset
X \times Y$ (each solution of the Maurer-Cartan equation in $A_{X
\times Y}^{0,*}(\Theta_{X \times Y}(-log\, \Gamma))$ defines a
deformation of both $\Gamma$ and  $X\times Y$).

Fix $M=A_{X \times Y}^{0,*}(\Theta_{X \times Y})$, $L=A_{X \times
Y}^{0,*}(\Theta_{X \times Y}(-log\, \Gamma))$ and $h$ the
inclusion:
$$
h:A_{X \times Y}^{0,*}(\Theta_{X \times Y}(-log\, \Gamma))
\hookrightarrow A_{X \times Y}^{0,*}(\Theta_{X \times Y}).
$$

\smallskip

In the general case, it does not suffice  to consider just the
DGLA $L$ or the morphism  $h :L \lrg M$, since they have no
control on the induced deformations on $ X \times Y$.

Therefore, we need to define a new functor: the deformation
functor associated with a pair of morphisms.

\smallskip

Given  morphisms of differential graded Lie algebras $h :L \lrg M$
and $g: N \lrg M$:
\begin{center}
$\xymatrix{ & L \ar[d]^h   \\
N\ar[r]^{g}   & M, \\ }$
\end{center}
we define the functor
$$
\Def_{(h,g)}: \Art \lrg \Set,
$$
$$
\Def_{(h,g)}(A) =\{(x,y,e^p) \in (L^1 \otimes m_ A)\times  (N^1
\otimes m_A ) \times \operatorname{exp}(M^0 \otimes m_A)  |
$$
$$
dx + \displaystyle\frac{1}{2} [x,x]=0,\   dy +\frac{1}{2}[y,y]=0,\
g(y)=e^p*h(x) \}/gauge,
$$
where this gauge equivalence is an extension of the previous ones
(see Definition~\ref{defin funtor Def_(h,g) non esteso}).

This functor is a generalization of the previous ones: by choosing
$N=0$ and $g=0$,  $ \Def_{(h,g)}$ reduces to $ \Def_h$; by
choosing $N=M=0$ and $h=g=0$, $ \Def_{(h,g)}$ reduces to $
\Def_L$.

\medskip

Consider the DGLA  $A^{0,*}_{X }(\Theta_X) \times
A^{0,*}_{Y}(\Theta_Y) $   and the morphism $g=(p^*,q^*):A^{0,*}_{X
}(\Theta_X) \times A^{0,*}_{Y}(\Theta_Y) \lrg A_{X\times
Y}^{0,*}(\Theta_{X\times Y})$, where $p$ and $q$ are the natural
projections of the product $X \times Y$ on $X$ and $Y$,
respectively: $p:X \times Y\lrg X$ and $q:X \times Y\lrg Y$.

We note that the solutions  $n=(n_1,n_2)$ of the Maurer-Cartan
equation in $N=A^{0,*}_{X  }(\Theta_X) \times
A^{0,*}_{Y}(\Theta_Y) $ correspond to  infinitesimal deformations
of both $X$ (induced by $n_1$) and  $Y$ (induced by $n_2$).
Moreover the image $g(n)$ satisfies the Maurer-Cartan equation in
$M=A^{0,*}_{X \times Y}(\Theta_{X \times Y})$ and so it is
associated with an infinitesimal deformation of $ X \times Y$,
that is exactly the one obtained as product of the deformations of
$X$ (induced by $n_1$) and of $Y$ (induced by $n_2$).

\smallskip

Therefore, this $g$ gives  exactly the   control  on the
deformations of the product that we are looking for.

Let $M=A_{X\times Y}^{0,*}(\Theta_{X\times Y}) $, $L=A_{X \times
Y}^{0,*}(\Theta_{X \times Y}(-log\, \Gamma)) $, $h:L \lrg M$ the
inclusion, $N=A^{0,*}_{X }(\Theta_X) \times A^{0,*}_{Y}(\Theta_Y)
$ and $g=(p^*,q^*):N \lrg M$. Then we are in the following
situation:
\begin{center}
$\xymatrix{ & & A_{X \times Y}^{0,*}(\Theta_{X \times Y}(-log\,
\Gamma))
\ar@{^{(}->}[d]^{\ h}   & \\
A_X^{0,*}(\Theta_X) \times A_Y^{0,*}(\Theta_Y)
\ar[rr]^{g=(p^*,q^*)}
 & & A_{X\times Y}^{0,*}(\Theta_{X\times
Y}).  &
 \\ }$
\end{center}
In conclusion, each deformation of the map $f$ corresponds to a
Maurer-Cartan solution  $l \in A_{X \times Y}^{0,*}(\Theta_{X
\times Y}(-log\, \Gamma))$, such that $h(l)$ induces a deformation
of  $ X \times Y$ isomorphic to a deformation induced by $g(n)$,
for some Maurer-Cartan solution $n \in A^{0,*}_{X  }(\Theta_X)
\times A^{0,*}_{Y}(\Theta_Y)$ (that is, $h(l)$ and $g(n)$ are
gauge equivalent in $A_{X\times Y}^{0,*}(\Theta_{X\times Y})  $).

Therefore, $\Def_{(h,g)}$ encodes all the geometric data of the
problem and the following theorem  is quite obvious (see
Theorem~\ref{teo Def_(h,g)  Def (f)}).

\begin{teos}[B]
Let $f:X \lrg Y$ be a holomorphic map of compact complex manifold.
Then, with the notation above, there exists an isomorphism of
functors
$$
\Def_{(h,g)}\cong \Def (f).
$$
\end{teos}

This theorem holds for the general case of infinitesimal
deformations of $f$, the other cases are obtained as
specializations of it.

For example, the deformations of $f$ with fixed domain or fixed
target, are obtained  by considering  $N=A^{0,*}_{Y }(\Theta_Y)$
or $N=A^{0,*}_{X  }(\Theta_X)$, respectively.

In particular, using path objects, for each choice of $h :L \lrg
M$ and $g:N\lrg M$, we are able to find a differential graded Lie
algebra $H_{(h,g)}$ such that $\Def_{\hil}\cong \Def_{(h,g)}$ (see
Theorem~\ref{teo EXT DEF_H=DEF(h,g)}).

Therefore, we give  an explicit description (more than the
existence) of a DGLA that controls the deformations of holomorphic
maps (Theorem~\ref{teo esiste hil governa DEF(f)}).

\bigskip

Finally, we apply these techniques to  study the obstructions to
deform holomorphic maps.

The idea is the following: if we have an infinitesimal deformation
of a geometric object, then we want to know if it is possible to
extend it.

More precisely, let $F:\Art \lrg \Set$ be a deformation functor. A
(complete) obstruction space for $F$ is a vector space $V$, such
that  for each surjection $B \lrg A$ in $\Art$ and each element $x
\in F(A)$, there exists an obstruction  element $v_x \in V$,
associated with $x$, that is zero if and only if $x$ can be lifted
to $F(B)$ (for   details see Section~\ref{sezione ostruzione}).

Therefore, we would like to control this obstruction space and
know when the associated obstruction element is zero.

In general, we just know a vector space that contains these
elements but we have no explicit  description of which elements
actually are  obstructions. Among other things, if $W$ is another
vector space which contains $V$, then also $W$ is an obstruction
space for $F$. Then, in some sense we are looking for the \lq \lq
smallest" obstruction space (see Remark~\ref{oss find smallest
obstruct theory}).

For example, the obstructions of the functor associated with a
DGLA $L$ are naturally contained in $H^2(L)$ (see
Section~\ref{section defi funtor DEF_L of DGLA}), but we do not
know which classes in $H^2(L)$ are indeed obstructions.

\medskip

In the case of a complex compact manifold $X$, an obstruction
space for the deformation functor $\Def_X$ is the second
cohomology vector space $H^2(X, \Theta_X)$ of the holomorphic
tangent bundle $\Theta_X$ of $X$  (Theorem~\ref{teo sernesi ostruz
DEF_X=H^2(X,T_X)}).

If $X$ is also K\"{a}hler, then A.~Beauville, H.~Clemens \cite{bib
clemens} and Z.~Ran \cite{bib Ran ostruzio} \cite{bib Ran
ostruzioII} proved that the obstructions are contained in a
subspace of $H^2(X, \Theta_X)$ defined as the kernel of a well
defined map. This is the so-called \lq \lq Kodaira's principle"
(see for example \cite[Theorem~10.1]{bib clemens},
\cite[Corollary~3.4]{bib mane COSTRAINT} or
\cite[Corollary~12.6]{bib fiorenz-Manet}, \cite[Theorem~0]{bib Ran
ostruzio} or \cite[Corollary~3.5]{bib Ran ostruzioII}).

\medskip

In the case of embedded deformations of a submanifold $Z$ in a
fixed manifold $X$, then the obstructions are naturally contained
in the first cohomology  $H^1(Z, N_{Z|X})$ of the normal bundle
$N_{Z|X}$ of $Z$ in $X$. In this case too, if $X$ is K\"{a}hler,
it is possible to define a map on $H^1(Z, N_{Z|X})$, called the
\lq \lq semiregularity map", that contains the obstructions in the
kernel. The idea of this map is due to S.~Bloch \cite{bib bloch}
and it is also studied, using the DGLA approach, by M.~Manetti
\cite[Theorem~0.1 and Section~9]{bib ManettiPREPRINT}.

\bigskip



In the case of deformations of a   holomorphic map $f:X \lrg Y$
with fixed codomain,  it was proved by E.~Horikawa  in \cite{bib
HorikawaI} (see Theorem~\ref{teo TANG e OSTR DEF(f) horikawa})
that the obstructions are contained in the second hypercohomology
group $ \mathbb{H}^2\left(X,\shO(\Theta_X) \stackrel{f_*}
{\lrg}\shO(f^*\Theta_Y)\right) $.

\medskip

Using the approach via DGLA that we have explained, we can give an
easy proof of this theorem (Proposition~\ref{prop mia uguale
horikawa}) but, maybe most important,  we can improve it in the
case of K\"{a}hler manifolds (Corollary~\ref{coroll
semiregolari-ostruzion}).

\medskip

Actually, let $n=dim X$, $p=dim Y - dim X$  and $\mathcal{H}$ be
the space of harmonic forms on $Y$ of type $(n+1,n-1)$. By
Dolbeault's theorem and Serre's duality we obtain the equalities
$\mathcal{H}^\nu=(H^{n-1} (Y,\Omega_Y^{n+1}))^\nu=
H^{p+1}(Y,\Omega_Y^{p-1})$.

Using the contraction $\contr$ of vector fields with differential
forms (see Sections~\ref{sec cotrction map e Lie deriv} and
\ref{sottosection semiregolarity}), for each $\omega \in
\mathcal{H}$ we can define the following map
$$
A_X^{0,*}(f^*\Theta_Y) \stackrel{\contr \omega}{\lrg}
A_X^{n,*+n-1}
$$
$$
\contr \omega( \phi f^* \chi)=\phi f^*(\chi \contr \omega)\in
A_X^{n,p+n-1}  \qquad \forall\  \phi f^* \chi \in
A_X^{0,p}(f^*\Theta_Y).
$$
By choosing $\omega $ such that  $f^* \omega=0$, we get the
following commutative diagram (see Corollary~\ref{coroll
semiregolari-ostruzion})
\begin{center}
$\xymatrix{A_X^{0,*}(f^*\Theta_Y) \ar[r]^{\ \ \contr \omega}
& A_X^{n,*+n-1} \\
A_X^{0,*}(\Theta_X) \ar[u]^{f_*} \ar[r]  & 0. \ar[u]     \\
}$
\end{center}
Then, for each $\omega$ we get  a morphism
$$
\mathbb{H}^2\left(X,\shO(\Theta_X)
\stackrel{f_*}{\lrg}\shO(f^*\Theta_Y)\right) \lrg
 H^{n}(X,\Omega^n_ X),
$$
which composed with the integration on $X$ gives
$$
\sigma : \mathbb{H}^2\left(X,\shO(\Theta_X)
\stackrel{f_*}{\lrg}\shO(f^*\Theta_Y)\right) \lrg
H^{p+1}(Y,\Omega_Y^{p-1}).
$$
Using $\sigma$ we get the following theorem (see
Corollary~\ref{coroll semiregolari-ostruzion}).

\begin{teos}[C]
Let $f:X \lrg Y$ be a holomorphic map of  compact K\"{a}hler
manifolds. Let $p=dim Y - dim X$. Then the obstruction space to
the infinitesimal deformations of $f$  with fixed $Y$  is
contained in the kernel of the map
$$
\sigma : \mathbb{H}^2\left(X,\shO(\Theta_X) \stackrel{f_*}
{\lrg}\shO(f^*\Theta_Y)\right) \lrg H^{p+1}(Y,\Omega_Y^{p-1}).
$$
\end{teos}

\bigskip

The structure of this work is as follows.

\medskip

\emph{Chapter I} contains the basic material about deformation
functors. In Section~\ref{sezione funotri di Artin} we define the
deformation functors of  Artin rings, tangent and obstruction
spaces and some related properties.

Section~\ref{sezio funtore def infinitesimali DEF  X} is devoted
to studying the deformation functor $\Def_X$ of the infinitesimal
deformations of a compact complex manifold $X$.

In the last Section~\ref{sezio DGLA e funtore DEF_L} we introduce
the  differential graded Lie algebras (DGLAs)   and two associated
functors: the Maurer-Cartan functor $\MC_L$ and the deformation
functor $\Def_L$ (for each DGLA $L$).

\medskip

In  \emph{Chapter II} we fix the notation  about complex
manifolds. We recall the notions of differential forms
(Section~\ref{sezio differential form}) of \v{C}ech and Dolbeault
cohomology (Section~\ref{sezio diff dolbeault su fibrati holo})
and some properties of K\"{a}hler manifolds (Section~\ref{sezio
kaler manifold}). We also study the maps $f_*$ and $f^*$ induced
by a  holomorphic map $f$ (Section~\ref{sezione f_* and f^*}). In
particular, in Section~\ref{sezio DGLA di kodaira-spencer}, we
define the Kodaira-Spencer differential graded Lie algebra
$A_X^{0,*} (\Theta_X)$ associated with a complex manifold $X$.

Section~\ref{sezio defo manifold defK=defX} contains the proof of
theorem A (Theorem~\ref{teo def_k =Def_X}): we prove the existence
of an isomorphism $\Def_{A_X^{0,*}(\Theta_X)}\cong\Def_X$ and so
the Kodaira-Spencer  algebra $A_X^{0,*} (\Theta_X)$ controls the
infinitesimal deformation of $X$.

\medskip

\emph{Chapter III} is the technical bulk of this thesis. In
Section~\ref{sezi def MC_(h,g) and DEF(h,g) no ext} we define the
Maurer-Cartan functor $\MC_{(h,g)}$  and  the deformation functor
$\Def_{(h,g)}$  associated with a pair of morphisms of DGLAs $h:L
\lrg M$ and $g:N \lrg M$. Sections~\ref{sez tangent e ostru
MC(h,g) DEF(h,g)} and \ref{sez propriet di MC(h,g) DEF(hg)} are
devoted to the study of some properties of these functors as for
example tangent and obstruction spaces.

In Section~\ref{sezio extended defo funtore} we introduce the
extended deformation functors to prove the existence of a DGLA
$H_{(h,g)}$ such that  $\Def_{(h,g)} \cong \Def_{H_{(h,g)}} $
(Theorem~\ref{teo EXT DEF_H=DEF(h,g)}).

\medskip

In \emph{Chapter IV} we study the infinitesimal deformations of
holomorphic maps.

Section~\ref{sezio definizio funtore DEF(f)} is devoted to
defining the deformation functor $\Def(f)$ of infinitesimal
deformations of a holomorphic map $f:X \lrg Y$ of compact complex
manifolds.

In Section~\ref{sez teo Def_(h,g)Def (f)}, we prove theorem B,
i.e., the existence of a pair of morphisms of DGLAs $ h: A_{X
\times Y}^{0,*}(\Theta_{X \times Y}(-log\, \Gamma))\hookrightarrow
A_{X\times Y}^{0,*}(\Theta_{X\times Y})$ and $g=(p^*,q^*):
A_X^{0,*}(\Theta_X) \times A_Y^{0,*}(\Theta_Y) \lrg A_{X\times
Y}^{0,*}(\Theta_{X\times Y}) $, such that $\Def_{(h,g)} \cong
\Def(f)$ (see Theorem~\ref{teo Def_(h,g)  Def (f)}).

\medskip

\emph{Chapter V} contains examples and applications of the
techniques described before. In Section~\ref{sezione fixed target}
we study the infinitesimal deformations of holomorphic maps with
fixed codomain and Section~\ref{sottosection semiregolarity}
contains the main result about the semiregularity map
(Corollary~\ref{coroll semiregolari-ostruzion}).

Then we study infinitesimal deformations of a holomorphic map with
fixed  domain and codomain (Section~\ref{sezio fixed target
domain}) and the infinitesimal deformations of an inclusion
(Section~\ref{sezio inclusion}).

\bigskip

{ \bf Acknowledgments.} It is a pleasure for me to show my deep
gratitude to the advisor of this thesis Prof. Marco Manetti, an
excellent   helpful professor who supports and encourages me every
time. I'm indebted with him for many useful discussions, advices
and suggestions. Several ideas of this work are grown under his
influence. I wish also to thank him for having introduced me to
the exciting subject of deformation  theory.

\smallskip

I can not forget here my graduate study advisor Prof. Marialuisa
J. de Resmini, who encouraged me to continue my studies and \lq
\lq sent" me to a summer school where I had my first meeting with
Algebraic Geometry. I am very grateful  for her improvements to
the final version of this thesis.

\smallskip

I am also in debt with all those people  with whom I have useful
mathematical  discussions or who make my staying in Rome pleasant.

In this set there  certainly are Paolo Antonini, Sara Azzali,
Francesco Esposito, Elena Martinengo, Silvia Montarani and Antonio
Rapagnetta.

I wish also to thank three not mathematical friends of mine
Arianna, Federica and Patrizia who support me in many and
different ways.

\smallskip

This thesis is dedicated to my family Maria Teresa, Giovanni
Michele, Enrico and especially to (the memory of) my grandmother
Mina.

\chapter{Functors of Artin rings}

In this chapter we collect  some  definitions and main properties
of deformation functors.

In the first section, we introduce the notions of functor of Artin
rings, of deformation functors and we define the tangent and
obstruction spaces.

Section~\ref{sezio funtore def infinitesimali DEF  X} is devoted
to the study of the deformation functor $\Def_X$ of infinitesimal
deformation of complex manifolds.

In Section~\ref{sezio DGLA e funtore DEF_L} we  introduce the
fundamental notions of differential graded Lie algebra (DGLA) $L$
and of deformation functor associated with a DGLA $\Def_L$.

The main references for this chapter are \cite{bib Fantec-Manet},
\cite{bib manetPISA}, \cite{bib manRENDICONTi}, \cite{bib Sernesi}
and \cite{bib Schlessinger }.

\section{Generalities on functors of Artin rings}
\label{sezione funotri di Artin}

Let $\K$ be a fixed field of characteristic  zero.

Consider the following categories:
\begin{itemize}
\item  $\Set$: the category of sets in a fixed universe
with $\{*\}$ a fixed set of cardinality 1;

\item   $\Art=\Art_\K$ : the category of local Artinian
$\K$-algebras with residue field $\K$ ($A/m_A=\K$);

\item  $\widehat{\Art}=\widehat{\Art_\K}$ : the category of
complete local Noetherian  $\K$-algebras with residue field $\K$
($A/m_A=\K$).

For each $S  \in \widehat{\Art}$ we also consider:

\item  $\Art_S$ : the category of local Artinian $S$-algebras
with residue field $\K$ (for   such an element $A$, the structure
morphism $S \lrg A$ induces a trivial extension of the residue
field $\K$);

\item  $\widehat{\Art_S}$ : the category of complete local
Noetherian  $S$-algebras with residue field $\K$.

\end{itemize}

\begin{oss}
We note that  $\Art_S \subset \widehat{\Art_S}$. Moreover, by
morphisms in a category of local objects we mean local morphisms
and we often use the notation $A \in \mathbf{C}$ instead of $A \in
ob(\mathbf{C})$, when $\mathbf{C}$ is a category.

\end{oss}
If $\beta: B \lrg A$ and $\gamma:C \lrg A$ are morphisms in
$\widehat{\Art_S}$ (in $\Art_S$, respectively), then
$$
B\times_A C=\{(b,c)\in B \times C \, |\,  \beta(b)=\gamma(c) \}
$$
is the \emph{fiber product} of $\beta$ and $\gamma$ and $B\times_A
C \in \widehat{\Art_S}$ (in $\Art_S$, respectively).

\begin{dhef}
A \emph{small extension } in $\widehat{\Art_S}$ (in $\Art_S$,
respectively) is a short exact sequence
$$
e: \ \ 0 \lrg J \lrg B \stackrel{\alpha}{\lrg} A \lrg 0,
$$
where $\alpha$ is a morphism in $\widehat{\Art_S}$ (in $\Art_S$,
respectively) and the kernel $J$ is an ideal of $B$ annihilated by
the maximal ideal $m_B$, $m_B \cdot J=(0)$. This implies that $J$
is $\K$-vector space.

A small extension is called \emph{principal} if  $J$ is a one
dimensional vector space ($J\cong \K$).
\end{dhef}

We will often say that a morphism $\alpha :B\lrg  A$ is a small
extension, meaning that $0 \lrg \ker(\alpha) \lrg B
\stackrel{\alpha}{\lrg} A \lrg 0$ is a small extension.

\begin{oss}\label{oss surje sono small extes}

Every surjective morphism  in $\Art_S$ can be expressed as a
finite composition of small extensions.

Actually, let $B \lrg A$ be a surjection with kernel $J$:
$$
0 \lrg J \lrg B \lrg A \lrg 0.
$$
Since $B$ is a local artinian ring, its maximal ideal $m$ is
nilpotent: there exists $n_0 \in \mathbb{N}$ such that $m^n=0$ for
each $n \geq n_0$; in particular, $m^n J=0$.

Therefore, it is sufficient to consider the sequence of small
extensions
$$
0 \lrg m^nJ/m^{n+1}J \lrg B/m^{n+1}J \lrg  B/m^{n}J\lrg 0.
$$

\end{oss}

\begin{oss}
In view of  Remark~\ref{oss surje sono small extes}, it will be
enough to check  the surjection for small extensions instead of
verifying the surjection for any morphism  in $\Art_S$.

\end{oss}

\bigskip

\begin{dhef}\label{def funtori Artin}
A \emph{functor of Artin rings} is a covariant functor $F:
 \Art_S  \lrg \Set$, such that $F(\K)=\{*\}$ is the one
point set.
\end{dhef}

The functors of Artin rings $F:  \Art_S  \lrg \Set$ and their
natural transformations build a category denoted by $\Fun_S$. A
natural transformation of functors $\gamma:F \lrg G$ is an
\emph{isomorphism of functors } if and only if $\gamma(A):F(A)
\lrg G(A)$ is bijective for each $A \in \Art_S $.

\begin{exe}
The trivial functor $F(A)=\{*\}$, for every $A \in \Art_S$.
\end{exe}

\begin{exe}
Let $R \in \widehat{\Art_S}$. We define
$$
h_R: \Art_S  \lrg \Set,
$$
such that
$$
h_R(A)=\Hom_S(R,A).
$$
$F \in \Fun_S$ is called \emph{pro-representable} if it is
isomorphic to $h_R$, for some $R \in \widehat{\Art_S}$.

If we can choose $R \in \Art_S$ then $F$ is called
\emph{representable}.

\end{exe}

Let $\kepsi$, with $\epsi^2=0$, be the ring of dual numbers over
$\K$. $\kepsi=\K \oplus \K\epsi$ is a $\K$-vector space of
dimension 2 and  has a trivial $S$-algebra  structure (induced by
$S \lrg \K \lrg \kepsi$).

\begin{dhef}
The set $t_F:=F(\K[\epsi])$ is called the \emph{tangent space} of
$F \in \Fun_S$.
\end{dhef}

Let
\begin{equation}\label{equa eta:F(AX_C B) ->F(A)XF(C) F(B)}
\eta: F(B\times_A C)\lrg F(B)\times_{F(A)}F(C)
\end{equation}
be the map induced by the fiber product in $\Art_S$:
\begin{center}
$\xymatrix{B\times_A C \ar[r]\ar[d] & C \ar[d] \\
           B \ar[r] & A.   \\ }$
\end{center}

\begin{dhef}
A functor $F$ is called \emph{homogeneous} if $\eta$ is an
isomorphism whenever $B \lrg A$ is surjective.
\end{dhef}

\begin{dhef}\label{defin deformation funtore}
A functor $F$ is called a  \emph{deformation functor } if
\begin{itemize}
  \item [i)] $\eta$ is surjective whenever $B \lrg A$ is surjective;
  \item [ii)] $\eta$ is an isomorphism whenever $A=\K$.
\end{itemize}

\end{dhef}


\begin{oss}
The deformation functors will play an important role in this work.

In particular, we will study the following four deformation
functors:
\begin{itemize}
  \item[{\bf 1)}] the functor $\Def_X$ of infinitesimal deformation of
  complex manifolds, in section
  \ref{sezio funtore def infinitesimali DEF  X};
  \item[{\bf 2)}] the functor $\Def_L$ associated with a differential
  graded Lie algebra $L$,
  in Section~\ref{sezio DGLA e funtore DEF_L};
  \item[{\bf 3)}] the functor $\Def_{(h,g)}$ associated with a pair
  of morphisms of differential graded Lie algebras
  $h:L \lrg M$ and $g:N \lrg M$, in section
  \ref{sezio non ext def(h,g) e MC(h,g)};
  \item[{\bf 4)}]  the functor $\Def_f$ associated with the infinitesimal
  deformations of a holomorphic map $f$, in section
  \ref{sezio definizio funtore DEF(f)}.
\end{itemize}
\end{oss}

\begin{exe}\label{exe DEF_Xsotto definizione funtore deformaz}
Let $X$ be an algebraic scheme over $\K$ (separated of finite type
over $\K$). Define the following functor
$$
\Def_X: \Art \lrg \Set,
$$
where $\Def_X(A)$ is the set of isomorphism classes of commutative
diagram:
\begin{center}
$\xymatrix{X \ar[r]^i\ar[d] & X_A\ar[d]^{p_A} \\
           \Spec(\K) \ar[r] & \Spec(A).   \\ }$
\end{center}
where $i$ is a closed embedding and $p_A$ is a flat morphism.
$\Def_X$ is a deformation functor (see \cite[Section~3.7]{bib
Schlessinger } or \cite[Prop.~III.3.1]{bib Sernesi}).

\end{exe}

\begin{prop}
Let $F$ be a deformation functor. Then $t_f$ has a natural
structure of $\K$-vector space and  every natural transformation
of deformation functors $\phi:F \lrg G$ induces a linear map
between tangent spaces.
\end{prop}
\begin{proof}
Since $F(\K ) $ is just one point and  the morphism  $\eta$
defined in (\ref{equa eta:F(AX_C B) ->F(A)XF(C) F(B)}) is an
isomorphism for $A=\K$, we have $F(\kepsi)\times F(\kepsi)\cong
F(\kepsi \times_{\K} \kepsi)$.

Consider the map
$$
+: \kepsi \times_{\K} \kepsi \lrg \kepsi,
$$
$$
\qquad (a+b\epsi,a+b'\epsi)\longmapsto a+ (b+b')\epsi.
$$
Then using the previous  isomorphism, the map $+$ induces the
addition on $F(\kepsi)$:
$$
F(\kepsi)\times F(\kepsi) \stackrel{\cong}{ \lrg} F(\kepsi
\times_{\K} \kepsi)  \stackrel{F(+)}{ \lrg} F(\kepsi).
$$
Analogously, for the multiplication by a scalar $k \in \K$ we
consider the map:
$$
k: \kepsi   \lrg \kepsi,
$$
$$
\qquad a+b\epsi \longmapsto a+ (kb)\epsi.
$$

\end{proof}

\begin{oss}
In the previous proposition we just used the fact that $F(\K)$ is
one point and that $\eta$ in (1) is an isomorphism for $A=\K$ and
$C=\kepsi$

\end{oss}

\begin{dhef}
A morphism $\phi:F \lrg G$   in $\Fun_S$ is:
\begin{itemize}
  \item[-] $unramified$ if $\phi: t_F \lrg t_G$ is injective;
  \item[-] $smooth$ if  the map
$$
F(B) \lrg G(B) \times_{G(A)}F(A)
$$
induced by the diagram
\begin{center}
$\xymatrix{F(B) \ar[r]\ar[d]_\phi &  F(A)  \ar[d]^\phi \\
          G(B) \ar[r] & G(A),   \\ }$
\end{center}
is surjective for every surjection $B \lrg A$ in $\Art_S$;
  \item[-] $\acute{e}tale$ if $\phi$ is both smooth and unramified.
\end{itemize}
\end{dhef}

\begin{oss}\label{oss smooth allora F(B)->>>G(B)}
If $\phi:F \lrg G$ is smooth then, by taking $A=\K$ in the
previous definition, we conclude that  $\phi:F(B) \lrg  G(B)$ is
surjective.

\end{oss}

\begin{oss}\label{oss etale allora iso su tangenti}
If $\phi:F \lrg G$ is an \'{e}tale  morphism, then it induces an
isomorphism   $\phi':t_F \lrg t_G$ on tangent spaces. Actually,
$\phi$ is unramified and so $\phi'$ is injective. By hypothesis,
$\phi$ is also smooth and so, applying  Remark~\ref{oss smooth
allora F(B)->>>G(B)} in the special case $B=\kepsi$, $\phi'$ is
surjective.

\end{oss}

\begin{dhef}
A \emph{functor $F$ is smooth} if  for every surjection $B \lrg A$
of $S$-algebras  $F(B) \lrg F(A)$ is surjective (that is, the
morphism $F \lrg *$ is smooth).

\end{dhef}

\begin{prop}\label{prop etal+homog allora iso}
Let $\phi:F \lrg G$ be an $\acute{e}tale$ morphism of deformation
functors. If $G$ is homogeneous then $\phi$ is an isomorphism.
\end{prop}
\begin{proof}
Since $\phi$ is smooth, $\phi$ is surjective. Therefore, it is
sufficient to prove the injectivity of $\phi$ using  that
$G$ is homogeneous and $\phi$ is unramified.\\
The proof of this fact is taken from \cite[Lemma~2.10]{bib
manetPISA}; it is given  for the sake of completeness.\\
We prove it by induction on the length of $A$.

If $A=\K$ then $F(\K)=G(\K)=\{*\}$ and so the statement is
obvious.

Let
$$
0 \lrg \K\epsi \lrg B \lrg A \lrg 0,
$$
be a principal small extension ($\epsi\cdot m_B=0$). By induction
$\phi: F(A) \lrg G(A)$ is injective.

Consider the following isomorphism of $S$-algebras:
$$
\varphi: B \times_\K \kepsi \lrg B \times_A B,
$$
$$
(b,\overline{b}+\beta \epsi) \longmapsto (b,b+\beta \epsi).
$$
Since $F$ is   a deformation functor,  $F(\varphi): F(B) \times
t_F \lrg F(B) \times_{F(A)} F(B)$ is surjective.

Since $G$ is also  homogeneous,   $G(\varphi): G(B) \times t_G
\lrg G(B) \times_{G(A)} G(B)$ is an isomorphism. We note that
$G(\varphi)(G(B) \times \{0\})=\Delta$ diagonal.

Now, suppose that $\phi(\xi)=\phi(\eta) \in G(B)$ for $\xi$ and
$\eta \in F(B)$.

Since $\phi$ is injective on $F(A)$,  $(\xi,\eta) \in
F(B)\times_{F(A)} F(B)$.

Moreover, the surjectivity of $F(\varphi)$ implies the existence
of an element $h \in t_F$ such that
$F(\varphi)(\xi,h)=(\xi,\eta)$. Then
$G(\varphi)(\phi(\xi),\phi(h))=(\phi(\xi),\phi(\eta))\in \Delta$
and so $\phi(h)=0$ ($G$ is  an isomorphism).

By hypothesis, $\phi$ is unramified, therefore, $h=0$ and so
$\xi=\eta$.
\end{proof}

\begin{oss}\label{oss su azione tangente funtore}
Let $F$ be a deformation functor and $f:B \lrg A$ a surjection
with $\ker f\cong\K$. Using the isomorphism $B \times_\K \kepsi
\cong B \times_A B$ of  Proposition~\ref{prop etal+homog allora
iso}, we obtain a commutative diagram
\begin{center}
$\xymatrix{F(B) \times t_F \ar[r]^\pi\ar[d]_\tau  &  F(B)  \ar[d] \\
          F(B) \ar[r]_{F(f)} & F(A).   \\ }$
\end{center}
$\pi$ is the projection and $\tau$ defines an action of $t_F$ on
$F(B)$ which restricts to a transitive action on each fiber of
$F(f)$ (see \cite[Lemma~2.12]{bib Fantec-Manet}).
\end{oss}

\begin{cor}\label{cor funtoreF banale sse t_F=0}
Let $F$ be a deformation functor, then $F=\{*\}$ if and only if
$t_F=(0)$.
\end{cor}

\begin{proof}
One implication is obvious. Next, let $F$ be a deformation functor
such that $t_F=(0)$. We prove that $F(A)=\{*\}$ by induction on
$dim_\K(A)$. If $A=\K$ then there is nothing to prove (by
definition of functor of Artin rings). Next, let $\pi:B \lrg A$ be
a small extension  and suppose that $F(A)=\{*\}$. By
Remark~\ref{oss su azione tangente funtore}, $t_F=(0)$ acts
transitively on the unique fiber $F(B)$ of the map $F(\pi)$. This
implies $F(B)=\{*\}$.
\end{proof}

\subsection{Obstruction theory}\label{sezione ostruzione}

\begin{dhef}\label{def ostruzione space}
Let F be a functor of Artin rings; an \emph{obstruction theory }
for $F$ is a pair $(V,v_e)$ such that:
\begin{itemize}
  \item  $V$ is a $\K$-vector space, called \emph{obstruction
  space};
  \item    for every small
extension  in $\Art_S$
$$
e: 0 \lrg J \lrg B \stackrel{\alpha}{\lrg} A \lrg 0
$$
$v_e: F(A) \lrg V \otimes_\K J$ is an obstruction map, satisfying
the following properties:
\begin{itemize}
  \item If $\xi \in F(A)$ can be lifted to $F(B)$ then
  $v_e(\xi)=0$.
  \item For any morphism $\phi: e_1 \lrg e_2$ of small
  extension, i.e.,
\begin{equation}\label{equa morfism small exte a:e_1 ->e_2}
 \xymatrix{e_1:& &0 \ar[r]& J_1\ar[r]\ar[d]^{\phi_J} &B_1 \ar[r]
\ar[d]^{\phi_B} & A_1 \ar[r] \ar[d]^{\phi_A} &0 & & \\
           e_2:& &0 \ar[r]& J_2\ar[r] &B_2 \ar[r] & A_2 \ar[r]
           &0& & \\ }
\end{equation}
we have $v_{e_2}(\phi_A(a))=(Id_V \otimes \phi_J)(v_{e_1}(a))$,
for every $a \in F(A_1)$.
\end{itemize}
\end{itemize}
\end{dhef}

\begin{dhef}
An obstruction theory for a functor is  \emph{complete} if the
lifting exists if and only if the obstruction  vanishes.
\end{dhef}

\begin{oss}
If $F$ has $(0)$ as complete obstruction space then $F$ is smooth.
In Proposition~\ref{prop f liscio se e solo se OF=0} we will prove
that the converse is also true for a deformation functor.
\end{oss}

\begin{oss}\label{oss v__G o f ostruzione per F}
Let $\psi:F \lrg G$ be a natural transformation of functors and
$(V,v_e)$ an obstruction theory for $G$. Then $(V,v_e' :=v_e \circ
\psi)$ is an obstruction theory for $F$.

Actually, consider the small extension  in $\Art_S$
$$
e: 0 \lrg J \lrg B \stackrel{\alpha}{\lrg} A \lrg 0
$$
and the map
$$
v_e':F(A)\stackrel{F(\psi)}{\lrg}G(A) \stackrel{v_e}{\lrg} V
\otimes _\K J.
$$
Let $\xi \in F(A)$ and suppose that it can be lifted to
$\tilde{\xi} \in F(B)$. Therefore, $\psi(\xi) \in G(A)$ can be
lifted to $\psi(\tilde{\xi}) \in G(B)$ and so
$v_e'(\xi)=v_e(\psi(\xi))=0$.

Now, let $\phi: e_1 \lrg e_2$ be a morphism of small extension as
in (\ref{equa morfism small exte a:e_1 ->e_2}). Then for each $a
\in F(A_1)$ we have
$$
(Id_V \otimes \phi_J)(v_{e_1}'(a))=(Id_V \otimes
\phi_J)(v_{e_1}(\psi(a)))=
$$
$$
v_{e_2}(\phi_A(\psi(a)))= v_{e_2}(\psi(\phi_A(a)))= v_{e_2}'
(\phi_A((a))).
$$

\medskip

Moreover, if the morphism $\psi$ is also smooth and $(V,v_e)$ is
complete for $G$, then $(V,v_e' )$ is complete for $F$.

Indeed, suppose that $\xi \in F(A)$ is such that $0=v_e' (\xi)=
v_e(\psi(\xi))$. Then there exists $\eta \in G(B)$ that lifts
$\psi(\xi) \in G(A)$ ($(V,v_e)$ is complete for $G$). Consider the
following diagram
\begin{center}
$\xymatrix{F(B) \ar[r]\ar[d] &  F(A)\ar[d]  \ni  \xi \ar[d] \\
          \eta \in G(B) \ar[r] & \ \ G(A)  \ni \psi(\xi).  \\ }$
\end{center}
Since $\psi$ is smooth, the map $ F(B)  \lrg
F(A)\times_{G(A)}G(B)$ is surjective and so there exists
$\tilde{\xi}\in F(B)$ that lifts $\xi$.
\end{oss}

\begin{oss}\label{oss find smallest obstruct theory}
If $V$ is complete  a obstruction theory for a functor $F$, using
embeddings of vector spaces we can construct infinitely many
complete obstruction theories. Therefore, the goal is to find a
\lq \lq smallest" complete obstruction theory. The main results in
this context is the following Theorem~\ref{teo exis ostruzio
fantech man}. First of all, we give a definition.
\end{oss}

\begin{dhef}
A morphism of obstruction theories $(V,v_e) \lrg (W,w_e)$ is a
linear map (of vector spaces) $\theta:V \lrg W$ such that
$w_e=\theta v_e$, for every small extensions $e$.

An obstruction theory $(O_F,ob_e)$ for $F$ is called
\emph{universal} if for any obstruction theory $(V,v_e)$ there
exists a unique morphism $(O_F,ob_e)\lrg (V,v_e)$.
\end{dhef}

\begin{teo}\label{teo exis ostruzio fantech man} (Fantechi, Manetti)
Let $F$ be  a deformation functor. Then there exists a universal
obstruction theory $(O_F,ob_e)$ for $F$. Moreover, the universal
obstruction theory is complete and every element of the vector
space $O_F$ is of the form $ob_e(\xi)$ for some principal
extension
$$
e: \ 0 \lrg \K\epsi \lrg B \lrg A \lrg 0.
$$
and some $\xi \in F(A)$.
\end{teo}
\begin{proof}
See \cite[Theorem~3.2 and Corollary~4.4]{bib Fantec-Manet}.
\end{proof}

Let $\phi:F \lrg G$ be a morphism of deformation functors and
$(V,v_e)$, $(W,w_e)$ obstructions theories for $F$ and $G$,
respectively. A linear map $\phi': V \lrg W$ is \emph{compatible}
with $\phi$ if $w_e \phi=\phi' v_e$ for every small extensions
$e$.

\begin{teo}\label{teo standard smoothn criterion}
Let $\phi:F \lrg G$ be a morphism of deformation functors and
$\phi':(V,v_e) \lrg (W,w_e)$ a compatible morphism of obstruction
theories. If $(V,v_e)$ is complete, $\phi'$ injective and $t_F
\lrg t_G$ surjective, then $\phi$ is smooth.

\end{teo}

\begin{proof}
See \cite[Proposition~2.17]{bib manetPISA}. We have to prove that
the map
$$
F(B) \lrg G(B) \times_{G(A)}F(A)
$$
is surjective, for all small extensions
$$
0 \lrg \K\epsi \lrg B \lrg A \lrg 0.
$$
Let $(b',a) \in G(B) \times_{G(A)}F(A)$ and $a'  \in G(A)$ their
common image, that is, $b' \in G(B)$ lifts $a' \in G(A)$ and so
$w_e(a')=0$.

By hypothesis,  $\phi'$ is injective and so $v_e(a)=0$
($0=w_e(a')=w_e(\phi(a))=\phi'(v_e(a))$). Therefore, there exists
$b \in F(B)$ that lifts $a$:
\begin{center}
$\xymatrix{b \in F(B) \ar[r]\ar[d] &  F(A)\ar[d]  \ni  a \ar[d] \\
          b' \in G(B) \ar[r] & G(A) \ni  a'.  \\ }$
\end{center}
In general $b$ does not lift $b'$. Let $b''= \phi(b) \in G(B)$;
then $(b'',b') \in G(B)\times_{G(A)} G(B)$.

As observed in the proof of Proposition~\ref{prop etal+homog
allora iso}, we have an isomorphism $B \times_\K \kepsi  \cong B
\times_A B$; since $G$ is a deformation functor,  there exists a
surjective morphism
$$
\alpha: G(B) \times t_G =G( B \times_A B) \lrg G(B)\times_{G(A)}
G(B).
$$
Therefore, there exists $h \in t_G$ such that $(b'',
\overline{b''}+h)$ is a lifting of $(b'',b''+h=b')$.

By hypothesis, $t_F \lrg t_G$ is surjective and so there exists a
lifting $k \in t_F$ of $h\in t_G$. Taking $k+b \in F(B)$ produces
a lifting of $a$ that maps on $b'$.

\end{proof}

\begin{oss}
Let $\phi: F \lrg G$ be a morphism of deformation functors and
$(O_F,v_e)$ and $(O_G,w_e)$ their universal obstruction theories.
Then $(O_G, w_e \circ \phi)$ is an obstruction theory for $F$.
Thus, by Theorem~\ref{teo exis ostruzio fantech man}, there exists
a morphism $o(\phi):(O_F,v_e)\lrg (O_G,w_e)$.

\smallskip

In conclusion, \emph{every morphism of deformation functors
induces a linear map both between tangent spaces and universal
obstruction spaces}.
\end{oss}

Next, we want to consider some useful properties between these
morphisms.

\begin{prop}\label{prop f liscio se e solo se OF=0}
Let $\phi: F \lrg G$ be a morphism of deformation  functors. Then
$\phi$ is smooth if and only if $t_F \lrg t_G$ is surjective and
$o(\phi): O_F \lrg O_G$ is injective. In particular,  $F$ is
smooth if and only if $O_F=0$.
\end{prop}
\begin{proof}
If  $t_F \lrg t_G$ is surjective and $o(\phi): O_F \lrg O_G$ is
injective then $\phi$ is smooth by Theorem~\ref{teo standard
smoothn criterion}.

Viceversa, suppose that $\phi$ is smooth, then by Remark~\ref{oss
smooth allora F(B)->>>G(B)}, $t_F=F(\kepsi) \lrg G(\kepsi)=t_G$ is
surjective. Let $B \lrg A$ be a small extension and $a$ an element
of $F(A)$, with obstruction  $ x \in O_F$   such that $o
(\phi)(x)=0 \in O_G$. By Theorem~\ref{teo exis ostruzio fantech
man}, $O_G$ is complete and so $\phi(a) \in G(A)$ can be lifted to
$b' \in G(B)$. Again, by Remark~\ref{oss smooth allora
F(B)->>>G(B)}, $F(B)\lrg G(B)$ is surjective and so there exists
$b \in F(B)$ that lifts $a$. Therefore, the obstruction of $a$ is
zero ($x=0$) and this proves that $o(\phi)$ is  injective.
\end{proof}

\begin{cor}\label{cor etale se bietti su tg e inj su ostru}
A morphism of deformation functors $\phi: F \lrg G$ is
$\acute{e}$tale if and only if $t_F \lrg t_G$ is bijective and
$o(\phi): O_F \lrg O_G$ is injective.
\end{cor}
\begin{proof}
If $\phi$ is \'{e}tale, then by Remark~\ref{oss etale allora iso
su tangenti} $t_F \lrg t_G$ is bijective. Since $\phi$ is also
smooth, $o(\phi): O_F \lrg O_G$ is injective.

Conversely, by  Proposition~\ref{prop f liscio se e solo se OF=0}
$\phi$ is smooth; by hypothesis,  $t_F \lrg t_G$ is injective, and
so $\phi$ is also unramified.
\end{proof}

\begin{cor}
Let $\phi:F \lrg G$ be a morphism of deformation functors with $G$
homogeneous. If $t_F \lrg t_G$ is bijective and $o(\phi): O_F \lrg
O_G$ is injective then $\phi$ is an isomorphism.
\end{cor}
\begin{proof}
Put together Corollary~\ref{cor etale se bietti su tg e inj su
ostru} and Proposition~\ref{prop etal+homog allora iso}.
\end{proof}

\begin{oss}
When $G$ is not homogeneous, $\phi$ could  be not  injective and
so we can just  conclude the surjectivity of $\phi$. Therefore, in
these cases we prove   the injectivity directly. This will happen
in Theorems~\ref{teo def_k =Def_X}  and~\ref{teo Def_(h,g) Def
(f)}.
\end{oss}

\begin{cor}\label{cor f:F_G liscio implica biezione ostruzio}
If $\phi:F \lrg G$ is smooth then $o(\phi): O_F \lrg O_G$ is
bijective.
\end{cor}
\begin{proof}
By Proposition~\ref{prop f liscio se e solo se OF=0}, we have just
to prove that  $o(\phi): O_F \lrg O_G$ is surjective. Let $B\lrg
A$ be a small extension  and $y \in O_G$ the obstruction to
lifting $a' \in G(A)$ to $b' \in G(B)$. Since $\phi$ is smooth,
there exists $a \in F(A)$, such that $\phi(a)=a'$, and $b \in
F(B)$, such that $\phi(b)=b'$. Therefore, the obstruction $x \in
O_F$  to lifting $a \in F(A)$ to $b \in F(B)$ is a lifting of $y
\in O_G$.

\end{proof}

\section{Deformation functor of complex manifolds}
\label{sezio funtore def infinitesimali DEF  X}

In this section we study the infinitesimal  deformation functor
associated with  a compact complex manifold. Thus we will work
over the complex numbers and so $\K=\C$ and $\Art=\Art_{\C}$.

The main references are \cite[Chapter~4]{bib kodaira libro},
\cite{bib Schlessinger }, \cite[Chapter~II]{bib Sernesi}.

\begin{dhef}\label{dhef infintesimal defor of X}
Let $X$ be a compact complex manifold and $A \in \Art$. An
\emph{infinitesimal deformation } of $X$ over $\Spec(A)$ is a
commutative diagram of complex spaces
\begin{center}
$\xymatrix{X \ar[r]^i\ar[d] & X_A\ar[d]^{\pi} \\
           \Spec(\C) \ar[r]^{a} & \Spec(A),   \\ }$
\end{center}
where   $\pi$ is a proper and flat holomorphic  map, $a \in
\Spec(A)$ is the closed point, $i$ is a closed embedding  and
$X\cong X_A \times_{\Spec(A)} \Spec(\C)$.

If $A=\K[\epsi]$ we call it a  \emph{first order} deformation of
$X$.
\end{dhef}

Sometimes, for an infinitesimal deformation $X_A$ over $\Spec(A)$,
we also use the short notation $(X_A,\pi,\Spec(A))$.

\begin{oss}
Let $X_A$ be an infinitesimal deformation of $X$. We note that, by
definition, it can be interpreted as a morphism of sheaves of
algebras $\shO_A\to \shO_X$ such that $\shO_A$ is flat over $A$
and $\shO_A\otimes_A \mathbb{C}\to \shO_X$ is an isomorphism.
\end{oss}

Given another deformation $X'_A$ of $X$ over $\Spec(A)$:
\begin{center}
$\xymatrix{X \ar[r] ^{i'}\ar[d] & X'_A\ar[d]^{\pi'} \\
           \Spec(\C) \ar[r]^{a} & \Spec(A),   \\ }$
\end{center}
we say that $X_A$ and $X'_A$ are \emph{isomorphic} if there exists
an isomorphism $\phi:X_A \lrg X'_A$ over $\Spec(A)$, that induces
the identity on $X$, that is, the following diagram is commutative

\begin{center}
$\xymatrix{ &X \ar[rd]^{i'} \ar[ld]_{i}& \\
X_A \ar[rr]^\phi\ar[dr]_{\pi} & & X'_A\ar[dl]^{\pi'} \\
           & \Spec(A). &   \\ }$
\end{center}

We note that for every $X$ we can always define the
\emph{infinitesimal product deformation}:
\begin{center}
$\xymatrix{X \ar[r]^{i\ }\ar[d] &  X\times \Spec(A)\ar[d]  \\
           \Spec(\C) \ar[r]^{a\ } & \Spec(A).   \\ }$
\end{center}

\begin{dhef}
An infinitesimal deformation of $X$ over $\Spec(A)$ is called
\emph{trivial} if it is isomorphic to the infinitesimal product
deformation.

$X$ is called \emph{rigid} if every infinitesimal deformation of
$X$ over $\Spec(A)$ (for each $A \in \Art$) is trivial.
\end{dhef}

For every deformation $X_A$ of $X$ over $\Spec(A)$ and every
morphism $A \lrg B$ in $\Art$ ($\Spec(B) \lrg \Spec(A)$), there
exists an associated deformation of $X$ over $\Spec(B)$, called
\emph{pull-back deformation}, induced by a basis  change:
\begin{center}
$\xymatrix{X \ar[r] \ar[d] & X_A\times_{ \Spec(A)}\Spec(B)\ar[d]  \\
           \Spec(\C) \ar[r]  & \Spec(B).   \\ }$
\end{center}

\begin{dhef}\label{defin funtore DEF_X}
The \emph{infinitesimal deformation functor} $\Def_X$ of a complex
manifold $X$ is defined as follows:
$$
\Def_X : \mathbf{Art}\to\mathbf{Set}, \qquad
$$
$$
 \ \ \ \ \  \ \ \ \ \ \qquad  A\longmapsto \Def_X(A)=
\left\{\begin{array}{c}  \mbox{  isomorphism   classes   of}\\
 \mbox{ infinitesimal   deformations }\\
 \mbox{ of  $X$   over   $\Spec(A)$}
\end{array} \right\}.
$$
\end{dhef}

\bigskip

\begin{prop}\label{prop infinites DEF_X is funtore deformazione}
$\Def_X$ is a deformation functor, i.e., it satisfies the
conditions of Definition~\ref{defin deformation funtore}.
\end{prop}

\begin{proof}
See \cite[Section~3.7]{bib Schlessinger } or
\cite[Proposition~III.3.1]{bib Sernesi}.
\end{proof}

\subsection{Tangent and obstruction spaces of $\Def_X$}
\label{subsction deformation of non singular varieties}

Let $X$ be a compact complex manifold and $\Theta_X$ its
holomorphic tangent bundle.

In this section we prove that the tangent space of $\Def_X$ is
$\check{H}^1(X,\Theta_X)$   and that the obstruction space is
naturally contained in $\check{H}^2(X,\Theta_X)$. First of all, we
recall a useful lemma.

\begin{lem}\label{lemm serneII.1.5 automor-derivazio}
Let $B_0$ be a $\C$-algebra and
$$
0 \lrg J \lrg B \lrg A \lrg 0
$$
a small extension in $\Art$. Then there is a 1-1 correspondence
$$
\left\{\begin{array}{c} \mbox{automorphisms of the trivial
deformation $B_0\otimes_\C B$} \\
\mbox{ inducing the identity on $B_0 \otimes_\C A$}
\end{array} \right\}
 \longleftrightarrow \operatorname{Der_\C(B_0,B_0)}\otimes J
$$
where the identity corresponds to the zero derivation, and the
composition of automorphisms corresponds to the sum of
derivations.
\end{lem}
\begin{proof}
See \cite[Lemma~II.1.5]{bib Sernesi}
\end{proof}

\begin{oss}\label{oss succ esatta lunga automorfismi derivazione}
Let $U$ be a Stein  open subset of a complex manifold $X$. Then
the previous lemma is equivalent to saying  that the following
sequence is exact:
$$
0 \lrg \Gamma(U_i, \Theta_X) \otimes J \lrg \Aut (\shO_X(U_i)
\otimes B) \lrg \Aut(\shO_X(U_i) \otimes A) \lrg 0.
$$
Moreover, we note that this is a central extension.

\end{oss}

\begin{teo}\label{teo sernesi first ordine=H^1(X,T_X)}
Let $X$ be a complex manifold. Then  there is a 1-1
correspondence:
$$
k: \, \frac{\{\mbox {first order deformations of X}\}
}{\mbox{isomorphism}} \lrg  \check{H}^1(X,\Theta_X),
$$
called the Kodaira-Spencer correspondence, where
$\Theta_X=\Hom(\Omega^1_X,\shO_X)=\operatorname{Der}_\K(\shO_X,\shO_X)$,
is the hololmorphic tangent bundle of $X$.

Moreover, $k(\xi)=0$ if and only if $\xi$ is the trivial
deformation class.
\end{teo}
\begin{proof}
For completeness we repeat this proof from
\cite[Proposition~II.1.6]{bib Sernesi} where all details are
available.

Let $X_{\epsi}$ be a first order deformation of $X$:
\begin{center}
$\xymatrix{X  \ar[r] \ar[d] & X_{\epsi}\ar[d]  \\
           \Spec(\C) \ar[r]  & \Spec(\C[\epsi]),   \\ }$
\end{center}
and $\mathcal{U}=\{U_i\}_{i\in I}$ be a Stein  open cover of $X$
such that $U_{ij}=U_i \cap U_j$ and $U_{ijk}=U_i \cap U_j \cap
U_k$ are Stein for every $i,j$ and $k \in I$.

For any open $U_i$  the deformations ${X_{\epsi}}_{|U_i} $ are
trivial, then for each $i \in I$ there exist isomorphisms of
deformations
$$
\theta_i :\, U_i \times \Spec(\C[\epsi]) \lrg {X_{\epsi}}_{|U_i}.
$$
Therefore, for all $i$ and $j \in I$, the composition
$$
\theta_{ij}={\theta_i}^{-1}\theta_j: U_{ij} \times
\Spec(\C[\epsi]) \lrg U_{ij} \times \Spec(\C[\epsi])
$$
is an automorphism of the trivial deformation $U_{ij} \times
\Spec(\C[\epsi])$ of the Stein open subset $U_{ij}$.

Applying Lemma~\ref{lemm serneII.1.5 automor-derivazio}, we
conclude that there exists an element $d_{ij} \in
\Gamma(U_{ij},\Theta_X)$ corresponding to $\theta_{ij}$, for all
$i$ and $j \in I$.

Moreover, on each $U_{ijk}$ the following equality holds
$$
\theta_{jk}{\theta_{ik}}^{-1} \theta_{ij}^{-1}=
{\theta_j}^{-1}\theta_k {\theta_k}^{-1}\theta_i
{\theta_i}^{-1}\theta_j=id_{|U_{ijk}\times \Spec(\kepsi)}.
$$
Therefore, applying again Lemma~\ref{lemm serneII.1.5
automor-derivazio} yields
$$
d_{jk}-d_{ik}+d_{ij}=0,
$$
that is, $\{d_{ij}\}$ is a \v{C}ech 1-cocycle and so it defines an
element in $\check{H}^1(X,\Theta_X)$.

\noindent It can be checked that this element does not  depend on
the choice of the open cover $\mathcal{U}$.

Let $ X_{\epsi}'$ be another first order deformation
\begin{center}
$\xymatrix{X  \ar[r] \ar[d] & X_{\epsi}'\ar[d]  \\
           \Spec(\C) \ar[r]  & \Spec(\C[\epsi]),   \\ }$
\end{center}
and $\phi$ an  isomorphism   of deformations: $\phi:X_\epsi\lrg
X'_\epsi$.

Then, for each $i \in I$, there exists an induced automorphism
$$
\alpha_i={\theta'_i}^{-1} \circ \phi_{| U_i}\circ \theta_i:
$$
$$
U_i \times \Spec(\C[\epsi])
\stackrel{\theta_i}{\lrg}{X_{\epsi}}_{|U_i}
\stackrel{\phi_{|U_i}}{\lrg}{X'_{\epsi}}_{|U_i}
\stackrel{{\theta'_i}^{-1}}{\lrg} U_i \times \Spec(\C[\epsi])
$$
and so a corresponding element $a_i \in \Gamma(U_i,\Theta_X)$.

Therefore, we have ${\theta'_i} \alpha_i=  \phi_{| U_i} \theta_i$
and
$$
({\theta'_i} \alpha_i)^{-1}({\theta'_j} \alpha_j)={\theta_i}^{-1}
 \phi_{| U_{ij}}^{-1}  \phi_{| U_{ij}} \theta_j={\theta_i}^{-1}
 \theta_j.
$$
This implies
$$
{\alpha_i}^{-1} {\theta_{ij}'}^{-1}\alpha_j= \theta_{ij}
$$
or equivalently
$$
d'_{ij}+a_j-a_i =d_{ij}.
$$
In conclusion, the \v{C}ech cocycles  $\{d_{ij}\}$ and
$\{d'_{ij}\}$ are cohomologous and so they represent the same
element in $\check{H}^1(X,\Theta_X)$.

\smallskip

Conversely, let $\theta \in \check{H}^1(X,\Theta_X)$ and
$\{d_{ij}\}\in Z^1(\mathcal{U},\Theta_X) $ be a representative of
$\theta$ with respect to an  open Stein cover $\{\mathcal{U}\}$.
By Lemma~\ref{lemm serneII.1.5 automor-derivazio}, we can
associate with each $d_{ij}$ an automorphism $\theta_{ij}$ of the
trivial deformation $ U_{ij} \times \Spec(\C[\epsi])$. Since the
element $\{d_{ij}\}$ satisfies the cocycle condition,
$\theta_{ij}$   satisfies this condition too:
$$
\theta_{jk}{\theta_{ik}}^{-1} \theta_{ij}^{-1}= id_{|U_{ijk}\times
\Spec(\C[\epsi])}.
$$
Using these automorphisms we can glue together the schemes $U_i
\times \Spec(\C[\epsi])$ (see \cite[pag.~69]{bib hartshorne}) to
obtain a scheme $X_\epsi$ that is a first order deformation of
$X$.

At this point  the last assertion is clear.

\end{proof}

\begin{teo}\label{teo sernesi ostruz DEF_X=H^2(X,T_X)}
$\check{H}^2(X,\Theta_X)$ is a complete obstruction space for
$\Def_X$.
\end{teo}
\begin{proof}
For completeness, we take the proof from
\cite[Proposition~II.1.8]{bib Sernesi} where  details are
available.

Let $\mathcal{U}=\{U_i\}_{i\in I}$ be an open  Stein  cover of $X$
such that $U_{ij} $ and $U_{ijk} $ are Stein for all, $i,j$ and $k
\in I$ and
$$
0 \lrg J \lrg B \lrg A \lrg 0
$$
be a small extension in $\Art$.

Let $X_A$ be an infinitesimal deformation of $X$ over $\Spec(A)$.
Then we have isomorphisms
$$
\theta_i:U_i\times \Spec(A) \lrg {X_A}_{|U_i}
$$
such that $\theta_{ij}:={\theta_i}^{-1} \theta_j$ are
automorphisms of the trivial deformations $U_{ij} \times \Spec(A)$
and $\theta_{jk}{\theta_{ik}}^{-1} \theta_{ij}=id_{|U_{ijk}\times
\Spec(A)}$.

To define a deformation $X_B$ that lifts  the deformation $X_A$ is
necessary and sufficient to give automorphisms
$\{\tilde{\theta}_{ij}\}$ of the trivial deformation $U_{ij}
\times \Spec(B)$ such that
\begin{itemize}
  \item[$i)$] $\{\tilde{\theta}_{ij}\}$  glues together:
  $\tilde{\theta}_{jk} {\tilde{\theta}_{ik}}^{-1}
   \tilde{\theta}_{ij}=id_{|U_{ijk}\times \Spec(B)}$,
  \item[$ii)$] $\{\tilde{\theta}_{ij}\}$  lifts $\{ \theta _{ij}\}$:
$\tilde{\theta}_{ij}$ restricts to $\theta_{ij}$ on $U_{ij} \times
\Spec(A)$.
\end{itemize}

Let us choose automorphisms $\{\tilde{\theta}_{ij}\}$ that satisfy
condition $ii)$. Then the  automorphisms
$$
\tilde{\theta}_{ijk}=\tilde{\theta}_{jk}
{\tilde{\theta}_{ik}}^{-1}
   \tilde{\theta}_{ij}
$$
are automorphisms of the trivial deformation that restrict  to the
identity on $U_{ijk}\times \Spec(A)$. By Lemma~\ref{lemm
serneII.1.5 automor-derivazio}, there exists $\{ \tilde{d}_{ijk}\}
\in \Gamma(U_{ijk},\Theta_X)\otimes J$ that corresponds to
$\tilde{\theta}_{ijk}$. An easy calculation shows that $\{
\tilde{d}_{ijk}\}$ is a \v{C}ech cocylce  and so $\{
\tilde{d}_{ijk}\} \in  Z^2(\mathcal{U},\Theta_X)\otimes J$.

Next, let $\{\overline{\theta}_{ij}\}$ be different automorphisms
of the trivial deformations $U_{ij} \times \Spec(B)$ that satisfy
condition $ii)$. As above, let $\overline{d}_{ijk}$ be the
derivations corresponding  to $\overline{\theta}_{ijk}$.

The automorphisms  $\overline{\theta}_{ij}
{\tilde{\theta}_{ij}}^{-1}$ of $U_{ij} \times \Spec(B)$  restrict
to the identity on $U_{ij} \times \Spec(A)$ for each $U_{ij}$ and
so, again by Lemma~\ref{lemm serneII.1.5 automor-derivazio}, they
correspond  to some $\{d_{ij}\} \in \Gamma(U_{ij},\Theta_X)\otimes
J$.

Therefore, $$
\overline{d}_{ijk}=\tilde{d}_{ijk}+d_{jk}-d_{ik}+d_{ij}.
$$
This implies that the \v{C}ech cocycles $\{ \tilde{d}_{ijk}\}$ and
$\{ \overline{d}_{ijk}\}$ are cohomologous and so their cohomology
classes coincide in a well defined element $v_e(X_A)$ in
$\check{H}^2 (X,\Theta_X)\otimes J$:
$$
v_e(X_A):=[\{ \tilde{d}_{ijk}\}]=[\{ \overline{d}_{ijk}\}] \in
\check{H}^2 (X,\Theta_X)\otimes J.
$$
Moreover, the class $v_e(X_A)$ is zero if and only if the
collection of automorphisms also satisfies condition $i)$ and it
is equivalent to the existence of a lifting $X_B$ of the
deformation $X_A$.

\end{proof}

\section{DGLAs  and deformation functor}
\label{sezio DGLA e funtore DEF_L}

In this section  we study the deformation functor associated with
a differential graded Lie algebra (DGLA).

In particular, we give the fundamental definition of  a DGLA
(Definition~\ref{definizio DGLA}).

We also introduce the Maurer-Cartan functor $\MC_L$
(Definition~\ref{def Maurer-C funcotr MC_L of L})  and the
deformation functor $\Def_L$ associated with a DGLA $L$
(Definition~\ref{dhef funtore defo DLGA DEF_L}).

We start by defining the differential graded vector spaces.

\subsection{Differential graded vector spaces}

Let $\K$ be a fixed field of characteristic 0. Unless otherwise
stated, all vector spaces, linear maps, tensor products etc. are
intended over $\K$.

\medskip

Every \emph{graded vector space} is a $\Z$-graded vector space
(over $\K$). If $V=\oplus_i V^i$ is a graded vector space and $a
\in V$ is a  homogeneous element, then we denote by $\deg_V(a)$
the degree of $a$ in $V$; we will also use the notation $\deg(
a)=\bar{a}$, when $V$ is clear from the context.

The morphisms of graded vector space are   degree preserving
linear maps.

Given two graded vector spaces $V$ and $W$, we define $\Hom_\K
^n(V,W)$ as the vector space of $\K$-linear maps $f:V \lrg W$,
such that $f(V^i) \subset W^{i+n}$, for each $i \in \Z$.

Let $V$ be a graded vector space, then $V[n]$ is the complex $V$
with degrees shifted by $n$. More precisely, for $\K[n]$ we have
$$
\K[n]^i=
  \begin{cases}
    \K & \text{if $i+n=0$}, \\
    0 & \text{otherwise}.
  \end{cases}
$$
Then $V[n]=\K[n] \otimes V$, which implies $V[n]^i=V^{i+n}$.

\begin{oss}
We note that there exist isomorphisms
$$
\Hom^n_\K(V,W)\cong \Hom^0_\K(V[-n],W)\cong \Hom^0_\K(V,W[n]).
$$

\end{oss}

A \emph{differential graded vector space } is a pair $(V,d)$ where
$V=\oplus V^i$ is a graded vector space and $d$ is a differential
of degree 1 ($d:V^i \lrg V^{i+1}$ and $d \circ d=0$).

\medskip

For every differential graded vector space $(V,d)$ we use the
standard notation $Z^i(V)=\ker(d\colon V^i\to V^{i+1})$,
$B^i(V)=\image(d\colon V^{i-1}\to V^{i})$ and
$H^i(V)=Z^i(V)/B^i(V)$.

\medskip

 A morphism of differential graded vector
spaces is a degree preserving linear map that commutes with the
differentials.

A morphism is a \emph{quasi-isomorphism} if it induces
isomorphisms in cohomology.

\begin{exe}
Given $(V,d )$, then for each $i \in \Z$, the shifted differential
graded vector space $(V[i],d_{[i]})$ is defined as:
$$
V[i]=\bigoplus_j V[i]^j=\bigoplus_j V^{i+j} \qquad \mbox { and }
\qquad d_{[i]}=(-1)^i d.
$$

\end{exe}

\begin{exe}\label{exe dhef Hom^*(V,W) grade vect spac}
If $(V,d_V)$ and $(W,d_W)$ are differential graded vector spaces,
then we can define a new differential graded vector space
$$
\Hom^*(V,W)=\bigoplus_{n\in \Z}\Hom_\K^n(V,W)
$$
with natural differential $d'$ given by
$$
d'(f):= d_W f -(-1)^{\deg(f)} f \,d_V.
$$
Moreover, for each $i$, there exists the following isomorphism
$$
H^i(\Hom^*(V,W))\cong \Hom^i(H^*(V),H^*(W)).
$$
\end{exe}

\begin{exe}\label{exe definizio Htp(V,W)}
Given $(V,d_V)$ and $(W,d_W)$, we can also define the following
differential graded vector space
$$
\Htp(V,W)=\Hom^*(V[1],W)=\bigoplus_i  \Htp^i(V,W),
$$
with
$$
\Htp^i(V,W)=\Hom^i(V[1],W)=\Hom^{i-1}(V,W)
$$
and differential $\delta$:
$$
\delta(f)=d_W(f)-(-1)^if d_{V[1]}=d_W f+(-1)^if d_{V} \qquad
\forall \, f \in \Htp^i(V,W).
$$
We will use these differential graded vector spaces in the last
chapter  (Section~\ref{sezio HTP}).

\end{exe}

\subsection{Differential graded Lie algebras (DGLAs)}

\begin{dhef}\label{definizio DGLA}
A \emph{differential graded Lie algebra} (DGLA for short) is a
triple $(L,[\, ,\, ],d)$, where $\displaystyle (L=\bigoplus_{i\in
\Z}L^i,d)$ is a differential graded vector space and  $[\, , \,]:
L\times L\to L$ is a bilinear map, called bracket,  satisfying the
following conditions:

\begin{itemize}

\item [1. ] the bracket $[\, ,\,]$ is homogeneous and graded
skewsymmetric; i.e., $[L^i,L^j]\subset L^{i+j}$ and
$[a,b]+(-1)^{\deg(a)\deg(b)}[b,a]=0$, for every homogeneous $a$
and $b$.

\item [2. ] Every triple of homogeneous elements satisfies the graded
Jacobi identity
$$
[a,[b,c]]=[[a,b],c]+(-1)^{\deg(a)\deg(b)}[b,[a,c]].
$$
\item [3. ] $d(L^i) \subseteq L^{i+1}$, $d\circ d=0$  and
$$
d[a,b]=[da,b]+(-1)^{\deg(a)}[a,db].
$$
\end{itemize}
\end{dhef}

The last property is called  Leibniz's rule and, in particular, it
implies that the bracket induces a structure of differential
graded Lie algebra (with zero differential) on the cohomology
$H^*(L)=\oplus_iH^i(L)$ of a DGLA $L$.

\begin{exe}
If $L=\oplus L^i$ is a DGLA, then $L^0$ is a Lie algebra in the
usual sense; vice-versa, every Lie algebra is a differential
graded Lie algebra concentrated to degree 0.

\end{exe}

\begin{oss}\label{oss [a,a]=0 pari [[a,a],a]=0 dispari}
If  the degree of $a$ is even then $[a,a]=0$ (graded
skew-symmetry). If the degree of $a$ is odd, then
$[a,[a,b]]=\frac{1}{2}[[a,a],b]$ for all $b \in L$ (graded
Jacobi). In particular, $[a,[a,a]]=0$; moreover, $[a,b]=[b,a]$ for
all $a$ and $b \in L$ of odd degree.
\end{oss}

\begin{dhef}
A morphism of differential graded Lie algebras $\varphi: L \lrg M$
is a linear map that preserves degrees and commutes with brackets
and differentials; written in details we have
\begin{itemize}
  \item[-] $\varphi(L^i)\subseteq M^i$, for each $i$;
  \item[-] $\varphi(d_L a)=d_M(\varphi(a))$, for each $a \in L$;
  \item[-] $\varphi([ a,b])=[\varphi(a),\varphi(b)]$,
            for each $a,b \in L$.
\end{itemize}

\end{dhef}

A \emph{quasi-isomorphism} of DGLAs is a morphism that induces
isomorphisms in cohomology. Two DGLAs $L$ and $M$ are
\emph{quasi-isomorphic} if they are equivalent under the
equivalence relation $\sim$ defined by: $L\sim M$ if there exists
a quasi-isomorphism $\phi:L \lrg M$.

A DGLA  $L$ is \emph{formal} if it is quasi-isomorphic to its
cohomology graded vector space $H^*(L)$.

\begin{oss}
The following DGLAs are isomorphic:
$$
(L,[\, ,\, ],d) \cong (L,-[\, ,\, ],d) \cong (L,[\, ,\, ],-d)
\cong (L,-[\, ,\, ],-d).
$$
Actually, the morphism $\varphi=-id$ gives an isomorphism between
$(L,[\, ,\, ],d)$ and $(L,-[\, ,\, ],d)$, whereas
$\varphi(-)=(-1)^{\deg(-)}id$ is an isomorphism between $(L,[\,
,\, ],d)$ and $(L,[\, ,\, ],-d)$.
\end{oss}

\begin{dhef}
A linear map $f:L \lrg L$ is called a \emph{derivation} of degree
$n$ if $f(L^i)\subset L^{i+n}$ and it satisfies the graded Leibniz
rule:
$$
f([a,b])=[f(a),b]+(-1)^{n\overline{a}}[a,f(b)].
$$
\end{dhef}

The graded  Leibnitz rule implies that the differential $d$ is a
derivation of degree $1$.

\begin{oss}
Let $a \in L^i$ and consider the operator
$$
[a,\,]:L \lrg L,
$$
$$
[a,\, ](x)=[a,x], \qquad \forall \ x \in L.
$$
Then, the graded Jacobi identity implies that $[a,\, ]$ is a
derivation of degree $i =deg(a)$.
\end{oss}

\begin{exe}\label{exe def Ncon N^1+B^1=L^1}
Let $L$ be a DGLA and consider the vector space decomposition
$L^1=N^1\oplus B^1(L)$. Then the graded vector space $N=\oplus
N^i$ with
$$
\begin{cases}
N^i=0   & \text{ for  $i \leq 0$}\\
N^1=N^1  &\text{ for  $i=1$}\\
N^i=L^i  &\text{ for  $i \geq 2$}
\end{cases}
$$
is a sub-DGLA of $L$.
\end{exe}
\begin{exe}\label{exe def DGLA N=L+Kd}
Given a DGLA $(L= \oplus L^i,[\, ,\, ],d)$ we can associate with
it a new DGLA $(L'=\oplus {L'}^i,[\, ,\, ]',d')$ where
$$
\begin{cases}
{L'}^i=L^i &\text{ for  $i \neq 1$}\\
{L'}^1=L^1\oplus \K\, d & \text{ for  $i=1$},
\end{cases}
$$
$$
[v+ad,w+bd]'=[v,w]+a\,d(w)-(-1)^{\deg(v)}b\,d(v)
$$
and
$$
 d'(v+a\,d)=[d,v+ad]'=dv,
$$
for any $v,w \in L^i$ and $a,b \in \K$.
\end{exe}

\begin{exe}\label{exe definizio M[t,dt]}
Let $M$ be a DGLA. Then $M[t,dt]=M \otimes \K[t,dt]$ is a DGLA,
where $\K[t,dt]$ is the differential graded algebra of  polynomial
differential forms over the affine line. More precisely,
$\K[t,dt]=\K[t] \oplus \K[t]dt$, where $t$ has degree $0$ and $dt$
has degree 1. As vector space $M[t,dt]$ is generated by elements
of the form $mp(t)+ nq(t)dt$, with $m,n \in M$ and $p(t),q(t)\in
\K[t]$. The differential and  the bracket on $M[t,dt]$ are defined
as follows:
$$
d(mp(t)+ nq(t)dt)=(dm)p(t) + (-1)^{\deg( m)} m p'(t)dt +
(dn)q(t)dt,
$$
$$
[mp(t),nq(t)]=[m,n]p(t)q(t), \ \ \
[mp(t),nq(t)dt]=[m,n]p(t)q(t)dt.
$$

\bigskip

For every $a \in \K$ define the \emph{evaluation morphism} in the
following way
$$
e_a:M[t,dt] \lrg M,
$$
$$
e_a(\sum m_it^i +n_it^i dt)=\sum m_i a^i.
$$
The evaluation morphism is a morphism of DGLAs which is a left
inverse of the inclusion and    is a surjective quasi-isomorphism
for each $a$.
\end{exe}

\begin{exe}\label{exe LxA is DGLA e^[a,] auto}
If $L$ is a DGLA and $B$ is a commutative $\K$-algebra then $L
\otimes B$ has a natural structure of DGLA, given by
$$
[l\otimes a , m \otimes b]=[l,m]\otimes ab;
$$
$$
d(l\otimes a)=dl \otimes a.
$$
\end{exe}
If $B$ is also nilpotent (for example $B=m_A$ the maximal ideal of
a local artinian $\K$-algebra $A$)  then $L\otimes B$ is a
nilpotent DGLA. Therefore,  for every $a \in L^0 \otimes B$, we
can define an automorphism of the DGLA $L \otimes B$:
$$
e^{[a,-]}:=\sum_{n\geq 0} \frac{ [a,-]^n}{n!}\ : L \otimes B \lrg
L \otimes B,
$$
where
$$
[a,-]:L\otimes B \lrg L \otimes B,
$$
$$
[a,-](b):=[a,b]
$$
is a nilpotent derivation  of degree zero (since $ [a,-]([b,c])=
[[a,-](b),c] +[b,[a,-](c)]$).

\begin{exe}
Let $(L,d)$ be a DGLA and $\operatorname{Der}^i(L,L)$ the space of
derivations of $L$ of degree $i$. Then
$\operatorname{Der}^*(L,L)=\bigoplus_i\operatorname{Der}^i(L,L)$
is a DGLA with bracket
$$
[f,g]=fg-(-1)^{\overline{f}\,\overline{g}}gf,
$$
and differential $d'$ given by
$$
d'(f)=[d,f].
$$

\end{exe}

\subsection{Maurer-Cartan functor associated with a DGLA}
\begin{dhef}\label{def Maurer-C funcotr MC_L of L}
The \emph{Maurer-Cartan equation} in a differential graded Lie
algebra  $L$ is
$$
dx+\displaystyle\frac{1}{2} [x,x]=0, \qquad x \in L^1.
$$
The solutions of this equation are called the \emph{Maurer-Cartan
elements} of the DGLA $L$.

\end{dhef}

\begin{oss}\label{oss d(-)+[x,-]^2=0 se x in MC}
Let $x\in L^1$ be an element of degree one in the DGLA $L$  and
consider the operator of degree one
$$
d(-)+[x,-]: L \lrg L,
$$
$$
\qquad (d(-)+[x,-])(a)=da+[x,a] \qquad \forall \ a \in L.
$$

Then, $d(-)+[x,-]$ is a differential ($(d(-)+[x,-])^2=0$) if and
only if   $x$ satisfies the Maurer-Cartan equation. Indeed, using
Remark~\ref{oss [a,a]=0 pari [[a,a],a]=0 dispari},
$$
(d(-)+[x,-])^2(a)=(d(-)+[x,-])(da+[x,a])=
$$
$$
d^2 a+d[x,a] +[x,da] +[x,[x,a]]=
$$
$$
[dx,a] + (-1)^{deg(x)}[x,da]+[x,da]+\frac{1}{2}[[x,x],a]=
$$
$$
[dx+ \frac{1}{2} [x,x],a].
$$
\end{oss}

\begin{oss}\label{oss MC(x)=0 [x+d,x+d]'=0}
Let ${L'}$ be the DGLA of Example~\ref{exe def DGLA N=L+Kd} (with
${L'}^1=L^1\oplus \K\,d$)  then
$$
dx+\displaystyle\frac{1}{2}  [x,x]=0 \quad \mbox{ if and only if }
\quad [x+d,x+d]'=0.
$$

\end{oss}

The  previous definition   led  to the following definition of the
Maurer-Cartan functor.

\begin{dhef}
Let $L$ be a DGLA; then the \emph{Maurer-Cartan functor associated
with $L$} is
$$
MC_L:\Art \lrg \Set,
$$
$$
MC_L(A)=\{ x \in L^1\otimes m_A \ |\ dx+ \displaystyle\frac{1}{2}
[x,x]=0 \}.
$$
\end{dhef}

\begin{oss}

A morphism  $\phi:L \lrg M$  of DGLA preserves bracket and
differential, therefore, it induces a  morphism of functors $\phi:
\MC_L \lrg \MC_M$.
\end{oss}

\begin{oss}\label{oss MC_L is homogeneous}
We note that $\MC_L$ is  a  \emph{homogeneous} functor, since
$\MC_L(B \times _A C) \cong \MC_L(B)\times_{\MC_L(A)}\MC_L(C)$ for
each pair of morphisms $\gamma:C \lrg A$ and $\beta:B \lrg A$,
with $\beta$ surjective.
\end{oss}

\begin{oss}\label{oss calcolo tangente MC_L}
By definition, the \emph{tangent space} of $\MC_L$ is:
$$
t_{\MC_L}:=\MC_L(\kepsi)=\{ x \in L^1 \otimes \K \epsi\, | \, dx+
\displaystyle\frac{1}{2}[x,x]=0 \}\cong
$$
$$
\{ x \in L^1 \, | \, dx  =0 \}=Z^1(L).
$$

\end{oss}

\begin{lem}\label{lem calcolo ostruzione MC_L=H^2(L)}
$H^2(L)$ is a complete obstruction space for $\MC_L$.
\end{lem}
\begin{proof}
Let
$$
e: \qquad 0 \lrg J \lrg B \stackrel{\alpha}{\lrg} A\lrg 0
$$
be a small extension in $\Art $ and $ x \in \MC_L(A)$.

We want to define a map $v_e:\MC_L(A) \lrg H^2(L) \otimes J$.

Let $ \tilde{x} \in L^1  \otimes m_B$ be a lifting of $x$ and
define
$$
h=d\tilde{x} +\displaystyle\frac{1}{2} [\tilde{x},\tilde{x}] \in
L^2  \otimes m_B.
$$
In general, $\tilde{x}$ does not satisfy the Maurer-Cartan
equation and so $h$ is in general different from zero.

It turns out that $\alpha(h)=dx+\displaystyle\frac{1}{2}[x,x]=0$
and so $h \in L^2 \otimes J$. Moreover,
$$
dh=d^2\tilde{x}+\displaystyle\frac{1}{2}[d\tilde{x},\tilde{x}]-
\displaystyle\frac{1}{2}[\tilde{x}, d\tilde{x}]=\footnote{We have
$- \displaystyle\frac{1}{2}[\tilde{x}, d\tilde{x}]=-
\displaystyle\frac{1}{2} (-(-1)^2[d\tilde{x},\tilde{x}] )=
\displaystyle\frac{1}{2}[d\tilde{x},\tilde{x}]. $}
$$
$$
=[d\tilde{x},\tilde{x}]=[h,\tilde{x}]- \displaystyle\frac{1}{2}
[[\tilde{x},\tilde{x}],\tilde{x}].
$$
By definition, $[h ,\tilde{x}] \in [L^2 \otimes J, L^1 \otimes
m_B]=0$  ($e$ is a small extension) and, by Remark~\ref{oss
[a,a]=0 pari [[a,a],a]=0 dispari},
$[[\tilde{x},\tilde{x}],\tilde{x}]=0$.

In conclusion, $dh=0$ and so $h \in H^2(L) \otimes J$.

We note that this class does not  depend on the choice of the
lifting $\tilde{x}$. Indeed, let $y \in L^1 \otimes m_B$ be
another lifting of $x$: $\alpha(y)=\alpha(\tilde{x})=x$. Then
$y=\tilde{x}+t$ for some $t \in L^1 \otimes J$. Using   $[L^1
\otimes J, L^1 \otimes m_B]=0$, we have
$$
h'=dy+ \displaystyle\frac{1}{2}[y,y]=d\tilde{x} +dt+
\displaystyle\frac{1}{2}[\tilde{x}+t,\tilde{x}+t]= d\tilde{x}+
\displaystyle\frac{1}{2}[\tilde{x},\tilde{x}]+dt=h+dt
$$
and so $h$ and $h'$ represent the same class in $H^2(L) \otimes
J$.

Therefore, the following obstruction map is well defined
$$
v_e:\MC_L(A) \lrg H^2(L) \otimes J,
$$
$$
x \longmapsto v_e(x)=[h].
$$
If $[h]=0$ then $h=dq$ for some $q \in L^1 \otimes J$. This
implies that $\overline{x}=\tilde{x}-q$ is a lifting of $x$ that
satisfies the Maurer-Cartan equation, i.e.,
$$
d\overline{x}+ \displaystyle\frac{1}{2}[\overline{x},
\overline{x}]= d\tilde{x}-dq+
\displaystyle\frac{1}{2}[\tilde{x}-q, \tilde{x}-q]=d\tilde{x}-dq+
\displaystyle\frac{1}{2}[\tilde{x},\tilde{x}]=h-dq=0.
$$

Thus, $v_e$ satisfies condition $1$ of Definition~\ref{def
ostruzione space} of obstruction theory. The other property
(change of basis) is an easy calculation.

If $x \in \MC_L(A)$ can be lifted to $x' \in MC_L(B)$ then
$[h]=0$.

In conclusion, $(H^2(L), v_e)$ is a complete obstruction theory
for $\MC_L$.
\end{proof}

\begin{oss} {(\bf About smoothness)}

If $H^2(L)=0$  then $\MC_L$ is smooth.

\smallskip

If $L$ is abelian then $\MC_L$ is \emph{smooth}. Actually, in this
case, $\MC_L(A)=Z^1(L)\otimes m_A$. Moreover, if $B
\twoheadrightarrow A$   then $ Z^1(L)\otimes m_A
\twoheadrightarrow Z^1(L)\otimes m_B$.

\end{oss}

\subsection{Gauge action} \label{sezio gauge action}

\begin{dhef}\label{definizio gauge L DGLA}
Two elements $x$ and  $y \in L^1\otimes m_A$ are said to be
$gauge$ $equivalent$ if there exists $a \in L^0\otimes m_A$ such
that
$$
y=e^a * x:=x+\sum_{n\geq 0}  \frac{ [a,-]^n}{(n+1)!}([a,x]-da).
$$
\end{dhef}
The operator $*$ is called the gauge action of the group $\exp(L^0
\otimes m_A)$ on $L \otimes m_A$; indeed $e^a*e^b*x=e^{a \bullet
b}*x$, where $\bullet$\footnote{$a \bullet b  =a+b +
\displaystyle\frac{1}{2}[a,b]+\frac{1}{12} [a,[a,b]]-\frac{1}{12}
[b,[b,a]]+ \cdots $ } is the Baker-Campbell-Hausdorff product in
the nilpotent DGLA $L \otimes m_A$.

\begin{oss}
For a better understanding of the gauge action, it is convenient
to consider the DGLA $L'$ of Example~\ref{exe def DGLA N=L+Kd}
(with ${L'}^1=L^1\oplus \K\,d$) and the affine embedding
$$
\phi:L^1 \lrg {L'}^1, \qquad \phi(x)=x+d \qquad \forall \  x \in
L^1.
$$
As already observed, $\displaystyle dx+\frac{1}{2}[x,x]=0$ if and
only if $[\phi(x),\phi(x)]'=0$.

As in Example~\ref{exe LxA is DGLA e^[a,] auto}, for each $A \in
\Art$ and  $a \in {L'}^0\otimes m_A$, we can consider the
exponential of the adjoint action $e^{[a,\,]'}:{L'}^1\otimes m_A
\lrg {L'}^1\otimes m_A$. Using the embedding $\phi$ this action
induces the gauge action of $L^0\otimes m_A$ on $L^1\otimes m_A$.
Actually, for any $a \in  L^0\otimes m_A$ and $x \in L^1\otimes
m_A$
$$
\phi^{-1} (e^{[a,\,]'}\phi(x))=e^{[a,]' }(x+d)-d =
$$
$$
 x + \sum_{n\geq 1}
\frac{ {[a,-]'}^n}{(n)!}(x+d)=x+  \sum_{n\geq 0} \frac{
{[a,-]'}^{n+1}}{(n+1)!}(x+d)=
$$
$$
x+  \sum_{n\geq 0} \frac{ {[a,-]'}^{n }}{(n+1)!}([a,x]-da)= e^a*x.
$$

\end{oss}

\begin{exe}\label{exe azione gauge se []=0}
Let $J$ be an ideal of $A \in \Art$ ($J \subset m_A$) such that
$J\cdot m_A=0$ (for example, $J$ may be the kernel of a small
extension).

If $x \in L^1 \otimes J$ then, for each $a \in L^0 \otimes m_A$ we
have
$$
e^a *x=x+\sum_{n=0}^{\infty} \frac{ [a,-]^n}{(n+1)!}([a,x]-da)=
x+\sum_{n=0}^{\infty} \frac{ [a,-]^n}{(n+1)!}( -da)=x+e^a*0,
$$
or, in general, if   $y \in L^1 \otimes m_A$ then
$e^a*(x+y)=x+e^a*y$.

If $a \in L^0\otimes J$ then, for each $x \in L^1 \otimes m_A$:
$$
e^a*x= x+\sum_{n=0}^{\infty} \frac{
[a,-]^n}{(n+1)!}([a,x]-da)=x-da.
$$
In general, if   $b \in L^0 \otimes m_A$,  then
$e^{a+b}*x=e^b*x-da.$
\end{exe}

We note that
\begin{equation}\label{equa e^a*x=x se e solo se [a,x]=da}
e^a*x=x \qquad \mbox{ if and only if } \qquad  [a,x]=da.
\end{equation}
Actually, $e^a*x=x$ if and only if $\displaystyle 0=
\frac{e^{[a,-]}-id}{[a,-]} ([a,x]-da)$. Applying the inverse of
the operator $\displaystyle  \frac{e^{[a,-]}-id}{[a,-]}$, yields
$e^a*x=x$ if and only if $ [a,x]-da=0$.

\begin{oss}
The solutions of the Maurer-Cartan equation are preserved under
the gauge action.

Actually, we have
$$
d(e^a*x)=d'(e^{[a,\, ]' }(d+x)-d)= [d,e^{[a,\, ]'
}(d+x)-d]'=[d,e^{[a,\, ]' }(d+x) ]'
$$
and using Remark~\ref{oss MC(x)=0 [x+d,x+d]'=0}
$$
[e^a*x,e^a*x]=[e^{[a,\, ]' }(d+x)-d,e^{[a,\, ]' }(d+x)-d]' =
$$
$$
[e^{[a,\, ]' }(d+x),e^{[a,\, ]' }(d+x)]'-2[d,e^{[a,\, ]' }(d+x)]'=
$$
$$
e^{[a,\, ]'}[d+x,d+x]'-2[d,e^{[a,\, ]' }(d+x) ]'=-2[d,e^{[a,\, ]'
}(d+x) ]'.
$$
Therefore, $$ d(e^a*x)+ \displaystyle\frac{1}{2} [e^a*x,e^a*x]= 0.
$$

\end{oss}

Finally, for each $x \in \MC_L (A) $, we define the
\emph{irrelevant stabilizer} of $x$:
$$
Stab_A(x)=\{e^{dh+[x,h]}|\, h \in L^{-1}\otimes m_A\} \subset
\exp(L^0\otimes A ).
$$
The name irrelevant stabilizer is due to the fact that
$e^{dh+[x,h]}*x=x$. Actually, $dh+[x,h]$ satisfies condition
(\ref{equa e^a*x=x se e solo se [a,x]=da}), i.e.,
$$
[dh+[x,h],x]=[dh,x]+[[x,h],x]=d[h,x]+[h,dx]+\frac{1}{2} [h,[x,x]]=
$$
$$
d[h,x]+[h,dx+\frac{1}{2}  [x,x]]=d[h,x]=d(dh+[x,h]).
$$

Moreover, we  observe that $Stab_A(x)$ is a subgroup of
$\exp(L^0\otimes A )$ and that for any $a \in L^0\otimes A $
$$
e^a Stab_A(x)e^{-a}=Stab_A(y),  \qquad \mbox{ with } \qquad
y=e^a*x.
$$

\subsection{Deformation functor associated with a
DGLA}\label{section defi funtor DEF_L of DGLA}

\begin{dhef}\label{dhef funtore defo DLGA DEF_L}
The \emph{deformation functor} associated with a differential
graded Lie algebra  $L$ is:
$$
\Def_L:\Art \lrg \Set,
$$
$$
\Def_L(A)=\frac{MC_L}{\exp\,(L^0 \otimes m_A)}.
$$
\end{dhef}

In this case too a morphism of DGLAs $\phi: L \lrg M$ induces a
morphism of the associated functors $\phi: \Def_L \lrg \Def_M$.

The name deformation functor is justified by the following
proposition.

\begin{prop}\label{prop DEF_L is funtore deformazione}
$\Def_L$  is a deformation funcotr, i.e., it satisfies the
conditions of Definition~\ref{defin deformation funtore}.
\end{prop}

\begin{proof}
If $A=\K$, then it is clear that $\Def_L(B \times C) = \Def_L(B)
\times \Def_L(C)$ and so condition $ii)$ of Definition~\ref{defin
deformation funtore} is satisfied.

Next, let $\beta:B \lrg A$ and $\gamma: C \lrg A$ be morphisms in
$\Art$ with $\beta$ surjective. Let $(l,m) \in \Def_L(B)
\times_{\Def_L(A)} \Def_L(C)$ and $\tilde{l} \in MC_L(B)$ and
$\tilde{m} \in \MC_L(C)$ be liftings of $l$ and $m$, respectively,
such that $\beta(\tilde{l})=\gamma(\tilde{m})\in \Def_L(A)$.
Therefore, there exists $a \in L^0 \otimes m_A$ such that
$e^a*\beta(\tilde{l})=\gamma(\tilde{m})$. Let $b \in L^0 \otimes
m_B$ be a lifting of $a$. By replacing $\tilde{l}$ with its gauge
equivalent element $l'=e^b*\tilde{l}$ we can suppose
$\beta(l')=\gamma(\tilde{m})$ in $\MC_L(A)$. By Remark~\ref{oss
MC_L is homogeneous}, $\MC_L$ is homogeneous and so there exists
$n \in \MC_L(B \times_A C)$ that lifts $(l',\tilde{m})$. This
implies that
$$
\Def_L(B \times_A C) \lrg  \Def_L(B) \times_{\Def_L(A)} \Def_L(C)
$$
is surjective. Hence condition $i)$ of Definition~\ref{defin
deformation funtore} also holds.

\end{proof}

\medskip
\begin{oss}\label{oss tangente DEF_L =H^1}
By definition, the \emph{tangent space} of $\Def_L$ is:
$$
t_{\Def_L}:=\Def_L(\kepsi)=\frac{\{ x \in L^1 \otimes \K \epsi\, |
\, dx  =0 \}}{\{da \, | a \in L^0 \otimes \K \epsi \} }\cong
$$
$$
H^1(L).
$$
In general, if $L \otimes m_A$ is abelian then $\Def_L(A)=H^1(L)
\otimes m_A$.

\end{oss}

\begin{lem}
The projection $\pi:\MC_L \lrg \Def_L$ is a  smooth morphism of
functors.
\end{lem}
\begin{proof}
Let $ \beta: B \lrg A$ be a surjection in $\Art $; we prove that
$$
\MC_L(B) \lrg \Def_L(B)\times_{\Def_L(A)}\MC_L(A),
$$
induced by
\begin{center}
$\xymatrix{\MC_L(B)  \ar[r]^\beta \ar[d]_\pi &  \MC_L(A)   \ar[d]^\pi \\
          \Def_L(B) \ar[r]_\beta  & \Def_L (A),   \\ }$
\end{center}

is surjective.

Let $( b,a) \in \Def_L(B)\times_{\Def_L(A)}\MC_L(A)$ and
$\tilde{b}\in \MC_L(B)$ be a lifting of $b$. Then
$\beta(\tilde{b})$ and $a$ have a common image in $\Def_L(A)$ and
so $\beta(\tilde{b}) =e^t*a$, for some $t \in L^0 \otimes m_A$.

Let $s \in L^0 \otimes m_B$ be a lifting of $t$ and define
$b'=e^{-s}*\tilde{b} \in \MC_L(B)$. Then $\beta(b')=
e^{-t}*\beta(\tilde{b})=a$ and $b'$ lifts $b$.
\end{proof}

Therefore, by Corollary~\ref{cor f:F_G liscio implica biezione
ostruzio}, $\pi$ induces an isomorphism between universal
obstruction theories.

In conclusion, Lemma~\ref{lem calcolo ostruzione MC_L=H^2(L)}
implies that  $H^2(L)$ \emph{is a complete obstruction space of}
$\Def_L$.

\begin{teo}\label{teo no EXTE iso H^i allro iso DEF}
Let $\phi:L \lrg M$ be a morphism of DGLAs and denote by
$$
H^i(\phi):H^i(L) \lrg H^i(M)
$$
the induced maps in cohomology.
\begin{itemize}
  \item[i)] If $H^1(\phi)$ is surjective (resp. bijective) and
  $H^2(\phi)$ injective, then the morphism $\Def_L \lrg \Def_M$ is
  smooth (resp. \'{e}tale).
  \item[ii)] If  in addition to $i)$  $H^0(\phi)$ is surjective,
  then the morphism $\Def_L \lrg \Def_M$ is an isomorphism.
\end{itemize}
\end{teo}
\begin{proof}
$i)$ follows from Proposition~\ref{prop f liscio se e solo se
OF=0} (resp. Corollary~\ref{cor etale se bietti su tg e inj su
ostru}).   For a proof of $ii)$ see \cite[Theorem~3.1]{bib
manetPISA} (it also follows  from the inverse function
Theorem~\ref{TEO funzione inversa} on the extended case).
\end{proof}

\begin{cor}\label{cor L q iso M allora Def_L=Def_M}
Let $L\lrg M$ be a quasi-isomorphism of DGLAs. Then the induced
morphism $\Def_L \lrg \Def_M$ is an isomorphism.
\end{cor}

%

\begin{cor}
If $H^0(L)=0$, then $\Def_L$ is homogenous.
\end{cor}
\begin{proof}
Let $N$ be the DGLA introduced in Example~\ref{exe def Ncon
N^1+B^1=L^1}. Then the natural inclusion $N \lrg L$ gives
isomorphisms $H^i(N)\lrg H^i(L)$ for each $i \geq 1$. Since
$H^0(L)=0$,   $H^0(N)\lrg H^0(L)$ is surjective. Therefore,
Theorem~\ref{teo no EXTE iso H^i allro iso DEF} $ii)$ implies that
$\Def_N \lrg \Def_L$ is an isomorphism, with $\Def_N\cong \MC_N$
  homogeneous.

\end{proof}

\begin{oss}\label{oss L DGLA governa funtore}
Let $F:\Art \lrg \Set$ be the functor of the infinitesimal
deformations of some algebro-geometric object defined over $\K$.

Then the guiding principle of Kontsevich (see \cite{bib
kontsevich}) affirms the existence of a DGLA $L$ such that $F
\cong \Def_L$. In spite of Corollary~\ref{cor L q iso M allora
Def_L=Def_M}, it is clear that this DGLA is defined only up to
quasi-isomorphism. In this case we say that $L$ \emph{controls}
the deformation functor $F$.

In  Section~\ref{sezio defo manifold defK=defX} we will prove the
existence of a DGLA that controls  the infinitesimal deformations
of $X$ (Theorem~\ref{teo def_k =Def_X}) and in Section~\ref{sez
teo Def_(h,g)Def (f)} the existence of a DGLA that controls the
infinitesimal deformations of a  holomorphic map $f$
(Theorem~\ref{teo esiste hil governa DEF(f)}).
\end{oss}

\chapter{Deformation of complex manifolds}

In the first part of this chapter   we fix notation  and recall
some known facts about complex manifolds   that will be useful in
the sequel.

Therefore,   any book on complex varieties is a good reference for
this chapter (for example \cite{bib Griffone}, \cite{bib
manRENDICONTi}, \cite{bib Voisin}, etc.).

In particular, we wish to recall the \v{C}ech cohomology and
Leray's theorem (Section~\ref{sottosezione cech e leray theorem})
and some properties of K\"{a}hler manifolds (Section~\ref{sezio
kaler manifold}). We also study the maps $f_*$ and $f^*$ induced
by a  holomorphic map $f$ (Section~\ref{sezione f_* and f^*}).

Moreover, we give the fundamental definition of the
\emph{Kodaira-Spencer differential graded Lie algebra}  $KS_X$
associated with a compact complex manifold $X$
(Definition~\ref{def kodaira spencer algebra}), of the
\emph{contraction map} $ \bi $ and of the \emph{holomorphic Lie
derivative} $ \bl $ (Section~\ref{sec cotrction map e Lie deriv}
).

\smallskip

In the second part (Section~\ref{sezio defo manifold defK=defX})
we prove that the functor  $\Def_X$ of the infinitesimal
deformations of a compact complex manifold $X$ (see
Definition~\ref{defin funtore DEF_X}) is isomorphic to the
deformation functor $\Def_{KS_X}$ associated with the
Kodaira-Spencer algebra $KS_X$ of $X$ (Theorem~\ref{teo def_k
=Def_X}).

\begin{teos}
Let $X$ be a complex compact manifold and $KS_X$ its
Kodaira-Spencer algebra.  Then there exists an isomorphism of
functors
$$
  \Def_{KS_X} \lrg \Def_X.
$$
\end{teos}

Therefore, in spite of Remark~\ref{oss L DGLA governa funtore} we
can say that the differential graded Lie algebra of
Kodaira-Spencer $KS_X$ controls the infinitesimal deformations of
a complex compact manifold $X$.

\smallskip

This theorem is well known and a proof based on the theorem of
Newlander-Nirenberg can be found in \cite{bib Catanesecime},
\cite{bib GoMil III}  or more recently in  \cite{bib
manRENDICONTi}. Here we are interested in a simpler proof that
avoids the use of this theorem.

\bigskip

Note that in this chapter  we will work over the complex numbers
and so $\K=\C$.

We also assume that  every variety $X$ is smooth (complex),
compact, and connected.

\section{Differential forms}\label{sezio differential form}

\medskip

Let $X$ be   such  a manifold of dimension $n$ and
$T_{X,\C}=T_X^{1,0} \oplus T_X^{0,1}$ its complex tangent bundle,
with $T_X^{1,0}=: \Theta_X$ the holomorphic tangent bundle and
$T_X^{0,1}=\overline{T_X^{1,0}}$.

This decomposition induces a dual decomposition on the sheaf of
differentiable forms
$$
\shA_X^1=\shA_X^{1,0}\oplus \shA_X^{0,1},
$$
with $\shA_X^{1,0}$ the sheaf of complex differentiable forms of
type $(1,0)$. If $z_1, \ldots , z_n$ are local holomorphic
coordinates on $X$, then $\shA_X^{1,0}$ is generated by the
$dz_i$:   each $\alpha \in \shA_X^{1,0}$ has the form $\alpha=
\sum_i \alpha_i dz_i$, with $\alpha_i \in \shA_X^{0,0}$ for any
$i$.

\noindent  In general, a  $(p,q)$-form $\alpha$ can be locally
written as $\alpha= \sum_{K,J} \alpha_{K,J} dz_K \wedge
d\overline{z}_J$ with $\alpha_{K,J} \in \shA_X^{0,0}$,  $K=(1 \leq
k_1<k_2 < \cdots < k_p \leq n)$ a multi-index of length $p$ and
$J=(1 \leq j_1<j_2 < \cdots < j_q \leq n)$ a multi-index of length
$q$, such that $ dz_K=dz_{k_1} \wedge dz_{k_2} \wedge \cdots
\wedge dz_{ k_p}$ and $d\overline{z}_J=d\overline{z}_{j_1} \wedge
d\overline{z}_{j_2} \wedge \cdots \wedge d\overline{z}_{ j_q}  $.

\smallskip

If $\alpha=f \in \shA_X^{0,0}$, then
$$
df=  \sum_{h=1}^n \frac{\de f }{\de z_h} dz_h + \sum_{h=1}^n
\frac{\de f }{\de \overline{z}_h}d \overline{z}_h= \de f + \debar
f,
$$
with $\de f \in \shA_X^{1,0}$ and $\debar f \in \shA_X^{0,1}$.

In general, for $\alpha=\sum_{K,J} \alpha_{K,J} dz_K \wedge
d\overline{z}_J \in \shA_X^{p,q}$, we have
$$
d\alpha = \sum_{K,J} d \alpha_{K,J}\wedge  dz_K \wedge
d\overline{z}_J = \de \alpha+ \debar \alpha,
$$
with
$$
\de \alpha =\sum_{I,J} \de \alpha_{I,J} \wedge  dz_K \wedge
d\overline{z}_J \in \shA^{p+1,q}
$$
and
$$
\debar \alpha =\sum_{I,J} \debar \alpha_{I,J} \wedge  dz_K \wedge
d\overline{z}_J \in \shA^{p,q+1}.
$$

Obviously, since $d^2=0$ we have $\de\ ^2={\debar\ }^2=\de \debar
+ \debar \de =0$.

\begin{prop}\label{prop esattezza locale debar su X}
Let $\alpha$ be a form of type $(p,q)$, with $q >0$, such that
$\debar \alpha=0$. Then there exists, locally on $X$, a form
$\beta$ of type $(p,q-1)$ such that $\debar \beta=\alpha$.
\end{prop}
\begin{proof}
See \cite[Proposition~2.31]{bib Voisin}.
\end{proof}

\begin{dhef}
$(\shA_X^{*,*},\wedge)$ is the sheaf of graded algebras of
differential forms of $X$, i.e., if $\shA_X^{(p,q)}$ is the sheaf
of differentiable $(p,q)$-forms then
$$
\shA_X^{*,*}:=\bigoplus_i \shA^i_X, \ \qquad \mbox{ with } \qquad
\shA^i_X=\bigoplus_{p+q=i} \shA_X^{p,q}.
$$
\end{dhef}

We use the notation $A_X^{p,q}=\Gamma(X,\shA_X^{p,q})$ for the
vector space of   global sections of  $\shA_X^{p,q}$.

\begin{dhef}
$Der ^*(\shA^{*,*}_X)$ is the sheaf of $\C$-linear derivations on
$\shA_X^{*,*}$; more precisely, if $Der^{a,b}(\shA^{*,*}_X)$ are
the derivations of bi-degree $(a,b)$ then
$$
Der ^*(\shA^{*,*}_X):=\bigoplus_k \bigoplus_{a+b=k}
Der^{a,b}(\shA^{*,*}_X).
$$
\end{dhef}

\noindent We note that $\de$ and $\debar$ are  global sections of
$Der^{1,0}(\shA^{*,*}_X)$ and $Der^{0,1}(\shA^{*,*}_X)$,
respectively.

\begin{oss}
$Der ^*(\shA^{*,*}_X)$ is a \emph{sheaf of differential graded Lie
algebras} with bracket and differential given by the following
formulas:
$$
[f,g]:=f \circ g -(-1)^{\deg(f)\,\deg(g)}g \circ f
$$
and
$$
d(f):=[\de+\debar,f]= \de f +\debar f-(-1)^{\deg(f)} (f \de + f
\debar).
$$
\end{oss}

\bigskip

In particular, fixing $p=0$, $(\shA^{0,*}_X,\wedge)$ is a sheaf of
graded algebras and $\displaystyle
Der^*(\shA_X^{0,*},\shA_X^{0,*}):=\bigoplus_p
Der^p(\shA_X^{0,*},\shA_X^{0,*})$ is a sheaf of DGLAs (in this
case the differential reduces to $d(f)=[\debar,f]=\debar
f-(-1)^{\deg(f)} f \debar $).

\section{K\"{a}hler manifolds}\label{sezio kaler manifold}

This section is devoted to the compact K\"{a}hler manifolds. For
definitions and properties of K\"{a}hler manifolds see for example
\cite{bib Griffone}, \cite{bib manRENDICONTi} or \cite{bib
Voisin}.

We include this section just to prove an important application
(Lemma~\ref{lemma kaler aciclico  q^* A_Y}) of the $\de
\debar$-Lemma (Lemma~\ref{de debar lemma}) that will be
fundamental in the obstruction calculus of the last chapter of
this thesis (Theorem~\ref{teo ostruzi in ker H^2}).


\begin{lem}[$\de \debar$-Lemma] \label{de debar lemma}
Let $X$ be a compact K\"alher manifold and consider the operators
$\de$ and $\debar$ on $A_X$. Then
$$
\image \debar \de =\ker \de \cap \image \debar= \ker \debar \cap
\image \de.
$$
\end{lem}

\begin{proof}
See for example   \cite[Theorem~6.37]{bib manRENDICONTi} and
\cite[Proposition~6.17]{bib Voisin}.

\end{proof}

Let $f:X \lrg Y$ be a  holomorphic map of compact complex
manifolds. Let $\Gamma\subset X \times Y$ be the graph of $f$ and
$p:X \times Y\lrg X$ and $q: X \times Y \lrg Y$ be the natural
projections.

\begin{lem}\label{lemma kaler aciclico  q^* A_Y}
If X and Y are compact $K\ddot{a}hler$, then the sub-complexes
$Im(\de)=\de A_{X \times Y}$, $\de A_\Gamma $, $\de A_{X \times Y}
\cap q^* A_Y$ and  $\de A_{X \times Y} \cap p^* A_X$ are acyclic.

\end{lem}
\begin{proof}
By hypothesis, $X \times Y$ is  K\"{a}hler. Then applying the $\de
\debar$-Lemma~\ref{de debar lemma} to $A_{X \times Y}$ we get
$$
\ker(\debar)\cap \image (\de)= \image(  \debar \de)
$$
and so $H^*_{\debar} (\de (A_{X \times Y}))=0$.  $\Gamma\subset X
\times Y$ is also K\"{a}hler  and so the same conclusion holds for
$\de A_\Gamma$: $\de A_\Gamma$ is acyclic.

Analogously, since $Y$ is K\"{a}hler  $\de A_Y$ and $q^* \de A_Y$
are acyclic. Therefore, to prove that $\de A_{X \times Y} \cap q^*
A_Y$ is acyclic it  suffices to prove that $\de A_{X \times Y}
\cap q^* A_Y  =q^* \de A_Y$.

The inclusion $\supseteq$ is obvious. Let $p \in q^* A_Y \cap \de
A_{X \times Y} $, then $p=q^* \phi=\de z$ with $\phi \in A_Y$ and
$z \in A_{X \times Y}$. The map $ q^* : H_{\de}(A_Y)\lrg
H_{\de}(A_{X \times Y})$ is injective. Therefore, $\de p= q^* \de
\phi= \de \de z=0$ and so $\phi $ is $\de$-closed ($ \phi \in
H_{\de}(A_Y)$); moreover, $q^*[\phi]=[\de z]=0$. Thus $\phi$ is
$\de $-exact, that is, $\phi=\de t$ with $t \in A_Y$. This implies
$p=q^* \de t \in q^* \de A_Y $.

The case $\de A_{X \times Y} \cap p^* A_X$ can be proved in the
same way.
\end{proof}

\begin{oss}\label{oss de-debar lemma per X,Y,XxY, Gamma}
In the previous lemma the  K\"{a}hler hypothesis on $X$ and $Y$
can be substituted by the validity of the $\de \debar$-lemma  in
$A_X$,$A_Y$, $A_{X \times Y}$ and $A_\Gamma$.
\end{oss}

\section{Holomorphic fiber bundle and Dolbeault's
cohomology}\label{sezio diff dolbeault su fibrati holo}

Let $E$ be a  holomorphic fiber bundle on $X$. Then the $\debar$
operator can be extended to the $Dolbeault$ operator
$$
\debar_E: \shA_X^{p,q}(E) \lrg \shA_X^{p,q+1}(E).
$$
If $e_1, \ldots, e_n$ is a local frame for $E$ then
$$
\debar_E(\sum_i \phi_i e_i)= \sum _i \debar (\phi) e_i.
$$
Since $E$ is a  holomorphic fiber bundle, this definition does not
depend on the choice of the local frame. By definition, $\debar_E$
satisfies the property ${\debar_E }^2=0$.

\smallskip

Let $A_X^{p,q}(E)=\Gamma(X,\shA_X^{p,q}(E))$ be the vector space
of global sections of the sheaf $\shA_X^{p,q}(E)$. Then we can
consider,  for each $p\geq 0$, the following complex:
$$
0 \lrg A_X^{p,0}(E)\stackrel{\debar_E}{\lrg}  A_X^{p,1}(E)
\stackrel{\debar_E}{\lrg}  \cdots \stackrel{\debar_E}{\lrg}
A_X^{p,q}(E) \stackrel{\debar_E}{\lrg} \cdots.
$$

The cohomology of this complex is   \emph{Dolbeault's cohomology}
$H_{\debar_E}^{p,*}(X,E)$ of $E$. We note that, for $p=0$, $\ker
(\debar _E : A_X^{0,0}(E) \lrg A_X^{0,1}(E))$ coincides with the
holomorphic sections of $E$ and
$$
H^q_{\debar_E}(X,E):=H_{\debar_E}^{0,q}(X,E)=\frac{\ker (\debar _E
: A_X^{0,q}(E) \lrg A_X^{0,q+1}(E)) }{\image (\debar _E :
A_X^{0,q-1}(E) \lrg A_X^{0,q}(E))}.
$$


\begin{prop}\label{prop local closed allora loc esatta coeffic E}
Let $\alpha$ be a differential form  with coefficients in $E$ of
type $(0,q)$ with $q>0$. If $\debar_E \alpha=0$, then there
exists, locally on $X$, a differential form $\beta$ of type
$(0,q-1)$, with coefficients in $E$,  such that $\debar_E
\beta=\alpha$.
\end{prop}
\begin{proof}
See \cite[Proposition~2.36]{bib Voisin}.
\end{proof}


\subsection{\v{C}ech cohomology and Leray's theorem}
\label{sottosezione cech e leray theorem}

We follow \cite[Section~1.3]{bib manRENDICONTi}.

\medskip

Let $E$ be a  holomorphic bundle on the complex manifold $X$. Let
$\mathcal{U}=\{U_i\}_{i \in I}$ be a locally finite open covering
of $X$ and denote $U_{i_0 \cdots i_k}= U_{i_0} \cap \cdots \cap
U_{i_k}$.

Define the \v{C}ech   $q$-chains of $E$:
$$
\check{C}^k(\mathcal{U},E)=\{f_{i_0 \cdots i_k} \, | \, f_{i_0
\cdots i_k}:U_{i_0 \cdots i_k}\lrg E \mbox{ is a  holomorphic
section}\}
$$
and the \v{C}ech differential
$$
\check{\delta} : \check{C}^k(\mathcal{U},E) \lrg
\check{C}^{k+1}(\mathcal{U},E),
$$
$$
 (\check{\delta}{f})_{i_0 \cdots i_{k+1}}=\sum_{j=0}^{k+1}
(-1)^j\, f_{i_0 \cdots \widehat{i_j} \cdots i_{k+1}}.
$$
A simple calculation shows that $\check{\delta}^2=0$ and so we can
define
 the $\C$-vector space of \v{C}ech cohomology
$$
\check{H}^k(\mathcal{U},E)=\frac{\ker (\check{\delta} :
\check{C}^k(\mathcal{U},E) \lrg \check{C}^{k+1}(\mathcal{U},E) }
{\image (\check{\delta} : \check{C}^{k-1}(\mathcal{U},E) \lrg
\check{C}^{k}(\mathcal{U},E)}.
$$

\bigskip

Next, define a morphism $\theta: \check{H}^k(\mathcal{U},E) \lrg
H_{\debar_E}^{0,k}(X,E)$.

\smallskip

Let $t_i : X \lrg \C$, with $i \in I$, be a partition of unity
subordinate to the cover $\mathcal{U}$, that is, $supp(t_i)
\subset U_i$, $\sum_i t_i=1$ and $\sum_i \debar t_i=0$.

\smallskip

For each $f \in \check{C}^k(\mathcal{U},E)$ and $i \in I$ we
define
$$
\phi_i(f)=\sum_{j_1 \cdots j_k} f_{ij_1 \cdots j_k} \debar
t_{j_1}\wedge \cdots \wedge \debar t_{j_k} \in \Gamma(U_i,
\shA^{0,k}(E))
$$
and then
$$
\phi(f)= \sum_i t_i \phi_i(f) \in \Gamma(X, \shA^{0,k}(E)).
$$
It can be proved that $\phi$ is a well defined  morphism of
complexes that induces a morphism $\theta $ in cohomology (for
details see \cite[Proposition~1.22]{bib manRENDICONTi}).

\begin{teo}\label{teo VOISINh^q=H^check q}
Let $ \mathcal{U}=\{U_i\}_{i \in I}$ be a locally finite countable
open covering of a complex manifold $X$ and $E$ a  holomorphic
vector bundle. If $H_{\debar_E}^{k-q}(U_{i_0 \cdots i_q},E)=0$ for
every $q<k$ and $i_0 \cdots i_k$, then the morphism $\theta$ is an
isomorphism
$$
\theta:\check{H}^k(\mathcal{U},E) \lrg H_{\debar_E}^{k}(X,E).
$$
\end{teo}
\begin{proof}
See \cite[Theorem~1.24]{bib manRENDICONTi} or
\cite[Theorem~4.41]{bib Voisin}.
\end{proof}

\begin{oss}
If the open sets $U_i$ of the cover $\mathcal{U}$ are
biholomorphic to open convex subsets of $\C^n$ then $\mathcal{U}$
satisfies the hypothesis of Theorem~\ref{teo VOISINh^q=H^check q}.

\end{oss}


\begin{oss}\label{oss associo Hcheck iso H dolebeautl}
It is convenient to give an explicit description of the inverse
map  of $\theta$, at least for $k=2$:
$$
\theta^{-1}:  H_{\debar_E}^{2}(X,E)\lrg
\check{H}^2(\mathcal{U},E).
$$

Let $h \in H_{\debar_E}^{2}(X,E)$. By applying
Proposition~\ref{prop local closed allora loc esatta coeffic E}
for each $i \in I$, there exists $\tau_i \in \Gamma(U_i,
\shA^{0,1}(E))$ such that $h_{|U_i}=\debar \tau_i$.

Define $\sigma_{ij}=(\tau_i -\tau_j)_{|U_{ij}} \in \Gamma
(U_{ij},\shA^{0,1}(E))$. $\sigma_{ij}$ is $\debar$-closed;
actually,
$$
\debar\sigma_{ij}= ( \debar \tau_i - \debar
\tau_j)_{|U_{ij}}=h_{|U_{ij}}- h _{|U_{ji}}=0.
$$

Therefore,  for each $U_{ij}$ there exists $\rho_{ij}\in \Gamma
(U_{ij},\shA^{0,0}(E))$, such that $\debar \rho_{ij}=\sigma_{ij}$.

We observe that $(  \sigma_{jk} -  \sigma_{ik} +\sigma_{ij})
_{|U_{ijk}}=0$; indeed,
$$
( \sigma_{jk} -  \sigma_{ik} +\sigma_{ij})_{|U_{ijk}}=
$$
$$
( ( \tau_j - \tau_k )-( \tau_i -\tau_k) +( \tau_i - \tau_j )
)_{|U_{ijk}}=0.
$$

Define $\alpha_{ijk}= (\rho_{jk}-\rho_{ik}+ \rho_{ij})_{|U_{ijk}}
\in \Gamma (U_{ijk},\shA^{0,0}(E)) $. First of all, we have that
$\debar \alpha_{ijk} =0$; actually,
$$
\debar \alpha_{ijk} = (\debar \rho_{jk} -\debar  \rho_{ik}  +
\debar \rho_{ij})_{|U_{ijk}}= ( \sigma_{jk} -  \sigma_{ik}
+\sigma_{ij})_{|U_{ijk}}=0
$$
This implies that $\alpha_{ijk} \in \Gamma(U_{ijk},E) $.

Moreover, $(\check{\delta}\alpha)_{ijkl}=0$; indeed,
$$
(\check{\delta}\alpha)_{ijkl}= (\alpha_{jkl } -\alpha_{ikl
}+\alpha_{ijl }-\alpha_{ijk })_{|U_{ijkl}}=
$$
$$
(( \rho_{kl}-\rho_{jl}+ \rho_{jk} ) -(\rho_{kl}-\rho_{il}+
\rho_{ik} ) +(\rho_{jl}-\rho_{il}+ \rho_{ij} ) -(
\rho_{jk}-\rho_{ik}+ \rho_{ij}  ) )_{|U_{ijkl}}=0.
$$

This implies that  $\alpha \in \check{H}^2(X,E)$.

$\alpha$ is independent of the choices of the lifitngs $\tau_i$.
Actually, if we choose $\overline{\tau}_i$, such that
$h_{|U_i}=\debar \overline{\tau}_i$, then
$\overline{\tau}_i=\tau_i +\debar t_i$ and this change does not
affect the choice of $\alpha_{ijk}$.

If we choose $\overline{\rho}_{ij} \in \Gamma
(U_{ij},\shA^{0,0}(E))$ such that $\debar
\overline{\rho}_{ij}=\sigma_{ij}$, then $\overline{\rho}_{ij}=
\rho_{ij}+s_{ij}$, with $s_{ij} \in \Gamma (U_{ij},\shA^{0,0}(E))$
such that $\debar s_{ij}=0$. This implies that $s_{ij} \in \Gamma
(U_{ij},E)$. Therefore, $\{\overline{\alpha}_{ijk }\}=\{
\alpha_{ijk } \}+ \{  \check{\delta}( s_{ij})\}$ and so
$\overline{\alpha}_{ijk } $ and $\alpha_{ijk }$ represent the same
class in cohomology.

In conclusion, we have defined a map
$$
\vartheta:H_{\debar_E}^{2}(X,E)\lrg \check{H}^2(\mathcal{U},E),
$$
$$
[h]  \longmapsto [\alpha].
$$

Finally, it can be proved that   this map $\vartheta$ is the
inverse of $\theta$ (for   details see \cite[Theorem~1.24]{bib
manRENDICONTi} or \cite[Theorem~4.41]{bib Voisin}).
\end{oss}

\section{The Kodaira-Spencer  algebra $KS_X$}
 \label{sezio DGLA di kodaira-spencer}

\begin{dhef}\label{def kodaira spencer algebra}
Let  $\Theta_X$ be the holomorphic tangent bundle of a complex
manifold $X$. The  \emph{Kodaira-Spencer} (differential graded
Lie) algebra of $X$ is
$$
KS_X=\bigoplus_i \Gamma(X,\shA_X^{0,i}(\Theta_X))=\bigoplus_i
A_X^{0,i}(\Theta_X).
$$
In particular, $KS_X^i$ is the vector space of the global sections
of the sheaf of germs of the differential $(0,i)$-forms with
coefficients in $\Theta_X$.

The differential $\tilde{d}$ on $KS_X$ is the opposite of
Dolbeault's differential, whereas the bracket is defined in local
coordinates as the $\overline{\Omega}^*$-bilinear extension of the
standard bracket on $\shA^{0,0}_X(\Theta_X)$
($\overline{\Omega}^*=\ker (\de : \shA_X^{0,*} \lrg \shA_X^{1,*})$
is the sheaf of antiholomorphic differential forms).

\end{dhef}
Explicitly, if $z_1, \ldots,z_n$ are local holomorphic coordinates
on $X$, we have
$$
\tilde {d}(f d\overline{z}_I \frac{\de}{\de z_i})=-\debar
(f)\wedge d\overline{z}_I \frac{\de}{\de z_i}.
$$
$$
[f\frac{\de}{\de z_i}\, d\overline{z}_I,g\frac{\de}{\de z_j} \,
d\overline{z}_J]=(f\frac{\de g}{\de z_i}\frac{\de}{\de z_j}-
g\frac{\de f}{\de z_j}\frac{\de}{\de z_i})\, d\overline{z}_I
\wedge d \overline{z}_J ,\qquad \forall  \ f,g \in \shA_X^{0,0},
$$
($\shA_X^{0,*}(T_X)$ is a sheaf of DGLAs).

\noindent We note that by Dolbeault theorem we have
$H^i(A^{0,*}_X(\Theta_X))\cong H^i(X,\Theta_X)$ for all $i$, then

$$
H^q(KS_X)=\frac{\ker (\debar   : A_X^{0,q}(\Theta_X) \lrg
A_X^{0,q+1}(\Theta_X)) }{\image (\debar  : A_X^{0,q-1}(\Theta_X)
\lrg A_X^{0,q}(\Theta_X))} \cong H^q_{\debar }(X,\Theta_X).
$$

In Theorem~\ref{teo def_k =Def_X}, we will prove that the DGLA
$KS_X$ controls the infinitesimal deformations of  $X$.

\section{Contraction map and  holomorphic Lie derivative}
\label{sec cotrction map e Lie deriv}

In general, for any   vector space $V$  and linear functional
$\alpha:V \lrg \C$, we can define the \emph{contraction operator}
$$
\alpha \contr: \bigwedge^k V \lrg \bigwedge^{k-1}V,
$$
$$
\alpha\contr (v_1 \wedge \ldots \wedge v_k)=\sum_{i=1}^k
(-1)^{i-1} \alpha(v_i)(v_1 \wedge \ldots \wedge \hat{v}_i \wedge
\ldots \wedge v_k),
$$
that is a derivation of degree $-1$ of the graded algebra
$(\bigwedge^k V ,\wedge)$.

Then considering the contraction  $\contr$ of the differential
forms with vector fields  we can define two injective morphisms of
sheaves:
\begin{itemize}
  \item[-] the \emph{contraction map}
$$
\bi : \shA^{0,*}_X(\Theta_X)\lrg Der^*(\shA_X^{*,*})[-1],
$$
$$
  a \longmapsto \bi_a \ \mbox{ with } \
\bi_a(\omega)=a \contr \omega
$$

\item[-] the \emph{holomorphic Lie derivative}
$$
\bl: \shA^{0,*}_X(\Theta_X)\lrg Der^*(\shA_X^{*,*}),
$$
$$
   a \longmapsto \bl_a=[\de,\bi_a]\ \mbox{ with } \
\bl_a(\omega)=\de(a \contr \omega)+ (-1)^{\deg(a)} a \contr \de
\omega
$$
\end{itemize}
for each $a \in \shA^{0,*}_X(\Theta_X)$ and $ \omega \in
\shA_X^{*,*}$.

\begin{lem}\label{lem formule cartan}
With the notation above, for every $a,b \in
\shA_X^{0,*}(\Theta_X)$ we have
$$
\bi_{\tilde{d}a}=-[\debar, \bi_a],\qquad \bi_{[a,b]}=
[\bi_a,[\de,\bi_b]]=[[\bi_a,\de],\bi_b], \qquad [\bi_a,\bi_b]=0.
$$
\end{lem}
\begin{proof}
See  \cite[Lemma~2.1]{bib mane COSTRAINT}. Let  $z_1,z_2, \ldots ,
z_n$ be local holomorphic coordinates on $X$. By linearity, we can
assume that $\displaystyle a= f d\overline{z}_I \frac{\de}{\de
z_i}$ and $\displaystyle  b= g d\overline{z}_J\frac{\de}{\de z_j}$
($i \neq j$), with $f,g \in \shA_X^{0,0}$.

All the expressions vanish on $\shA_X^{0,*}$ and $\shA_X^{*,*}$ is
generated as $\C$-algebra by $\shA_X^{0,0}\oplus \shA_X^{0,1}
\oplus \shA_X^{1,0}$. Therefore, it is sufficient  to verify the
equalities on the $dz_h$ (that generate $\shA_X^{1,0}$).

Moreover, we note that $\debar dz_h=\de dz_h=\bi_a \bi_b
dz_h=\bi_b \bi_a  dz_h=0$. Therefore, $ [\bi_a,\bi_b]=0$ and the
other equalities follow from the easy calculations below.

Let $\omega=dz_h$ and $\displaystyle \tilde {d}(a)=-\debar
(f)\wedge d\overline{z}_I \frac{\de}{\de z_i}$. Then
$$
\bi_{\tilde{d}a}(dz_h) =\tilde{d}a \contr dz_h =
  \begin{cases}
     0 & \text{$h\neq i$}, \\
   -\debar (f)\wedge
d\overline{z}_I  & \text{$h=i$}.
  \end{cases}
$$
On the other hand,
$$
-[\debar, \bi_a](\omega)=(- \debar
\bi_a+(-1)^{\overline{a}-1}\bi_a \debar)(dz_h)=- \debar
\bi_a(dz_h)=
$$
$$
- \debar ( a \contr dz_h)=
  \begin{cases}
     0 & \text{$h\neq i$}, \\
   -\debar (f d\overline{z}_I) & \text{$h=i$}.
  \end{cases}
$$
Then the first equality holds.

As to $\bi_{[a,b]}$, we have
$$
[a,b]=(f\frac{\de g}{\de z_i}\frac{\de}{\de z_j}- g\frac{\de
f}{\de z_j}\frac{\de}{\de z_i})\, d\overline{z}_I \wedge d
\overline{z}_J
$$
and then
$$
\bi_{[a,b]}(dz_h) =
  \begin{cases}
0 & \text{$h\neq i,j$}, \\
\displaystyle f\frac{\de g}{\de z_i}\, d\overline{z}_I \wedge
  d\overline{z}_J   & \text{$h=j$},\\
\displaystyle -\,g\, \frac{\de f}{\de z_j}\, d\overline{z}_I
\wedge d\overline{z}_J   &
  \text{$h=i$}.
  \end{cases}
$$
On the other hand,
$$
 [\bi_a,[\de,\bi_b]]=[\bi_a,\de \bi_b
-(-1)^{\overline{b}-1}\bi_b \de]=
$$
$$
 \bi_a\de \bi_b
-(-1)^{\overline{b}-1} \bi_a\bi_b \de
-(-1)^{(\overline{a}-1)\overline{b}}(\de \bi_b \bi_a
-(-1)^{\overline{b}-1}\bi_b \de \bi_a ).
$$
Then
$$
[\bi_a,[\de,\bi_b]](dz_h)= \bi_a\de \bi_b(dz_h)-(-1)^{
\overline{a} \overline{b}} \bi_b \de \bi_a (dz_h)=
$$
$$
  \begin{cases}
0 & \text{$h\neq i,j$}, \\
\displaystyle \bi_a\de \bi_b(dz_j)=f\frac{\de g}{\de z_i}
d\overline{z}_I \wedge
  d\overline{z}_J   & \text{$h=j$},\\
\displaystyle -(-1)^{ \overline{a} \overline{b}} \bi_b \de \bi_a
(dz_i)= -(-1)^{ \overline{a} \overline{b}}g\frac{\de f}{\de z_j}
d\overline{z}_J \wedge d\overline{z}_I &  \text{$h=i$}.
  \end{cases}
$$
\end{proof}
The previous set of equalities is referred to as $Cartan$
$formulas$.

\begin{dhef}
Let $L$ and $M$ be two differential graded Lie algebras and let
$d'$ be the differential on the graded vector space $\Hom^*(L,M)$.
A linear map $i \in \Hom^{-1}(L,M)$ is called a $Cartan$
$homotopy$ if
$$
i([a,b])=[i(a),d'i(b)] \qquad \mbox{ and } \qquad[i(a),i(b)]=0
\qquad \forall \, a,b \in L.
$$
\end{dhef}
We recall that, by definition (see Example~\ref{exe dhef
Hom^*(V,W) grade vect spac}), we have
$$
d'i(a)=d_M(i(a))+i(d_L (a)).
$$
\begin{cor}
$\bi$ is a Cartan homotopy and the Lie derivative $\bl$ is a
morphism of sheaves of DGLAs.
\end{cor}
\begin{proof}
Using Cartan formulas we get
$d'(\bi_b)=[d,\bi_b]+\bi_{\tilde{d}b}=[\de +\debar,
\bi_b]-[\debar,\bi_b]=[\de,\bi_b]$. Then $ \bi_{[a,b]}=
[\bi_a,[\de,\bi_b]]=[\bi_a,d'(\bi_b)]$. Moreover, by
Lemma~\ref{lem formule cartan}, $[\bi_a,\bi_b]=0$ and so $\bi$ is
a Cartan homotopy.

As to $\bl$, we have
$$
\bl_{\tilde{d}a}=[\de,\bi_{\tilde{d}a}]=-[\de,[\debar, \bi_a]].
$$
Moreover,
$$
-[\de,[\debar, \bi_a]]=-[\de,\debar \bi_a -(-1)^{\deg (\bi_a)}
\bi_a \debar]=
$$
$$
 -\de \debar \bi_a+(-1)^{\deg (\bi_a)} \de\bi_a
\debar -(-1)^{\deg(\bi_a)} \debar \bi_a \de+ \bi_a \debar \de.
$$
Therefore,
$$ [d,\bl_a]=[\de + \debar
,\bl_a]=[\de,\bl_a]+[\debar,\bl_a]=
$$
$$
=[\de,[\de,\bi_a]]+[\debar,[\de,\bi_a]]=-[\de,[\debar,
\bi_a]]=\bl_{\tilde{d}a}.
$$
As to the equality $\bl_{[a,b]}=[\bl_a,\bl_b]$, we can prove it as
follows
$$
\bl_{[a,b]}= [\de,\bi_{[a,b]}]= [\de,[\bi_a,[\de,\bi_b]]]=
$$
$$
\de [\bi_a,[\de,\bi_b]]-(-1)^{\overline{a}-1+\overline{b}}
[\bi_a,[\de,\bi_b]]\de=
$$
$$
\de(\bi_a [\de,\bi_b] -  (-1)^{(\overline{a}-1)\overline{b}}
[\de,\bi_b]\bi_a) -(-1)^{\overline{a}-1+\overline{b}} ( \bi_a
[\de,\bi_b]  -(-1)^{(\overline{a}-1)\overline{b}}
[\de,\bi_b]\bi_a) \de=
$$
$$
\de \bi_a [\de,\bi_b] -  (-1)^{ \overline{a} \overline{b}}
[\de,\bi_b]\de\bi_a -(-1)^{\overline{a}-1} \bi_a \de [\de,\bi_b] -
(-1)^{ \overline{a} +\overline{a}\overline{b}}[\de,\bi_b] \bi_a
\de=
$$
$$
(\de \bi_a -(-1)^{\overline{a}-1}\bi_a \de)  [\de,\bi_b]
-(-1)^{\overline{a} \overline{b}} [\de,\bi_b] (\de \bi_a
+(-1)^{\overline{a} }\bi_a \de)=
$$
$$
[\de,\bi_a][\de,\bi_b]-(-1)^{\overline{a} \overline{b}}
[\de,\bi_b][\de,\bi_a]=[[\de,\bi_a],[\de,\bi_b]]=[\bl_a,\bl_b].
$$

\end{proof}

In particular, we have an injective morphism of sheaves
$$
\bl : \shA_X^{0,*}(\Theta_X) \lrg Der^*(\shA_X^{0,*},\shA_X^{0,*}
)
$$
\begin{equation}\label{equa def l_a}
a \longmapsto \bl_a \ \mbox{ with }\ \bl_a(\omega)=(-1)^{\deg(a)}
a \contr \de \omega.
\end{equation}

Explicitly, in local holomorphic coordinates $z_1,z_2, \ldots ,
z_n$, if $\displaystyle a=g\  d\overline{z}_I \frac{\de}{\de z_i}$
and $\omega= f d\overline{z}_J$, then
$$
\bl_a(\omega)= (-1)^{\deg (a)} \frac{\de f}{\de z_i} \, g \,
d\overline{z}_I \wedge d\overline{z}_J.
$$

Using $\bl$, for each $(A,m_A) \in \Art$, we can define the
following morphism:
$$
\bl : \shA_X^{0,*}(\Theta_X)\otimes A  \lrg
Der^*(\shA_X^{0,*}\otimes A,\shA_X^{0,*} \otimes A).
$$
In particular, for each solution of the Maurer-Cartan equation in
$KS_X$ we have the fundamental lemma below.

\begin{lem}\label{lem x sol MC sse (debar+l_x )^2=0}
$x \in \MC_{KS_X}(A)$ if and only if
$$
\debar + \bl_x: \shA_X^{0,*}\otimes A \lrg \shA_X^{0,*+1}\otimes A
$$
is a differential of degree 1 on $\shA_X^{0,*}\otimes A$.
\end{lem}
\begin{proof}
Since $\bl$ is a morphism of DGLAs, we have
$$
(\debar+\bl_x)^2=\debar \bl_x+\bl_x \debar +
\bl_x^2=[\debar,\bl_x]+ \displaystyle\frac{1}{2} [\bl_x,\bl_x]=\bl
(\tilde{d} x+ \displaystyle\frac{1}{2}[x,x]).
$$
\end{proof}

Moreover, using $\bl$, we can also define, for any $(A,m_A) \in
\Art$ and  $a \in \shA_X^{0,0}(\Theta_X) \otimes m_A$, an
automorphism $e^a$ of $\shA _X^{0,*} \otimes A$:
\begin{equation}\label{equa def di e^a su shA^(0,*)tensor A}
e^a : \shA _X^{0,*} \otimes A \lrg \shA _X^{0,*} \otimes A, \ \ \
\ \ e^a(f)=\sum^\infty_{n=0}\frac{\bl_a^n}{n!}(f).
\end{equation}

\begin{lem}\label{lemma e^a(d+l_x)e^-a = d+e^a *x}
For every local Artinian $\C$-algebra $(A,m_A)$,   $a\in
\shA_X^{0,0}(\Theta_X)\otimes m_A$ and  $x
 \in \MC_{KS_X}(A)$  we have
\begin{equation}\label{equa COMM e^a o d +x o e^-a=e^a*x}
e^a\circ (\debar +\bl_x) \circ e^{-a}=\debar + e^a*\bl_x\ :\
\shA_X^{0,0} \otimes A \lrg \shA_X^{0,1} \otimes A,
\end{equation}
where $*$ is the gauge action (and $e^a*\bl_x \in
\shA_X^{0,1}(\Theta_X)\otimes m_A$ acts on $\shA_X^{0,0}\otimes A$
as defined in (\ref{equa def l_a})). In particular,
$$
\ker(\debar + e^a * \bl_x: \shA_X^{0,0}\otimes A \lrg
\shA_X^{0,1}\otimes A)=e^a(\ker(\debar +  \bl_x:
\shA_X^{0,0}\otimes A \lrg \shA_X^{0,1}\otimes A)).
$$

\end{lem}
\begin{proof}
It follows from the definition of gauge action. More precisely,
since $e^a \circ e^b \circ e^{-a}=e^{[ a, \, ]}(b)$, we have
$$
e^a\circ (\debar +l_x ) \circ e^{-a}=e^{[ a, \, ]'}(\debar +l_x)=
$$
$$
\sum^\infty_{n=0}\frac{ {([ a, \, ]')} ^n}{n!}(\debar +l_x)=\debar
+ l_x + \sum^\infty_{n=1}\frac{([ a, \, ]')^n}{n!}(\debar+l_x)=
$$
$$
=\debar + l_x + \sum^\infty_{n=0}\frac{{([ a, \,
]')}^{n+1}}{(n+1)!}(\debar+l_x)=\debar + l_x +
\sum^\infty_{n=0}\frac{  {([ a, \, ]')}
^n}{(n+1)!}([a,\debar]'+[a,l_x])=
$$
$$
\debar + l_x + \sum^\infty_{n=0}\frac{ {([ a, \, ]')}
^n}{(n+1)!}([a,l_x]-\debar a)=\debar + e^a *l_x.
$$
\end{proof}


\begin{oss}
Let  $\phi_i$ be an automorphism of the $A$-module $\shA _X^{0,*}
\otimes A$ whose specialization to the residue field $\C$ is the
identity. Let $\phi=\sum_i \phi_i=id+\eta$ with $\eta \in
\Hom^0(\shA _X^{0,*},\shA _X^{0,*}) \otimes m_A$. Since we are in
characteristic zero, we can take the logarithm so that $\phi=e^a$
with $a \in   \Hom^0(\shA _X^{0,*},\shA _X^{0,*}) \otimes m_A$

\end{oss}

\subsection{The DGLA of a submanifold}
\label{sezio submanifold defz L'}

Let $X$ be a complex manifold and $i:X \hookrightarrow Y$ be the
inclusion of a submanifold $X$. Let $i^*:\shA_X^{0,*}\lrg
\shA_Y^{0,*}$ be the restriction morphism (of sheaves of DGLAs).
Finally, denote by $\Theta_Y$ the holomorphic tangent bundle of
$Y$ and by $N_{X|Y}$ the normal bundle of $X$ in $Y$. Define the
sheaf  $\mathcal{L'}=\oplus_i \mathcal{L'}^i$ such that
$$
0 \lrg \mathcal{L'}  \lrg \shA_Y^{0,*}(\Theta_Y) \lrg \shA_X^{0,*}
(N_{X|Y})  \lrg 0.
$$

Let $z_1,\ldots,z_n$ be holomorphic coordinates on $Y$ such that
$Y\supset X=\{z_{t+1}=\cdots=z_n=0\}$. Then $\eta \in
\mathcal{L'}^i $ if and only  if $\displaystyle \eta=\sum_{j=1}^n
\omega_j \frac{\de}{\de z_j}$, with $\omega_j \in \shA_Y^{0,i}$
such that $\omega_j \in\ker i^*$ for $j \geq t$. In particular,
$\mathcal{L'}^0$ is the sheaf of differentiable vector fields on
$Y$ that are tangent to $X$.

\begin{lem}\label{lemm L manetti is DGLA}
$\mathcal{L'}$ is a sheaf of differential graded Lie subalgebras
of $\shA_Y^{0,*}(\Theta_Y)$ such that $\bl_a (\ker i^*)\subset
\ker i^*$ if and only if $a \in \mathcal{L'}\subset
\shA_Y^{0,*}(\Theta_Y)$.
\end{lem}
\begin{proof}
See \cite[Section~5]{bib ManettiPREPRINT}. It is an easy
calculation in local holomorphic coordinates.

\end{proof}

Moreover,  consider the  automorphism  $e^a$ of $\shA _X^{0,*}
\otimes A$ defined in (\ref{equa def di e^a su shA^(0,*)tensor
A}): if $a \in  \mathcal{L'}^0\otimes m_A$ then $e^a(\ker
(i^*)\otimes A)=\ker (i^*) \otimes A$.

Let $L'$ be the DGLA of the global sections of $\mathcal{L'}$:
$$
0 \lrg L' \lrg A_Y^{0,*}(\Theta_Y) \stackrel{\pi'} \lrg
A_X^{0,*}(N_{X|Y})\lrg 0.
$$
In the literature, the notation $L'=A_Y^{0,*}(\Theta_Y(- log\ X))$
can also be found.

In Section~\ref{sezio inclusion}   we will prove that $L'$
controls the embedded deformations of the inclusion $i:X
\hookrightarrow Y$ (Corollary~\ref{cor L governa def coppia
(XcY)}).

\section{Induced maps $f_* $ and $f^*$ by a holomorphic map
$f$}\label{sezione f_* and f^*}

This section is devoted to the study  of the maps $f_*$ and $f^*$
induced by a  holomorphic map $f$. In particular, we prove a
property of these maps (Lemma~\ref{lemma f^*(x-|w)=f_*x-|f^*w})
that will be used in the last chapter (Section~\ref{sottosection
semiregolarity}).

Let $f:X \lrg Y$ be a  holomorphic map of compact complex
manifolds.

 Let $\mathcal{U}=\{U_i\}_{i \in I}$ and
$\mathcal{V}=\{V_i\}_{i \in I}$ be finite Stein open  covers of
$X$ and $Y$, respectively, such that $f(U_i) \subset V_i$ ($U_i$
is allowed to be empty). Then $f$ induces morphisms
$$
f^*: \check{C}^p(\mathcal{V},\Theta_Y) \lrg
\check{C}^p(\mathcal{U},f^*\Theta_Y),
$$
and
$$
f_*: \check{C}^p(\mathcal{U},\Theta_X) \lrg
\check{C}^p(\mathcal{U},f^*\Theta_Y).
$$
Explicitly, for each $i \in I$ and local holomorphic coordinate
systems $z=(z_1,z_2,\ldots ,z_n)$ on $U_i$ and $w=(w_1,w_2, \ldots
w_m)$ on $V_i$ such that $f(z_1,z_2,\ldots ,z_n ) =(f_1(z),\ldots
, f_m(z))$, we have
$$
f^*:\Gamma(V_i,\Theta_Y) \lrg \Gamma(U_i,f^*\Theta_Y),
$$
$$
f^*(\sum_j g_j(w)  \frac{\de}{\de w_j})=\sum_j  g_j(f(z))
\frac{\de}{\de w_j}
$$
and
$$
f_*:\Gamma(U_i,\Theta_X) \lrg \Gamma(U_i,f^*\Theta_Y),
$$
$$
f_*(\sum_k h_k(z) \frac{\de}{\de z_k})=\sum_{k,j}h_k(z) \frac{\de
f_j(z)}{\de z_k} \frac{\de}{\de w_j}.
$$

Moreover, $f_*$ and $f^*$ commute with the \v{C}ech differential
and they do not  depend on the choice of the cover. Therefore, we
get linear maps in cohomology :
$$
f^*: \check{H}^p(Y,\Theta_Y) \lrg \check{H}^p(X,f^*\Theta_Y)
$$
and
$$
f_* : \check{H}^p(X,\Theta_X) \lrg \check{H}^p(X,f^*\Theta_Y).
$$

\bigskip

Analogously, $f$ induces morphisms
$$
f^*: A_Y^{p,q}(\Theta_Y) \lrg A_X^{p,q}(f^*\Theta_Y)
$$
and
$$
f_*: A_X^{p,q}(\Theta_X) \lrg A_X^{p,q}(f^*\Theta_Y).
$$
Let $\mathcal{U}$ and  $\mathcal{V} $ be  Stein covers and
$z=(z_1,z_2,\ldots ,z_n)$ on $U_i$ and $w=(w_1,w_2, \ldots w_m)$
on $V_i$  local holomorphic coordinate systems  as above. Let
$K=(1 \leq k_1<k_2 < \cdots < k_p \leq n)$ be a multi-index of
length $p$ and $J=(1 \leq j_1<j_2 < \cdots < j_q \leq n)$ a
multi-index of length $q$. Then
$$
f_*: A ^{p,q}(U_i,\Theta_X) \lrg A_X^{p,q}(U_i,f^*\Theta_Y)
$$
$$
f_*\left(  h (z )dz_K \wedge d\overline{z}_J \frac{\de}{\de
z_i}\right)= h (z)  dz_K \wedge d\overline{z}_J \sum_j\frac{\de
f_j(z)}{\de z_i} \frac{\de}{\de w_j},
$$
and
$$
f^*: A ^{p,q}(V_i,\Theta_Y) \lrg A_X^{p,q}(U_i,f^*\Theta_Y),
$$
$$
f^*\left(  g (w )  dw_K \wedge d\overline{w}_J \frac{\de}{\de
w_i}\right) =   g(f (z) )\de f_K \wedge \debar f_J \frac{\de}{\de
w_i},
$$
where
$$
\de f_K=\sum_{h=1}^n \frac{\de f_{k_1}}{\de z_h} dz_h\wedge \cdots
\wedge \sum_{h=1}^n  \frac{\de f_{k_p}}{\de z_h }dz_h
$$
and
$$
\debar f_J=\sum_{h=1}^n \frac{\de f_{j_1}}{\de \overline{z}_h}
d\overline{z}_h\wedge \cdots \wedge \sum_{h=1}^n  \frac{\de
f_{j_q}}{\de\overline{z}_h }d\overline{z}_h.
$$

We note that $f^*$ and $f_*$ commute with $\de$ and $\debar$.

Moreover, for each $k$, there exist the following commutative
diagrams
\begin{center}
$\xymatrix{\check{C}^q(\mathcal{V},\Theta_Y)  \ar[r]^{\phi}
\ar[d]^{f^*} & A_Y^{0,q}(\Theta_Y)\ar[d]^{f^*}  \\
\check{C}^q(\mathcal{U},f^*\Theta_Y) \ar[r]^{\phi} &
A_X^{0,q}(f^*\Theta_Y)   \\ }$
\end{center}
and
\begin{center}
$\xymatrix{\check{C}^q(\mathcal{U},\Theta_X)  \ar[r]^{\phi}
\ar[d]^{f_*} & A_X^{0,q} (\Theta_X)\ar[d]^{f_*}  \\
\check{C}^q(\mathcal{U},f^*\Theta_Y) \ar[r]^{\phi} &
 A_X^{0,q}(f^*\Theta_Y),   \\ }$
\end{center}
where $\phi $ is the map defined in Section~\ref{sottosezione cech
e leray theorem}. Therefore,  $f_*\phi =\phi f_*$ and $f^*\phi
=\phi f^*$.

\begin{lem}\label{lemma f^*(x-|w)=f_*x-|f^*w}
Let $f:X \lrg Y$ be a  holomorphic map of complex manifolds. Let
$\chi \in \shA_Y^{0,*}(\Theta_Y)$  and $\eta \in
\shA_X^{0,*}(\Theta_X)$ such that $f^*\chi =f_* \eta \in
\shA_X^{0,*}(f^*\Theta_Y)$. Then for any  $\omega \in
\shA_Y^{*,*}$
$$
f^*(\chi \contr \omega)= \eta \contr f^* \omega.
$$
\end{lem}

\begin{proof}
Let $\mathcal{U}=\{U_i\}_{i \in I}$ and $\mathcal{V}=\{V_i\}_{i
\in I}$ be finite open Stein covers of $X$ and $Y$, respectively,
as above. For each $i \in I$, let $z$ be local holomorphic
coordinate  systems on $U_i$ and $w$ on $V_i$ such that $f(z)
=(f_1(z),\ldots , f_m(z))$.

Let
$$
\shA_X^{0,r}(\Theta_X)\ni \eta=\sum_{i=1}^n h_i(z) d\overline{z}_I
\frac{\de  }{\de z _i }
$$
and
$$
\shA_Y^{0,r}(\Theta_Y)\ni \chi=\sum_{h=1}^m \varphi_h(w)
d\overline{w}_{H} \frac{\de  }{\de w_h },
$$
with $I=(1 \leq i_1<i_2 < \cdots < i_r \leq n)$   and $H=(1 \leq
h_1<h_2 < \cdots < h_r \leq n)$  multi-indexes  of length $r$.

Therefore,
$$ f_*\eta=\sum_{i=1}^n h_i(z) d\overline{z}_I
\sum_{h=1}^m \frac{\de f_{h}}{\de z_i} \frac{\de }{\de w_h
}=\sum_{h=1}^m \left( \sum_{i=1}^n h_i(z)  \frac{\de f_{h}}{\de
z_i} d\overline{z}_I \right)\frac{\de }{\de w_h }
$$
and
$$
f^*\chi=\sum_{h=1}^m \varphi_h(f(z)) \debar f_{H} \frac{\de }{\de
w_h }.
$$
By hypothesis, $ f^*\chi =f_* \eta \in \shA_X^{0,*}(f^*\Theta_Y)$,
then
\begin{equation}\label{equa f^*=f_* esplicita}
\varphi_h(f(z)) \debar f_{H}= \sum_{i=1}^n h_i(z)  \frac{\de
f_{h}}{\de z_i} d\overline{z}_I, \qquad \forall \ h=1, \ldots, m.
\end{equation}
Next, let
$$
\shA^{p,q}(V_i) \ni \omega =g(w )\, dw_K \wedge d\overline{w}_J,
$$
with $K=(1 \leq k_1<k_2 < \cdots < k_p \leq n)$ a multi-index of
length $p$ and $J=(1 \leq j_1<j_2 < \cdots < j_q \leq n)$ a
multi-index of length $q$. Then
$$
f^*\omega= g(f(z))\,  \de f_K \wedge \debar f_J
$$
and
$$
\chi \contr \omega= \sum_{h=1}^m \varphi_h(w) g(w) \,
d\overline{w}_{H} \wedge \left(  \frac{\de  }{\de w_h }\contr dw_K
\right)  \wedge d\overline{w}_J=
$$
$$
\sum_{h=1}^p (-1)^{h-1}\varphi_{k_h}(w) g(w)\,  d\overline{w}_{H}
\wedge    dw_{K-\{k_h\}} \wedge d\overline{w}_J,
$$
with $dw_{K-\{k_h\}} =dw_{k_1}\wedge \ldots \wedge \widehat{
dw_{k_h}} \wedge   \ldots \wedge dw_{k_p}$.

Therefore,
$$ f^*(\chi \contr \omega)= \sum_{h=1}^p
(-1)^{h-1}\varphi_{k_h}(f(z)) g(f(z))\,  \debar f_{H} \wedge \de f
_{K-\{k_h\}} \wedge \debar f_J
$$
and using (\ref{equa f^*=f_* esplicita}) we get
$$
f^*(\chi \contr \omega)= \sum_{h=1}^p (-1)^{h-1}  g(f(z)) \left(
\sum_{i=1}^n h_i(z)  \frac{\de f_{{k_h}}}{\de z_i} d\overline{z}_I
\right) \wedge \de f _{K-\{k_h\}} \wedge \debar f_J.
$$
On the other hand,
$$
\eta \contr f^* \omega =
$$
$$
\sum_{i=1}^n h_i(z) g(f(z)) \,  d\overline{z}_I \left(\frac{\de
}{\de z _i }\contr \de f_K \right)\wedge \debar f_J=
$$
$$
\sum_{i=1}^n h_i(z) g(f(z)) \,  d\overline{z}_I \wedge\left(
\sum_{h=1}^p (-1)^{h-1} \frac{\de f_{{k_h}}}{\de z_i} \de f
_{K-\{k_h\}} \right)\wedge \debar f_J=
$$
$$
f^*(\chi \contr \omega).
$$

\end{proof}

\section{Deformations of complex manifolds}
\label{sezio defo manifold defK=defX}

In this section we prove that the infinitesimal deformations of a
compact complex  manifold $X$ are controlled by the differential
graded Lie algebra of Kodaira-Spencer $KS_X$, i.e.,
$\Def_{KS_X}\cong \Def_X$.

We start with some lemmas and we postpone the proof  to
Section~\ref{sez dimo teo DEF_K=DEF_X}, where we also give an
explicit description of the isomorphism (see Theorem~\ref{teo
def_k =Def_X}).

\begin{lem}\label{lem x in MC local gauge 0}
Let $A \in \Art$ and  $ x \in MC_{KS_X}(A)$, then there exists a
cover $ \mathcal{U}=\{U_i\}$ of $X$, such that  $x_{|U_i} \sim 0$
for each $i$.
\end{lem}

\begin{proof}
By Proposition~\ref{prop local closed allora loc esatta coeffic
E}, there exists a cover $ \mathcal{U}=\{U_i\}$ such that
$H^1(U_i, \Theta_X)=0$ for all $i$. Moreover, by Remark~\ref{oss
tangente DEF_L =H^1}, $H^1(X, \Theta_X)$ is the tangent space of
the deformation functor  $\Def_{KS_X}$. Therefore, by
Corollary~\ref{cor funtoreF banale sse t_F=0}, $\Def_{KS_X}$ is
locally trivial and so any $ x \in MC_{KS_X}(A)$ is locally gauge
equivalent to zero.

\end{proof}

In Section~\ref{sezio DGLA di kodaira-spencer} we  defined a
morphism of sheaves
$$
\bl : \shA_X^{0,*}(\Theta_X)\otimes A  \lrg
Der^*(\shA_X^{0,*}\otimes A,\shA_X^{0,*} \otimes A),
$$
$$
a \longmapsto \bl_a \ \mbox{ with }\ \bl_a(\omega)=(-1)^{\deg(a)}
a \contr \de \omega.
$$
Let $x \in MC_{KS_X}(A)$. Explicitly, in local holomorphic
coordinates $z_1,z_2, \ldots , z_n$, if $ \displaystyle
x=\sum_{i,j}x_{ij}d\bar{z}_i\frac{\de~}{\de z_j} $  and $\omega= f
d\overline{z}_J$, then
$$
\bl_{x}(f)=- \sum_{i,j}x_{ij}\frac{\de f}{\de z_j}d\bar{z}_i\wedge
d\overline{z}_J .
$$
We also proved  that for each $x \in MC_{KS_X}(A)$
$$
\debar + \bl_x: \shA_X^{0,*}\otimes A \lrg \shA_X^{0,*+1}\otimes A
$$
is a differential (Lemma~\ref{lem x sol MC sse (debar+l_x )^2=0}).

Define $\shO_A(x)$ as the kernel of $\debar + \bl_x:
\shA_X^{0,0}\otimes A \lrg \shA_X^{0, 1}\otimes A$. Then we have
the complex
$$
0 \lrg \shO_A(x) \lrg \shA^{0,0}_X\otimes A
\xrightarrow{\debar+\bl_x}\shA^{0,1}_X\otimes A
\xrightarrow{\debar+\bl_x} \cdots
 \xrightarrow{\debar+\bl_x}\shA^{0,n}_X\otimes  A \lrg  0.
$$
In Section~\ref{sezio DGLA di kodaira-spencer} we   also defined,
for each $s \in \shA_X^{0,0}(\Theta_X)\otimes m_A$, an
automorphism $e^s$ of $\shA^{0,*}_X\otimes A$.

\begin{lem}\label{lem sollevamento ISO da O_A a complexes}
Let $F,G: \Art \lrg \Set$ be the following  functors
$$
F(A):=\{\mbox{isomorphisms of complexes } e^s : (\shA_X^{0,*}
\otimes A ,\debar +\bl_x) \lrg (\shA_X^{0,*} \otimes A,\debar
+\bl_y)
$$
$$
with \ s \in \shA_X^{0,0}(\Theta_X)\otimes m_A \ that\ specialize\
to\  identity \}
$$
$$
G(A):=\{\mbox{isomorphisms of sheaves of $A$-module } \psi :
\shO_A(x) \lrg \shO_A(y) $$
$$
that\  specialize\  to\  identity \}.
$$
Then the restriction  morphism   $\phi:F \lrg G$ is surjective.
\end{lem}

\begin{proof}
We proceed by induction on $d=dim_\C A$.

If $A=\C$, then $G(\C)=\{\,\mbox {identity}\}$   and so it can be
lifted.

Assume that  $d \geq 2$  and let
$$
0 \lrg J \lrg B \stackrel{\alpha}{\lrg} A \lrg 0
$$
be a small extension; by  induction, each element in $G(A)$ can be
lifted to $F(A)$.\\
Let $\psi$ be an isomorphism between $\shO_B(x)$ and $\shO_B(y)$
($\psi \in G(B)$); we want to lift it to an isomorphism $e^s$.

$\alpha(x)$ and $\alpha(y)$ are in $MC_{KS_X}(A)$ and $\psi$
induces an isomorphism of sheaves of $A$-modules
$\overline{\psi}:\shO_A(\alpha(x)) \lrg \shO_A(\alpha(y)) $.
Therefore, by the induction hypothesis we can lift
$\overline{\psi}$ to an isomorphism of complexes $e^{\overline{s}}
$, i.e., ${e^{\overline{s}}}^{-1}\circ (\debar +\bl_{\alpha(x)})
\circ e^{\overline{s}}=\debar +\bl_{\alpha(y)}$ with $\overline{s}
\in \shA_X^{0,0}(\Theta_X)\otimes m_A$.

Then we can suppose that $\alpha(x)=\alpha(y) \in
 A^{0,1}(\Theta_X)\otimes m_A$ and  that $e^{\overline{s}}$ is
the identity.

This implies the existence of an element $p \in
\shA_X^{0,1}(\Theta_X)\otimes J$ such that  $x=y+p$. Since $x$ and
$y$ satisfy the Maurer-Cartan equation,  $\debar p=0$. Indeed
$$
0 =dx + \displaystyle\frac{1}{2}[x,x]=d(y+p) +
\displaystyle\frac{1}{2}[y+p,y+p]=dy+dp +
\displaystyle\frac{1}{2}[y,y]=dp.
$$
Therefore, by Proposition~\ref{prop local closed allora loc esatta
coeffic E}, there exists a Stein cover $\mathcal{U}=\{U_i\}_{i \in
I}$ of $X$ such that $p$ is locally $\debar$-exact, i.e., for each
$i\in I$ there exists $t_i \in \shA_X^{0,0}(\Theta_X)\otimes J$
such that $\debar t_i =p_{|U_i}$. Then
$$
y_{|U_i}=(x-p)_{|U_i}=x_{|U_i}-\debar t_i=e^{t_i}*x_{|U_i},
$$
where we use the fact that $J \cdot m_B=0$ as in Example~\ref{exe
azione gauge se []=0}.

In particular, by Lemma~\ref{lemma e^a(d+l_x)e^-a = d+e^a *x},
$e^{t_i}:   (\shA ^{0,*}(U_i) \otimes B ,\debar +\bl_x) \lrg (\shA
^{0,*} (U_i) \otimes B,\debar +\bl_y) $ is an isomorphism of
complexes, that lifts the isomorphism $e^{t_i}: \shO_B(x)(U_i)
\lrg \shO_B(y)(U_i)$. We note that $e^{t_i}$ restricts to the
identity on $\shO_A(x)(U_i)$.

On the other hand, by   Lemma~\ref{lem x in MC local gauge 0}, the
Maurer-Cartan element $x$ is locally gauge equivalent to zero.
Then for each $i\in I$ there exists $a_i \in \shA^{0,0}(U_i,
\Theta_X) \otimes m_B$ such that $e^{a_i}*x_{|U_i}=0$. As before,
Lemma~\ref{lemma e^a(d+l_x)e^-a = d+e^a *x} implies that $e^{a_i}:
(\shA ^{0,*}(U_i) \otimes B ,\debar +\bl_x) \lrg (\shA ^{0,*}
(U_i) \otimes B,\debar )  $ is an isomorphism of complexes, that
lifts the isomorphism $e^{a_i}: \shO_B(x)(U_i) \lrg
\shO_X(U_i)\otimes B$.

Next, consider the isomorphism
$$
\varphi_{|U_i}: \shO_X(U_i)\otimes B \lrg \shO_X(U_i)\otimes B
$$
defined as follows:
$$
\varphi_{|U_i}=e^{a_i} \circ e^{-t_i}\circ \psi_{|U_i}\circ
e^{-a_i}:
$$
$$
\shO_X(U_i)\otimes B \stackrel{e^{-a_i} }{\lrg} \shO_B(x)(U_i)
\stackrel{ \psi_{|U_i} }{\lrg}  \shO_B(y)(U_i) \stackrel{ e^{-t_i}
}{\lrg}  \shO_B(x)(U_i)  \stackrel{e^{a_i} }{\lrg}
\shO_X(U_i)\otimes B.
$$
Then, $\varphi_{|U_i}$ is an automorphism of  $\shO_X(U_i)\otimes
B$ that restricts to the identity on $\shO_X(U_i)\otimes A$.
Therefore, Lemma~\ref{lemm serneII.1.5 automor-derivazio} implies
the existence of $q_i \in \Gamma(U_i, \Theta_X)\otimes J$ such
that $ \varphi_{|U_i}=e^{q_i}$. In particular, $e^{a_i} \circ
e^{-t_i}\circ \psi_{|U_i}\circ e^{-a_i}=e^{q_i}$; by
Remark~\ref{oss succ esatta lunga automorfismi derivazione}, the
automorphism $e^{q_i}$ commutes with the other automorphisms and
so
$$
\psi_{|U_i}=e^{t_i+q_i}.
$$
Let $s_i=t_i+q_i \in \shA^{0,0}(U_i, \Theta_X)\otimes J$, then
$e^{s_i}=\psi_{|U_i}$ and so we have locally lifted the
isomorphism $\psi$.

Next, we  prove that the automorphisms $e^{s_i}$ can be glued
together to obtain an automorphism $e^s$ of $\shA_X^{0,*}\otimes
A$ that lifts $\psi$. Consider the intersection $U_{ij}$, then the
isomorphisms coincide on  $\shO_B(x)(U_{ij})$, i.e.,
$e^{s_i}_{|U_{ij}}=\psi
_{|U_{ij}}=e^{s_j}_{|U_{ij}}:\shO_B(x)(U_{ij}) \lrg
\shO_B(y)(U_{ij})$. Therefore, the isomorphism $e^{s_i-s_j}$ is
the identity on $\shO_B(x)(U_{ij})$. Since the action of
$\shA_X^{0,0}(\Theta_X) \otimes m_B$  on $A_X^{0,0}(U_{ij})
\otimes B$ is faithful on $\shO_X(U_{ij})\otimes B$, it follows
that $(s_i-s_j)_{|U_{ij}}=0$.

\end{proof}

\subsection{$KS_X$ controls the infinitesimal deformations of $X$}
\label{sez dimo teo DEF_K=DEF_X}

This section is devoted to prove that the Kodaira-Spencer algebra
of a complex manifold $X$ controls the infinitesimal deformations
of $X$.

\bigskip

\begin{teo} \label{teo def_k =Def_X}
Let $X$ be a complex compact manifold and ${KS_X}$ its
Kodaira-Spencer algebra. Then there exists an isomorphism of
functors
$$
 \gamma': \Def_{KS_X} \lrg \Def_X,
$$
defined in the following way: given a local Artinian $\C$-algebra
$(A,m_A)$ and a solution of the Maurer-Cartan equation $x \in
A^{0,1}_X(\Theta_X) \otimes m_A$ we set
$$
\shO_A(x)=\ker(\shA^{0,0}_X\otimes
A\xrightarrow{\debar+\bl_x}\shA^{0,1}_X\otimes A),
$$
and the map $\shO_A(x) \lrg \shO_X$ is induced by the projection
$\shA^{0,0}_X \otimes A \lrg \shA^{0,0}_X \otimes
\C=\shA^{0,0}_X$.
\end{teo}

\begin{oss}
As observed in Section~\ref{subsction deformation of non singular
varieties}, a deformation of $X$ can be interpreted as a morphism
of sheaves of algebras $\shO_A\lrg \shO_X$ such that $\shO_A$ is
flat over $A$ and $\shO_A\otimes_A \mathbb{C}\lrg \shO_X$ is an
isomorphism.

Consequently, the first part of the following proof consists of
showing the $A$-flatness of $\shO_A(x)$ and the existence of an
isomorphism $\shO_A(x)\otimes_A \mathbb{C}\cong\shO_X$.
\end{oss}

\begin{proof}
For each $(A,m_A) \in \Art$ and $x \in MC_{KS_X}(A)$, we have
defined
$$
\shO_A(x)=\ker(\shA^{0,0}_X\otimes
A\xrightarrow{\debar+\bl_x}\shA^{0,1}_X\otimes A).
$$
First of all, we observe that the projection $\pi$ on the residue
field $A/m_A\cong \C$ induces the following commutative diagram
\begin{center}
$\xymatrix{0 \ar[r]  & \shO_A(x) \ \ \ar[r] \ar[d]^{\pi} &
\shA^{0,0}_X \otimes A  \ar[r]^{\debar} \ar[d]^{\pi} &   \cdots
\ar[r]^{\debar} & \shA^{0,n}_X \otimes A \ar[r] \ar[d]^{\pi} &0 \\
& \shO_A(x)\otimes_A \C \ar[d] \ar[r] & \shA^{0,0}_X \ar[d]^{id}
\ar[r]^{\debar} &
\cdots   \ar[r]^{\debar} & \shA^{0,n}_X  \ar[d]^{id} \ar[r]& 0\\
0 \ar[r] & \shO_X  \ar[r] & \shA^{0,0}_X \ar[r]^{\debar}  & \cdots
  \ar[r]^{\debar} & \shA^{0,n}_X \ar[r]& 0.
\\ }$
\end{center}
Then $\pi$ induces the morphism  $\shO_A(x) \lrg \shO_X$.

\smallskip

Using Lemma~\ref{lem x in MC local gauge 0}, the Maurer-Cartan
solution $x$ is locally gauge equivalent to zero, therefore, there
exist a cover $\ \mathcal{U}=\{U_i\}$ and elements $a_i \in
 A^{0,0}(U_i,\Theta_X)\otimes m_A$ such
that $e^{a_i}*x_{|U_i}=0$, for each $i$. Thus, by Lemma~\ref{lemma
e^a(d+l_x)e^-a = d+e^a *x}, $e^{a_i} \circ(\debar +\bl_{x_{|U_i}})
\circ e^{-a_i}=e^{a_i}*(\debar
+\bl_{x_{|U_i}})=\debar+e^{a_i}*\bl_{x_{|U_i}}=\debar$ and so we
have the following commutative diagram
\begin{center}
$\xymatrix{0 \ar[r] & \shO_A(x)(U_i) \ar[r] \ar[d]^{e^{a_i}} &{
\shA^{0,0}_X(U_i)\otimes A }\ar[r]^{\ \ \ \ \ \ \ \ \
\debar+\bl_{x_{|U_i}}} \ar[d]^{e^{a_i}} & \cdots
\ar[r]^{\debar+\bl_{x_{|U_i}}\ \ \ \ \ \ \ } &
\shA^{0,n}_X(U_i)\otimes A
\ar[r] \ar[d]^{e^{a_i}} &0 \\
0 \ar[r] & \shO_X(U_i)\otimes A \ar[r] & \shA^{0,0}_X(U_i)\otimes
A \ar[r]^{\ \ \ \  \ \ \  \debar} &   \cdots   \ar[r]^{\debar\ \ \
\ \ \ \ \ } & \shA^{0,n}_X(U_i)\otimes A
\ar[r]& 0, \\
}$
\end{center}
where the vertical arrows are isomorphisms.

This implies that the deformation $\shO_A(x)$  is locally trivial,
i.e., $\shO_A(x)(U_i) \cong \shO_X(U_i) \otimes A$. Since
$\shO_X(U_i) \otimes A$ is flat over $A$,   $\shO_A(x)(U_i)$ is
also $A$-flat. Since flatness is a local property, $\shO_A(x)$ is
$A$-flat.

 Using  the isomorphism $\shO_A(x)(U_i) \cong \shO_X(U_i)
\otimes A$ we can also conclude that  $ \shO_A(x)(U_i) \otimes_A
\C \cong \shO_X(U_i)$  and so    $\shO_A(x) \otimes_A \C \lrg
\shO_X$ is    an isomorphism.


\medskip

Then it is well defined the following  morphism of functors of
Artin rings
$$
\gamma : \MC_{KS_X} \lrg \Def_X,
$$
such that
$$
\gamma(A): \MC_{KS_X}(A) \lrg \Def_X(A)
$$
$$
  x\longmapsto \shO_A(x).
$$
\smallskip

Next, we prove that the deformations $\shO_A(x)$ and $\shO_A(y)$
are isomorphic if and only if  $x,y \in \MC_{KS_X}(A)$ are gauge
equivalent.

Actually, if $\shO_A(x) \cong \shO_A(y)$,  applying
Proposition~\ref{lem sollevamento ISO da O_A a complexes},  we can
lift this isomorphism to an isomorphism of complexes $e^s$ with $s
\in A_X^{0,0}\otimes m_A$, i.e.,
\begin{center}
$\xymatrix{0 \ar[r] & \shO_A(x) \ar[r] \ar[d]^{\cong} &
\shA^{0,0}_X\otimes A \ar[r]^{\ \ \  \ \debar+\bl_x} \ar[d]^{e^s}
& & \cdots & \ar[r]^{\debar+\bl_x\ \ \ \ } & \shA^{0,n}_X\otimes A
\ar[r]
\ar[d]^{e^s} &0\\
0 \ar[r] & \shO_A(y) \ar[r]  & \shA^{0,0}_X\otimes A \ar[r]^{\ \ \
\ \debar+\bl_y} & & \cdots & \ar[r]^{\debar+\bl_y \ \ \ \ } &
\shA^{0,n}_X\otimes A \ar[r]& 0. \\ }$
\end{center}
The commutativity of the diagram  and Lemma~\ref{lemma
e^a(d+l_x)e^-a = d+e^a *x} imply that  $\debar + \bl_y
=e^{-s}\circ (\debar +\bl_x) \circ e^{s}= \debar + e^s*\bl_x$.
Therefore, $e^s*x=y$.

\medskip

In conclusion, the map $\gamma'$, induced by  $\gamma$ on
$\Def_{KS_X}=MC_{KS_X} / gauge$, is a well defined injective
morphism:
$$
\gamma': \Def_{KS_X} \hookrightarrow \Def_X.
$$
To conclude that $\gamma'$ is an isomorphism we prove that
$\gamma'$ is   \'{e}tale (and so $\gamma'$ is surjective).

Using Corollary~\ref{cor etale se bietti su tg e inj su ostru}, we
need to prove  that:

$1)$ $\gamma'$ induces a bijective map on the tangent spaces;

$2)$ $\gamma'$ induces an injective map on the obstruction spaces.

\medskip

As to  $\Def_{KS_X}$, by Remark~\ref{oss tangente DEF_L =H^1}, the
tangent space is isomorphic to $ H^1_{\debar }(X,\Theta_X)$ and
Lemma~\ref{lem calcolo ostruzione MC_L=H^2(L)} implies that the
obstructions are naturally contained in $ H^2_{\debar
}(X,\Theta_X)$. As to $\Def_X$, Theorems~\ref{teo sernesi first
ordine=H^1(X,T_X)} and~\ref{teo sernesi ostruz DEF_X=H^2(X,T_X)}
show that the tangent space is isomorphic to
$\check{H}^1(X,\Theta_X)$ and the obstructions are naturally
contained in $\check{H}^2(X,\Theta_X)$.

Then we will prove that the maps induced by $\gamma'$ coincide
with the Leray isomorphisms (see Theorem~\ref{teo
VOISINh^q=H^check q} and Remark~\ref{oss associo Hcheck iso H
dolebeautl}).

\bigskip

\emph{1) Tangent Spaces.} Let us  prove that the map
$\gamma_{\varepsilon} '$ induced by $\gamma'$ on the tangent space
$$
\gamma_\varepsilon ': \Def_{KS_X}(\C[\epsi]) \lrg
\Def_X(\C[\epsi])
$$
is the Leray isomorphism.

By Remark~\ref{oss tangente DEF_L =H^1}, we have
$\Def_{KS_X}(\C[\epsi])=H^1 (X,\Theta _X) $. Proceeding as in
Remark~\ref{oss associo Hcheck iso H dolebeautl}, there exists a
Stein cover $\mathcal{U }=\{U_i\}$ so that  we can associate with
each $x \in \Def_{KS_X}(\C[\epsi])$ an element
$[\sigma]=[\{\sigma_{ij}=(a_i- a_j)_{|U_{ij}}\}]\in
\check{H}^1(X,\Theta_X)\otimes \C\epsi$, with $ x_{|U_i}=\debar
a_i$, that does not  depend on the choice of $a_i$.

\medskip

Next, let $\gamma(x)=\shO_{\C[\epsi]}(x)$ be the deformation
associated with $x$, i.e.,
$$
0 \lrg \shO_{\C[\epsi]}(x) \lrg \shA^{0,0}_X\otimes \C[\epsi]
\xrightarrow{\debar+\bl_x}\shA^{0,1}_X\otimes \C[\epsi] \lrg
\cdots
 .
$$
As before, the deformation $\shO_{\C[\epsi]}(x)$ is locally
trivial; thus,  there exists $b_i \in \shA^{0,0}(U_i, \Theta_X)
\otimes \C\epsi$  such that ${e^{b_i}*x}_{|U_i}=0$ and so
\begin{center}
$\xymatrix{\shO_{\C[\epsi]}(x)(U_i)  \ar[rr]^{e^{b_i}}_{\cong} &
 & \shO_X(U_i) \otimes \C[\epsi].
 \\ }$
\end{center}
Proceeding as in the proof of Theorem~\ref{teo sernesi first
ordine=H^1(X,T_X)}, for all $i$ and $j$
$$
\varphi_{ij}:= e^{b_i-b_j}: \shO_X(U_{ij}) \otimes \C[\epsi] \lrg
\shO_X(U_{ij}) \otimes \C[\epsi]
$$
is an automorphism of the trivial deformation $\shO_X(U_{ij})
\otimes \C[\epsi]$ that restricts to the identity.

\noindent Applying Lemma~\ref{lemm serneII.1.5 automor-derivazio},
the class $[\{\tau_{ij}\}]=[\{(b_i-b_j)_{|U_{ij}}\}]\in
\check{H}^1(X,\Theta_X) \otimes \C\epsi $ is the Check 1-cocycle
associated with the deformation $\gamma'_\epsi(x)$.

Since ${e^{b_i}*x}_{|U_i}=0$,   $\debar b_i = {x}_{|U_i}=\debar
a_i$ and so $b_i=a_i +c_i$ with $c_i \in \Gamma(U_i,
\Theta_X)\otimes \C\epsi$.

Therefore, $[\tau_{ij}]=[\{(b_i-b_j)_{|U_{ij}}\}]=[\sigma]$. This
shows that $\gamma'_\epsi $ coincides with the Leray's
isomorphism.




\bigskip

\emph{2) Obstruction}

Let
$$
0 \lrg J \lrg B \stackrel{\alpha}{\lrg} A \lrg 0
$$
be a small   extension.


First we consider the  obstruction class  $[h]$ of $ x \in
\Def_{KS_X}$.

Let $x \in \Def_{KS_X}(A)$, and $\tilde{x}\in {KS_X}^1 \otimes
m_B$ be a lifting of $x$. The  obstruction class associated with
$x$ is $[h] \in H^2({KS_X}) \otimes J$ with
$$
h=\debar \tilde{x} + \displaystyle\frac{1}{2}[\tilde{x},\tilde{x}]
$$
and this class does not  depend on the choice of the lifting
$\tilde{x}$ as it was shown in Lemma~\ref{lem calcolo ostruzione
MC_L=H^2(L)}.

Proceeding as in Remark~\ref{oss associo Hcheck iso H dolebeautl},
there exists a Stein cover $\mathcal{U}=\{U_i\}$ such that the
class $[\{\alpha_{ijk}\}]=[\{\rho_{jk}-\rho_{ik}+ \rho_{ij}\}]\in
\check{H}^2(X,\Theta_X)\otimes J$ is the class associated with $h$
via the Leray's isomorphism, where $h_{|U_i}=\debar \tau_i$ and
$\debar\rho_{ij}= (\tau_i - \tau_j)_{|U_{ij}}$. In particular, we
note that $ \rho_{ij} \in  \Gamma (U_{ij},\shA^{0,0}(\Theta_X))
\otimes J$.



Since $h_{|U_i}=\debar \tau_i$, $x$ can be locally lifted to a
solution of the Maurer-Cartan equation ,
$$
\overline{x}_i=\tilde{x}_{|U_i}-\tau_i \in A^{0,1}(U_i,\Theta_X)
\otimes m_B.
$$
Indeed, on $U_i$  we have

$$
\alpha(\overline{x})= \alpha (\tilde{x})=x
$$
and
$$
 \debar \overline{x} + \displaystyle\frac{1}{2}
[\overline{x},\overline{x}]= \debar \tilde{x} -\debar \tau_i +
\displaystyle\frac{1}{2}[\tilde{x},\tilde{x}]=h_{|U_i}-\debar
\tau_i=0.
$$

Moreover,   ${e^{\rho_{ij}}*\overline{x}_j}_{|U_{ij}}=
{\overline{x}_i}_{|U_{ij}}$; more precisely (see Example~\ref{exe
azione gauge se []=0}), we have
$$
{e^{\rho_{ij}}*\overline{x}_j}_{|U_{ij}}=
e^{\rho_{ij}}*(\tilde{x}-\tau_j)_{|U_{ij}}=(-\debar \rho_{ij}+
\tilde{x}-\tau_j)_{|U_{ij}}=(\tilde{x}-\tau_i)_{|U_{ij} }=
{\overline{x}_i}_{|U_{ij}}.
$$

As above,   $x$ is locally equivalent to zero; therefore, for each
$i$ there exists $a_i \in A^{0,0}(U_i,\Theta_X)\otimes m_A$, such
that $e^{a_i} *x_{|U_i}=0$.

Analogously, for each $i$, there exists $b_i \in
A^{0,0}(U_i,\Theta_X) \otimes m_B$ that is a lifting of $a_i$ such
that
$$
e^{b_i} *\overline{x}_i =0.
$$

\medskip

Next, let $\gamma'(x)=\shO_A(x)$ be the deformation of $X$ induced
by $x$. As above the deformation is locally trivial and so  there
exist a cover $ \mathcal{U}=\{U_i\}$  and $a_i \in
\shA^{0,0}(U_i,T_X)\otimes m_A$ such that
$$
\shO_A(x)(U_i)\stackrel{e^{a_i}}{\cong}  \shO_X(U_i)\otimes A.
$$
Let $\varphi_{ij}$ be the following isomorphism
$$
\varphi_{ij}: \shO_X(U_{ij}) \otimes A \stackrel{e^{-a_j}}{\lrg}
\shO_A(x)(U_{ij})\stackrel{e^{a_i}}{\lrg} \shO_X(U_{ij}) \otimes
A.
$$

Proceeding as in the proof of Theorem~\ref{teo sernesi ostruz
DEF_X=H^2(X,T_X)}, since $b_i \in \shA^{0,0}(U_i,\Theta_X)\otimes
m_B$ are  liftings   of $a_i$,   $\tilde{\varphi_{ij}}= e^{-b_j
}e^{\rho_{ij}}e^{  b_i}\in \Aut(\shO_X(U_{ij})\otimes B)$, defined
as
$$
\tilde{\varphi}_{ij}: \shO_X(U_{ij}) \otimes B
\stackrel{e^{-b_j}}{\lrg}
\shO_B(\overline{x}_j)(U_{ij})\stackrel{e^{\rho_{ij}}}{\lrg}
\shO_B(\overline{x}_i)(U_{ij}) \stackrel{e^{b_i}}{\lrg}
\shO_X(U_{ij}) \otimes B,
$$
is a lifting of $\varphi_{ij}$.

By Remark~\ref{oss succ esatta lunga automorfismi derivazione},
the automorphisms $e^{\rho_{ij}}$, $e^{\rho_{ik}}$ and
$e^{\rho_{jk}}$ commute with the other automorphisms. Then
$\Phi_{ijk} =\tilde{\varphi_{jk}}{ \tilde{\varphi_{ik}}}^{-1}
\tilde{\varphi_{ij}}=e^{\rho_{jk}- \rho_{ik}+\rho_{ij}}$ is an
automorphism of the trivial deformation  that restricts to the
identity $({\Phi_{ijk}}_{|\shO_X(U_{ijk}) \otimes A}=id)$.
Therefore, by Lemma~\ref{lemm serneII.1.5 automor-derivazio}, the
element $[\{\alpha_{ijk}\}]=[\{\rho_{jk}-\rho_{ik}+
\rho_{ij}\}]\in \check{H}^2(\mathcal{U},\Theta_X)\otimes J$ is the
obstruction class associated with $\gamma'(x)$. In conclusion, in
the obstruction case too, the map induced by $\gamma'$ coincides
with the Leray's isomorphism.

\end{proof}

\subsection{Deformations of a product}

As an application of  Theorem~\ref{teo def_k =Def_X} we study
deformations of a product of compact complex manifolds $X$ and
$Y$. The following remark will be used in Section~\ref{sez teo
Def_(h,g)Def (f)}.

\begin{oss}\label{oss DEF(XxY)> DEF(X) x DEF(Y)}
In general, not all the deformations of the product $X \times Y$
are  products  of deformations of $X$ and of $Y$.
\end{oss}

The first example of this fact was given by Kodaira and Spencer in
their work (\cite[pag.~436]{bib Kodaira SpencerII}) where they
provided one of the first examples of obstructed varieties. More
precisely, they considered the product of the projective line and
the complex torus of dimension $q\geq2$ $\mathbb{P}^1 \times
\C^q/G$ and  proved that it is obstructed though the two manifolds
are unobstructed.

A sufficient and necessary condition to have an isomorphism
between products of deformations and deformations of the product
is  given by the lemma below.

\begin{lem}
The morphism
$$
F:\Def_{ X}\times \Def_{ Y} \lrg \Def_ {X \times Y}
$$
is an isomorphism if and only if $H^1(\shO_X)\otimes
H^0(\Theta_Y)=H^1(\shO_Y)\otimes H^0(\Theta_X)=0$.
\end{lem}
\begin{proof}
By Theorem~\ref{teo def_k =Def_X}, it suffices to prove that the
morphism
$$
F:\Def_{KS_X}\times \Def_{KS_Y} \lrg \Def_{KS_{X \times Y}}
$$
is an isomorphism.

Let $p$ and $q$ be the natural projections of the product $X
\times Y$   onto $X$ and $Y$, respectively. Consider the morphism
of DGLAs
$$
F: KS_X \times KS_Y \lrg KS_{X \times Y},
$$
$$
(n_1,n_2)\longmapsto p^* n_1 +q^*n_2.
$$
Thus, we have to prove that $F$ induces an isomorphism of the
associated deformation functors. Denote by $H^i(F)$ the induced
map on cohomology.

By Theorem~\ref{teo no EXTE iso H^i allro iso DEF},  if $H^0(F)$
is surjective, $H^1(F)$ is bijective and $H^2(F)$ is injective,
then $F$ induces an isomorphism.

$H^0(KS_{X \times Y})=H^0(X\times Y,\Theta_{X \times Y})$ and so,
by the Kunneth's formula, $H^0(KS_{X \times Y})=
H^0(\Theta_X)\oplus H^0(\Theta_Y)= H^0( KS_X \times KS_Y)$. This
implies that $H^0(F)$ is surjective.

Again by  Kunneth's formula we get
$$
H^1(KS_{X \times Y})\cong H^1(\Theta_{X \times Y})\cong
$$
$$
H^1(\Theta_X) \oplus H^1(\Theta_Y)\oplus H^1(\shO_X)\otimes
H^0(\Theta_Y)\oplus H^1(\shO_Y)\otimes H^0(\Theta_X),
$$
and
$$
 H^1( KS_X \times KS_Y)=H^1(\Theta_X) \oplus H^1(\Theta_Y).
$$
Then the hypothesis implies that $H^1(F)$ is bijective.

Finally, reasoning as above we get that
$$
H^2(F): H^2(\Theta_X) \oplus H^2(\Theta_Y)\lrg H^2(\Theta_{X\times
Y})
$$
is injective. This implies that $F$ induce an  isomorphism of
deformation functors.

On the other hand, if $F$ is an  isomorphism of deformation
functors, then $H^1(F)$ is a bijection onto the tangent spaces and
so $H^1(\shO_X)\otimes H^0(\Theta_Y)=H^1(\shO_Y)\otimes
H^0(\Theta_X)=0$.
\end{proof}

\chapter{Deformation functor of a pair of morphisms of
DGLAs}

In this   chapter we give the key definition of deformation
functor associated with a pair of morphisms of differential graded
Lie algebras.

In the first section   we define the (not extended) functors of
Artin rings $\MC_{(h,g)}$ (Definition~\ref{def funtore MC_(h,g)
non estesi}) and  $\Def_{(h,g)}$ (Definition~\ref{defin funtor
Def_(h,g) non esteso}) associated with a pair $h :L \lrg M$ and
$g: N \lrg M$. These functors will play an important role in the
study of infinitesimal deformations of holomorphic maps in the
next chapter.

Then Section~\ref{sezio extended defo funtore} is devoted
introducing the extended deformation functors
(Definition~\ref{dhef extended defo functor}).
In particular, we define the functors $\widetilde{\MC}_{(h,g)}$
and $ \widetilde{\Def}_{(h,g)} $ which are a generalization of the
previous $\MC_{(h,g)}$   and  $\Def_{(h,g)}$.

We introduce the extended  functors, since using their properties,
we can show the existence of a DGLA $\hil$ such that
${\Def}_{\hil}\cong  {\Def}_{(h,g)} $ (Theorem~\ref{teo EXT
DEF_H=DEF(h,g)}).




\section{Functors $\MC_{(h,g)}$ and $\Def_{(h,g)}$}
\label{sezio non ext def(h,g) e MC(h,g)}

In this section  we introduce  the   functors $\MC_{(h,g)}$ and
$\Def_{(h,g)}$ (Section~\ref{sezi def MC_(h,g) and DEF(h,g) no
ext}) associated with a pair $h :L \lrg M$ and $g: N \lrg M$ and
we study some of their properties (Section~\ref{sez tangent e
ostru MC(h,g) DEF(h,g)}). First of all, we  recall  the definition
of the mapping cone associated  with morphisms of complexes.

\bigskip

\noindent{\bf Note.} In this section, we suppose that $M$ is
concentrated in non negative degree. Since in the main application
$M$ will be the Kodaira-Spencer algebra of a manifold, this extra
hypothesis is not restrictive. Anyway, in Section~\ref{sezio
extended defo funtore} we will remove this hypothesis.

\subsection{The mapping cone of a pair of
morphisms}\label{sezione triplo cono}

The \emph{suspension of the mapping cone of a morphism  of
complexes} $h:(L,d ) \lrg (M,d )$ is the differential graded
vector space $(C^\cdot_h,\delta)$, where
$$
C_h^i=L^i\oplus M^{i-1}
$$
and the differential $\delta$ is
$$
\delta(l,m)=( d  l, -d m+h(l)) \in L^{i+1}\oplus M^i, \qquad
\forall \, (l,m) \in L^i\oplus M^{i-1} .
$$
Actually, $\delta^2(l,m)=\delta(  d  l, -d  m+h(l)) =(d ^2 l, d ^2
m-d(h(l))+h(d(l)))=(0,0)$.

\begin{oss}\label{oss ho mappa da C_h ->coker (h)}
Consider the projection $\pi:M \lrg \coker(h)$. Then there exists
a morphism of complexes
$$
\varphi: (C^\cdot_h,\delta)\lrg (\coker(h)[-1],d_{[-1]}),
$$
with
$$
\varphi(l,m)= \pi(m), \qquad \forall \ (l,m)\in C_h^i.
$$
Indeed,
\begin{center}
$\xymatrix{(l, m) \ar[r]^{\varphi} \ar[d]_\delta  &
\pi(m) \ar[d]^{d_{[-1]}}\\
(d l, -d m + h(l)) \ar[r]^\varphi & -  d \pi(m). \\
}$
\end{center}
If $h$ is injective then $(C^\cdot_h,\delta)$ and
$(\coker(h)[-1],d_{[-1]})$ are quasi-isomorphic.

\end{oss}

\begin{oss}
Let $f:X \lrg Y$ be a  holomorphic map of complex compact
manifolds and consider the induced map $ f_*: A_X^{p,q}(\Theta_X)
\lrg A_X^{p,q}(f^*\Theta_Y) $ (see Section~\ref{sezione f_* and
f^*}). Then the cohomology groups $H^i(C^\cdot _{f_*})$ of the
associated cone $(C^\cdot _{f_*},\delta)$ are isomorphic to the
hyper-cohomology groups $\mathbb{H}^1\left(X,\shO(\Theta_X)
\stackrel{f_*}{\lrg}\shO(f^*\Theta_Y)\right)$ (see \cite[Section 5
of Chapter 3]{bib Griffone}.
\end{oss}

\bigskip

Next, suppose that $h:(L,d ) \lrg (M,d )$ and $g:(N,d ) \lrg (M,d
)$ are morphisms of complexes:
\begin{center}
$\xymatrix{  & L \ar[d]^h\\
N\ar[r]^g & M. \\
}$
\end{center}

\begin{dhef}
The \emph{suspension of the mapping cone of a  pair of morphisms
 $(h,g)$ } is the differential graded vector space $(\cil^\cdot,D)$, where
$$
\cil^i=L^i \oplus N^i \oplus M^{i-1}
$$
and the differential $D $ is defined as follows
$$
L^i \oplus N^i \oplus M^{i-1} \ni (l,n,m) \stackrel{D}{\lrg} ( d
l, d n,-d m-g(n)+h(l)) \in L^{i+1} \oplus N^{i+1} \oplus M^i.
$$

Actually, $D^2(l,n,m)=D( d  l,d n,-d m-g(n)+h(l))=(d^2 l,d^2 n,d^2
m+dg(n)-d h(l)-g(d n)+h(d l) )=(0,0,0)$.
\end{dhef}

\begin{oss}
By definition, $(\cil^\cdot, D)$ coincides with  the suspended
mapping cone associated with the morphism of complexes $h-g: L
\oplus N \lrg M$ (such that $(h-g)(l,n)=h(l)-g(n)$).

\smallskip

Moreover, the projection $\cil^\cdot \lrg L^\cdot \oplus N^\cdot$
is a morphism of complexes and so there exists the following exact
sequence
\begin{center}
$ \xymatrix{0 \ar[r] & (M^{\cdot-1},-d ) \ar[r]   &
(C^\cdot_{(h,g)},D) \ar[r]   & (L^\cdot \oplus N^\cdot,d) \ar[r] &
0
\\ }$
\end{center}
that induces
\begin{equation}\label{success esatta lunga omolog cil(fg)}
\cdots \lrg H^i(\cil^\cdot) \lrg H^i(L^\cdot\oplus N^\cdot) \lrg
H^i(M^\cdot) \lrg H^{i+1}(\cil^\cdot) \lrg \cdots .
\end{equation}
\end{oss}

\begin{oss}
The complexes $\cil^\cdot$ and $\operatorname{C}^\cdot_{(g,h)}$
are
  isomorphic. Actually, let $\gamma :\cil^\cdot\lrg
\operatorname{C}^\cdot_{(g,h)}$ be defined as
$$
\cil^i \ni (l,n,m) \stackrel{\gamma}{\longmapsto} (-l,-n,m) \in
\operatorname{C}^i_{(g,h)}.
$$
Then
\begin{center}
$\xymatrix{(l,n,m) \ar[r]^{\gamma} \ar[d]_D  & (-l,-n,m)
\ar[d]^{\delta} \\
(d l, d  n, -d m -g(n) + h(l)) \ar[r]^\gamma     &
(-d l,-dn,-d m -g(n) + h(l)), \\
}$
\end{center}
and so $\gamma$ is a well defined morphism of complexes that is a
quasi-isomorphism ($\gamma^2=id$).

\end{oss}

\begin{lem}\label{lem h inj then C_(h,g)=C_pi dopo g}
Let $g:N \lrg M$ and $h:L \lrg M$ be morphisms of complexes with
$h $ injective, i.e., there exists the exact sequence of complexes
\begin{center}
$\xymatrix{0 \ar[r] &L  \ar[r]^h &  M \ar[r]^{\pi\ } & \coker
(h)\ar[r] & 0 \\
& & N.\ar[u]_g & & \\
}$
\end{center}
Then $(\cil^\cdot,D)$ is quasi-isomorphic to $(C^\cdot_{\pi \circ
g},\delta)$.

\end{lem}
\begin{proof}

Let $\gamma:\cil^\cdot \lrg C^\cdot_{ \pi \circ g }$ be defined as
$$
\cil^i \ni (l,n,m) \stackrel{\gamma}{\longmapsto} (-n,\pi(m)) \in
C^i_{ \pi \circ g },
$$
then
\begin{center}
$\xymatrix{(l,n,m) \ar[r]^{\gamma} \ar[d]_D  & (-n,\pi (m))
\ar[d]^{\delta} \\
(d l, d  n, -d m -g(n) + h(l)) \ar[r]^\gamma     &
(-d n,  -d \pi(m)-\pi(g(n))). \\
}$
\end{center}

Therefore, $\gamma$ is a well defined morphism of complexes that
is a quasi-isomorphism. The fact that   $\gamma$  induces
isomorphisms in cohomology is an easy calculation but we state it
for completeness.
\smallskip

Denote by the same $\gamma$ the map induced in cohomology.

$\gamma$ is injective. Suppose that $\gamma[(l,n,m)]=
[(-n,\pi(m))]$ is zero in $H^i( C^\cdot_{\pi \circ g})$. Then
$dl=dn=-dm-g(n)+l =0$ and there exists $(b,c) \in N^{i-1} \times
\coker(h)^{i-2}$ such that $-n=db$ and $\pi(m)= -dc +\pi \circ
g(b)$. Let $m'\in M^{i-2}$ be  a lifting of $c$, i.e.,
$\pi(m')=c$, $n'=-b \in N^{i-1}$ and $l'=m+dm'+g(n') \in M^{i-1}$.
Then $l' \in L^{i-1}$, actually, $\pi(l')=\pi(m)+\pi(dm')
+\pi\circ g(n')=-dc +\pi\circ g(b) +dc -\pi\circ g(b)=0$.
Therefore, $D(l',n',m')= (dl',  dn', -dm'-g(n')+l') = (dm-g(db),
-db,-dm'+g(b)+( m+dm'+g(n') )) =(l,n,m)$.

\smallskip

$\gamma$ is surjective. Let $[(n,t)] \in H^i( C^\cdot_{\pi \circ
g})$; then $dn=0$ and $-dt+ \pi \circ g (n)=0$. Let $m\in M^{i-1}$
be a lifting of $ t$, i.e., $\pi(m)=t$ and $l=-g(n)+dm \in M^i$.
Then $l \in L^i$; in fact,  $\pi(l)=-\pi \circ g(n) +dt=0$.
Moreover, $(l,-n,m) \in H^i(\cil^\cdot)$ ($dl=dn=0$ and
$-dm+g(n)+l=0$) and $\gamma [ (l,-n,m) ] =[n,\pi(m)]=[(n,t)]$.

\end{proof}

\begin{oss}

Let $L,M$ and $N$ be DGLAs and $h :L \lrg M$ and $g: N \lrg M$
morphisms of DGLAs. No canonical DGLA structure  on $\cil^\cdot$
can be defined such that the projection $\cil ^\cdot \lrg L^\cdot
\oplus N^\cdot$ is a morphism of DGLAs.

\end{oss}

\subsection{Definition of  $\MC_{(h,g)}$ and $\Def_{(h,g)}$ }
\label{sezi def MC_(h,g) and DEF(h,g) no ext}

\begin{dhef}\label{def funtore MC_(h,g) non estesi}
Let $h :L \lrg M$ and $g: N \lrg M$ be morphisms of differential
graded Lie algebras:
\begin{center}
$\xymatrix{ & L \ar[d]^h   \\
N\ar[r]^{g}   & M. \\ }$
\end{center}
For each $(A,m_A) \in \Art$  the \emph{Maurer-Cartan functor
associated with the pair $(h,g)$ }   is defined as follows
$$
\MC_{(h,g)}: \Art \lrg \Set,
$$
$$
\MC_{(h,g)}(A)=\{(x,y,e^p) \in (L^1 \otimes m_ A)\times  (N^1
\otimes m_A ) \times \operatorname{exp}(M^0 \otimes m_A)  |
$$
$$
dx +\frac{1}{2}[x,x]=0,\   dy +\frac{1}{2}[y,y]=0,\  g(y)=e^p*h(x)
\}.
$$
\end{dhef}

\begin{oss}
In \cite[Section~2]{bib ManettiPREPRINT}, M.~Manetti defined the
functor $\MC_h$ associated with a morphism $h:L \lrg M$ of DGLAs:
$$
\MC_{h}: \Art \lrg \Set,
$$
$$
\MC_{h}(A)=
$$
$$
\{(x, e^p) \in (L^1 \otimes m_ A)  \times \operatorname{exp}(M^0
\otimes m_A)  |\  dx +\frac{1}{2}[x,x]=0,\, e^p*h(x)=0 \}.
$$
Therefore, if we take $N=0$ and $g=0$, the new functor
$\MC_{(h,g)}$ reduces to the old one $\MC_h$.

By choosing $N=0$ and $h=g=0$, $\MC_{(h,g)}$ reduces to the
Maurer-Cartan functor  $\MC_L$ associated with the DGLA $L$
(Definition~\ref{def Maurer-C funcotr MC_L of L}).
\end{oss}

\begin{oss}
As in the  case of a differential graded Lie algebra (see
Remark~\ref{oss MC_L is homogeneous}), $\MC_{(h,g)} $ is  a
\emph{homogeneous} functor, since $\MC_{(h,g)}(B \times _A C)
\cong \MC_{(h,g)}(B)\times _{\MC_{(h,g)}(A)}\MC_{(h,g)}(C)$.
\end{oss}

As in the case of a differential graded Lie algebra (see
Definition~\ref{definizio gauge L DGLA}), we can define for each $
(A,m_A) \in \Art$ a gauge action over $\MC_{(h,g)}(A)$.

\begin{dhef}\label{def gauge action MC_(h,g) NO EXT}
The gauge action  of $\exp(L^0 \otimes m_A) \times \exp(N^0
\otimes m_A)$ over $\MC_{(h,g)}(A)$ is given by:
$$
(e^a,e^b)*(x,y,e^p)= (e^a*x,e^b*y,e^{g(b)}e^pe^{-h(a)}).
$$
This is well defined since
$$
e^{g(b)}e^pe^{-h(a)}*h(e^a*x)= e^{g(b)}e^p*h(x)=e^{g(b)} *g(y)=
g(e^b*y).
$$
and so $(e^a,e^b)*(x,y,e^p) \in \MC_{(h,g)}(A)$.
\end{dhef}

In conclusion, it makes sense to consider the following functor.

\begin{dhef}\label{defin funtor Def_(h,g) non esteso}
The \emph{deformation functor  associated with a  pair $(h,g)$ }
of morphisms of  differential graded Lie algebras is:
$$
\Def_{(h,g)}:\Art \lrg \Set,
$$
$$
\Def_{(h,g)}(A)= \frac{\MC_{(h,g)}(A)}{  \exp(L^0 \otimes m_A)
\times \exp(N^0 \otimes m_A) }.
$$
\end{dhef}

\begin{oss}\label{oss DEf(hg) reduces to DEF_h e DEF L}
In \cite[Section~2]{bib ManettiPREPRINT}, M.~Manetti defined the
functor $\Def_h$ associated with a morphism $h:L \lrg M$ of DGLAs:
$$
\Def_h: \Art \lrg \Set,
$$
$$
\Def_h(A)= \frac{\MC_{h}(A)}{ \exp(L^0 \otimes m_A) \times
\exp(dM^{-1} \otimes m_A)},
$$
where  the gauge action of $\exp(L^0 \otimes m_A) \times
\exp(dM^{-1} \otimes m_A)$ is given by the following formula
$$
(e^a,e^{dm})*(x,e^p)= (e^a*x,e^{dm} e^pe^{-h(a)}),\qquad \forall \
a \in  L^0 \otimes m_A, m \in  M^{-1} \otimes m_A.
$$

Therefore, if we take $N=0$ and $g=0$, the new functor
$\Def_{(h,g)} $ reduces to the old one $\Def_h$ ($M$ is
concentrated in non negative degree).

By choosing $N=M=0$ and $h=g=0$, $\Def_{(h,g)}$ reduces to the
Maurer-Cartan functor  $\Def_L$ associated with the DGLA $L$.
\end{oss}

The name deformation functor is justified by the theorem below.

\begin{teo}
$\Def_{(h,g)} $ is a deformation functors, i.e., it  satisfies the
conditions of Definition~\ref{defin deformation funtore}.
\end{teo}
\begin{proof}
If $A=\K$, then $\Def_{(h,g)}(\K)=\{\mbox{one element}\}$;
therefore, $ii)$ of Definition~\ref{defin deformation funtore}
holds, i.e., $\Def_{(h,g)}(B\times C)= \Def_{(h,g)}(B) \times
\Def_{(h,g)}(C)$.

\smallskip

Let $\beta:B \lrg A$ and $\gamma:C \lrg A$ be morphisms in $\Art$
and $(v,w) \in \Def_{(h,g)}(B) \times_{\Def_{(h,g)}(A)}
\Def_{(h,g)}(C)$. Then we are looking for a lifting $z  \in
\Def_{(h,g)}(B\times_A C)$, whenever $\beta: B \lrg A$ is
surjective.

Let $(x,y,e^p) \in \mc(B)$ and $(s,t,e^r)\in \mc(C)$ be liftings
for $v$ and $w$, respectively.

By hypothesis, $\beta(v)=\gamma(w) \in \Def_{(h,g)}(A)$;
therefore, $\beta(x,y,e^p)$ and $\gamma(s,t,e^r)$ are gauge
equivalent in $\mc(A)$, i.e., there exist $a \in L^0 \otimes m_A$
and $b \in N^0 \otimes m_A$ such that
$$
e^a * \beta(x)=\gamma(s) \ \ \ \ \ e^b*\beta(y)= \gamma(t)\ \ \ \
e^{g(b)} e^{\beta(p)}e^{-h(a)}=e^{\gamma(r)}.
$$
Let $c \in L^0 \otimes m_B$ such that $\beta(c)=a$ and  $d \in N^0
\otimes m_B$ such that $\beta(d)=b$.

Up to substituting $(x,y,e^p)$ with its gauge equivalent element
$(e^c*x,e^d*y,e^{g(d)}e^{p} e^{-h(c)})$ (they both lift $v$), we
can assume that\footnote{$\beta(e^c*x,e^d*y,
e^{g(d)}e^{p}e^{-h(c)})=(e^{\beta(c)}*\beta(x),e^{\beta(d)}*
\beta(y), e^{g(\beta(d))}e^{\beta(p)}e^{-h(\beta(c))})=(e^a *
\beta(x), e^b*\beta(y), e^{g(b)} e^{\beta(p)}e^{-h(a)} )
=\gamma(s,t,e^r)$ } $\beta(x,y,e^p)=\gamma(s,t,e^r) \in \mc(A)$.

Since $\MC_{(h,g)}(B\times_A C)=\MC_{(h,g)}(B)
\times_{\MC_{(h,g)}(A)} \MC_{(h,g)}(C)$,  a lifting $z \in
\MC_{(h,g)}(B\times_A C)$  is well defined and so it is sufficient
to take its class   $[z] \in \Def_{(h,g)}(B\times_A C)$.

\end{proof}

\begin{oss}\label{oss proj varrho induce no EXT Def(h,g)->Def_N}
Consider the functor $ \Def_{(h,g)}$. Then the projection
$\varrho$ on the second factor:
$$
\varrho: \Def_{(h,g)} \lrg \Def_N,
$$
$$
\Def_{(h,g)}(A) \ni (x,y,e^p)\stackrel{\varrho}{\lrg} y \in
\Def_N(A)
$$
is a morphism of deformation functors.
\end{oss}
\begin{oss}\label{oss DEF_(h,g) con h iniettivo}
If the morphism \emph{$h$ is injective}, then  for each $(A,m_A)
\in \Art$  the functor $\MC_{(h,g)}$ has the following form:
$$
\MC_{(h,g)}(A)=\{(x,e^p) \in  (N^1 \otimes m_A)\times
\operatorname{exp}(M^0 \otimes m_A)  |
$$
$$
dx +\frac{1}{2}[x,x]=0,\  e^{-p}*g(x) \in L^1 \otimes m_A \}.
$$
In this case the gauge equivalence $\sim$ is given by
$$
(x,e^p)\sim (e^b*x,e^{g(b)} e^p e^{a}), \qquad \mbox{ with } a \in
L^0 \otimes m_A \mbox{ and } b \in N^0 \otimes m_A.
$$

\end{oss}

\begin{lem}
The projection $\pi:\MC_{(h,g)}  \lrg \Def_{(h,g)}$ is a  smooth
morphism of functors.
\end{lem}
\begin{proof}
Let $ \beta: B \lrg A$ be a surjection in $\Art $; we prove that
$$
\MC_{(h,g)}(B) \lrg \Def_{(h,g)}(B)\times_{\Def_{(h,g)}(A)}\MC
_{(h,g)}(A)
$$
is surjective.

Let $((x,y,e^p), (l,n,e^m)) \in
\Def_{(h,g)}(B)\times_{\Def_{(h,g)}(A)}\MC _{(h,g)}(A)$; that is,
the class of $(l,n,m)$ and $\beta  (x,y,e^p)$ are the same element
in $\Def_{(h,g)}(A)$.

Then there exists $(a,b) \in \exp(L^0 \otimes m_A) \times \exp(N^0
\otimes m_A)$ such that
$$
(l,n,e^m)= (e^a,e^b)*(\beta(x,y,e^p))=(e^a * \beta(x),e^b *
\beta(y), e^{g(b)} e^{\beta(p)}e^{-h(a)}).
$$
Let $c \in L^0 \otimes m_B$ be a lifting of $a$ and $d \in N^0
\otimes m_B$ be a lifting of $b$. Then
$$
t= (e^c*x,e^d*y,e^{g(d)} e^{ p }e^{-h(c)})
$$
lies in $\MC_{(h,g)}(B) $ and   is a lifting of $((x,y,e^p),
(l,n,e^m))$.

Actually, $t$ is gauge equivalent to $(x,y,e^p)$ and
$$
\beta(t)=(e^{\beta(c)}*\beta(x),e^{\beta(d)}* \beta(y),
e^{g(\beta(d))}e^{\beta(p)}e^{-h(\beta(c))})=
$$
$$
(e^a * \beta(x), e^b*\beta(y), e^{g(b)} e^{\beta(p)}e^{-h(a)}
)=(l,n,e^m).
$$
\end{proof}

\subsection{Tangent and obstruction spaces of $\MC_{(h,g)}$
and $\Def_{(h,g)}$}\label{sez tangent e ostru MC(h,g) DEF(h,g)}

By Definition~\ref{def funtore MC_(h,g) non estesi}, \emph{the
tangent space of $\MC_{(h,g)}$ } is
$$
\MC_{(h,g)}(\kepsi)=
$$
$$
=\{(x,y,e^p) \in (L^1 \otimes \K\epsi) \times  (N^1 \otimes
\K\epsi) \times  \operatorname{exp}(M^0 \otimes \K\epsi)|
$$
$$
dx=dy=0, h(x)-g(y)-dp=0 \}
$$
$$
\cong \{(x,y, p) \in L^1 \times  N^1  \times
 M^0  |\ dx=dy=0, g(y)=h(x)-dp \}=
$$
$$
\ker (D:\cil^1 \lrg \cil^2).
$$

By Definition~\ref{defin funtor Def_(h,g) non esteso}, \emph{the
tangent space of $\Def_{(h,g)}$ } is
$$
\Def_{(h,g)}(\kepsi)=
$$
$$
\frac{\{(x,y, p) \in L^1 \times  N^1  \times
 M^0  |\ dx=dy=0, g(y)=h(x)-dp \} }{
\{ (-da,-db,g(b)-h(a))|\ a \in (L^0 \otimes \kepsi), b \in (N^0
\otimes \kepsi)  \}}.
$$
$$
\cong H^1(\cil^\cdot ).
$$
\begin{oss}
In the last equality we use the extra hypothesis that $M^{-1}=0$.

\end{oss}

As to the obstruction  space of $\MC_{(h,g)}$, we prove below that
it is naturally contained in $H^2 (\cil^\cdot)$. Since
$\pi:\MC_{(h,g)} \lrg \Def_{(h,g)}$ is smooth, Corollary~\ref{cor
f:F_G liscio implica biezione ostruzio} implies that the
obstruction space of $\Def_{(h,g)}$ is also contained in $H^2
(\cil^\cdot)$.

\begin{lem}\label{lem calcolo ostruzione MC_{(h,g)}=H^2(cil)}
$H^2(\cil^\cdot)$ is a complete obstruction space for
$\MC_{(h,g)}$.
\end{lem}
\begin{proof}
Let
$$
0 \lrg J \lrg B \stackrel{\alpha}{\lrg} A\lrg 0
$$
be a small extension and $(x,y,e^p) \in \MC_{(h,g)}(A)$.

Since $\alpha$ is surjective there exist  $\tilde{x} \in L^1
\otimes m_B $ that lifts $x$, $\tilde{y}\in N^1 \otimes m_B $ that
lifts $y$, and  $q \in M^0 \otimes m_B$ that lifts  $p$. Let
$$
l=d\tilde{x}+ \displaystyle\frac{1}{2}[\tilde{x},\tilde{x}] \in
L^2 \otimes m_B
$$
and
$$
k=d\tilde{y}+ \displaystyle\frac{1}{2}[\tilde{y},\tilde{y}] \in
N^2 \otimes m_B.
$$
As in Lemma~\ref{lem calcolo ostruzione MC_L=H^2(L)}, we can
easily  prove that $\alpha(l)=\alpha(k)=dl=dk=0$; then $l \in
H^2(L) \otimes J$ and $k \in H^2(N) \otimes J$.

Let $r =-g(\tilde{y})+e^q * h(\tilde{x})\in M^1\otimes m_B$; in
particular, $\alpha(r)=0$ and so $r \in M^1\otimes J$.

Next, we prove that $(l,k,r) \in Z^2(\cil ^\cdot)\otimes J$.

Since $dl=dk=0$, it remains to prove that $-dr-g(k)+h(l)=0$.

By definition, $h(\tilde{x})=e^{-q}*(g(\tilde{y})+r)=r +
e^{-q}*g(\tilde{y})$ (in the last equalities we use
Example~\ref{exe azione gauge se []=0}) and so
$$
h(l)=d(h(\tilde{x}))+ \displaystyle\frac{1}{2}
[h(\tilde{x}),h(\tilde{x})]= dr + d(e^{-q}*g(\tilde{y}))+
\displaystyle\frac{1}{2}
[e^{-q}*g(\tilde{y}),e^{-q}*g(\tilde{y})].
$$
Let $A=d(e^{-q}*g(\tilde{y}))$ and $B=
[e^{-q}*g(\tilde{y}),e^{-q}*g(\tilde{y})]$. Therefore,  it is
sufficient to prove  that  $A+ \displaystyle\frac{1}{2} B=g(k)$.

We have
$$
B=[e^{-q}*g(\tilde{y}),e^{-q}*g(\tilde{y})]= [e^{[-q, \,
]'}(d+g(\tilde{y}))-d, e^{[-q, \, ]'}(d+g(\tilde{y}))-d]'=
$$
$$
[e^{[-q, \, ]'}(d+g(\tilde{y})), e^{[-q, \, ]'}(d+g(\tilde{y}))
]'+[e^{[-q, \, ]'}(d+g(\tilde{y})),-d]' +[-d,e^{[-q, \,
]'}(d+g(\tilde{y}))]'=
$$
$$
e^{[-q, \, ]'}[d+g(\tilde{y}),d+g(\tilde{y})]'- 2[d,e^{[-q, \, ]'}
(d+g(\tilde{y}))]'=
$$
$$
e^{[-q, \, ]'}(2dg(\tilde{y})+[g(\tilde{y}),g(\tilde{y})])
-2[d,e^{[-q, \, ]'} (d+g(\tilde{y}))-d]' =2e^{[-q,]}(g(k))-2A.
$$
By assumption, $g(k) \in M^2 \otimes J$ and so $e^{[-q,]}
(g(k))=g(k)$.

Therefore,
$$ A+ \displaystyle\frac{1}{2} B=A +
\displaystyle\frac{1}{2}( 2g(k)-2A)=g(k)
$$
and so $[(l,k,r)] \in H^2(\cil ^\cdot)\otimes J$.

This class does not  depend on the choice of the liftings.

Actually, let $\tilde{x}'\in L^1 \otimes m_B $ be another lifting
of $x$. Then $\tilde{x}'=\tilde{x}+j_x$, for some $j_x \in L^1
\otimes J$ and so
$$
l'=d\tilde{x}'+ \displaystyle\frac{1}{2}[\tilde{x}',\tilde{x}']=
d\tilde{x} +dj_x+ \displaystyle\frac{1}{2}
[\tilde{x},\tilde{x}]=l+dj_x.
$$
Analogously, if $\tilde{y}'$ is another lifting of $y$  then there
exists $j_y \in N^1 \otimes J$ such that
$$
k'=k +dj_y.
$$
Moreover,
$$
r'=-g(\tilde{y}')+e^q*h(\tilde{x}')=-g(\tilde{y})-g(j_y)+
e^q*h(\tilde{x})+h(j_x)=r-g(j_y)+h(j_x).
$$
Therefore, $[(l',k',r')]=[(l+dj_x,k+dj_y,r-g(j_y)+h(j_x))]
=[(l,k,r)] \in H^2(\cil ^\cdot)\otimes J$.

In conclusion, $[(l,k,r)] \in H^2(\cil ^\cdot)\otimes J$ is the
obstruction class associated with the element $(x,y,e^p) \in
\MC_{(h,g)}(A)$.

\smallskip

If this class vanishes, then there exists $(u,v,z) \in
\cil^1\otimes J$ such that $( du,dv,-dz-g(v)+h(u))=(l,k,r)$. In
this case, define $ \overline{x}=\tilde{x}-u$, $
\overline{y}=\tilde{y}-v$ and $ \overline{z}=q-z$. Then $(
\overline{x},\overline{y},e^{\overline{z}}) \in \MC_{(h,g)}(B)$
and it lifts  $ (x,y,e^p)$. Actually,
$$
d \overline{x} +[ \overline{x}, \overline{x}]= d\tilde{x}-du
+[\tilde{x},\tilde{x}]=l-du=0,
$$
$$
d  \overline{y}+[ \overline{y}, \overline{y}]= d\tilde{y}-dv
+[\tilde{y},\tilde{y}]=k-dv=0
$$
and
$$
g(\overline{y})-e^{\overline{z}}*h(\overline{x})=
g(\tilde{y})-g(v) -e^{q-z}*(h(\tilde{x})-h(u))=\footnote{See
Example~\ref{exe azione gauge se []=0}.}
$$
$$
g(\tilde{y})-g(v) - e^{q-z}*h(\tilde{x}) + h(u) =
$$
$$
(g(\tilde{y})- e^{q}*h(\tilde{x}))-dz-g(v)+h(u)=-r+r=0.
$$

\end{proof}

\subsection{Properties}\label{sez propriet di MC(h,g) DEF(hg)}

\begin{lem}\label{lemma L M N abelian ->DEF(h,g)liscio}
Let $h:L \lrg M$ and $g:N \lrg M$ be morphisms of abelian DGLAs.
Then the functor $\MC_{(h,g)}$ is smooth.
\end{lem}
\begin{proof}
We have to prove that for every  surjection $  \varphi : B \lrg A
\in  \Art$ the   map
$$
\mc(B) \lrg \mc(A)
$$
is surjective. By hypothesis,  we have
$$
\mc (A)=\{(x,y,e^p) \in (L^1 \otimes m_A) \times  (N^1 \otimes
m_A) \times \exp(M^0 \otimes m_A)|
$$
$$
dx=dy=0, \qquad  g(y)=e^p*h(x)=h(x)-dp\}.
$$
This implies that the Maurer-Cartan equation reduces to the linear
equations  $dx=dy=0, g(y)+dp-h(x)=0$; then $\mc(A)=Z^1(\cil^\cdot
\otimes m_A)=Z^1(\cil^\cdot)\otimes m_A$.

Since $Z^1(\cil^\cdot \otimes m_A)  \twoheadrightarrow
Z^1(\cil^\cdot \otimes m_B)$, the lifting exists.




\end{proof}

\begin{oss}
Every commutative diagram of morphisms of DGLAs
\begin{center}
$\xymatrix{ & L\ar[d]_h
\ar[r]^{\alpha'} & P\ar[d]^\eta \\
& M\ar[r]^\alpha & Q \\
 N\ar[ur]^g \ar[r]^{\alpha''} & R \ar[ur]_\mu & \\}$
\end{center}
induces a \emph{morphism $\varphi^\cdot$ of complexes }
$$
\cil^i \ni (l,n,m) \stackrel{\varphi ^i}{\longmapsto} (\alpha'(l),
\alpha''(n),\alpha(m)) \in \operatorname{C}^i_{(\eta,\mu)}
$$
and a natural \emph{transformation $F$ }of the associated
deformation functors:
$$
F: \Def_{(h,g)} \lrg \Def_{(\eta,\mu)}.
$$

\end{oss}

Then, we obtain the following proposition that is a generalization
of \cite[Proposition~2.3]{bib ManettiPREPRINT}.
\begin{prop}\label{prop DEF_>DEF liscio,ostruz in ker tra H^2 }
Let
\begin{center}
$\xymatrix{ & L\ar[d]_h
\ar[r]^{\alpha'} & P\ar[d]^\eta \\
& M\ar[r]^\alpha & Q \\
N\ar[ur]^g \ar[r]^{\alpha''} & R \ar[ur]_\mu & \\}$
\end{center}
be a commutative diagram of differential graded Lie algebras. If
the functor $\Def_{(\eta,\mu)}$ is smooth, then the obstruction
space of $\, \Def_{(h,g)}$ is contained in the kernel of the map
$$
H^2(\cil^\cdot) \lrg H^2(\operatorname{C}_{(\eta,\mu)}^\cdot).
$$
\end{prop}
\begin{proof}
The morphism $F:\Def_{(h,g)} \lrg \Def_{(\eta,\mu)}$ induces a
linear map between obstruction spaces. If $\Def_{(\eta,\mu)}$ is
smooth, then its obstruction space is zero (Proposition~\ref{prop
f liscio se e solo se OF=0}).

\end{proof}

\begin{teo}
\label{teo NO exte quasi iso C(h,g)-C(n,m)then DEf iso DEF} If
$\varphi^\cdot:\cil^\cdot \lrg
\operatorname{C}^\cdot_{(\eta,\mu)}$ is a quasi-isomorphism of
complexes then $F:\Def_{(h,g)} \lrg \Def_{(\eta,\mu)}$ is an
isomorphism of functors.
\end{teo}

The proof of this theorem is postponed to the end of
Section~\ref{sezio definizione EXTENDED MC_(f,g)}.

\bigskip

Next, let
\begin{center}
$\xymatrix{ & H\ar[d]_\beta
\ar[r]^{\alpha} & L\ar[d]^h \\
& N\ar[r]^g & M  \\}$
\end{center}
be a commutative diagram of morphisms of DGLAs. Thus, it induces a
morphism of complexes
$$
H \stackrel{\varphi^\cdot}{\lrg} \cil ^\cdot \qquad \mbox{ with }
\qquad \varphi(x)=(\alpha(x),\beta(x),0),
$$
and a morphism of functors
$$
\Def_H \stackrel{F}{\lrg} \Def_{(h,g)} \qquad \mbox{ with } \qquad
F(x)= (\alpha(x),\beta(x),0).
$$

\begin{teo}\label{teo NO EXT H prodotto fibrato Q.I then ISO fun}
With the notation above, if the morphism $\varphi^\cdot$ is a
quasi-isomorphism then $F$ is an isomorphism of functors.
\end{teo}
The proof of this theorem too is postponed to the end of
Section~\ref{sezio definizione EXTENDED MC_(f,g)}.

\section{Extended Deformation Functors}\label{sezio extended defo
funtore}

In this section we study the extended deformation functors. In
particular, we are interested in the  functors
$\widetilde{\MC}_{(h,g)}$ and $ \widetilde{\Def}_{(h,g)} $ which
are a generalization of the functors ${\MC}_{(h,g)}$  and
${\Def}_{(h,g)} $ introduced in Section~\ref{sezio non ext
def(h,g) e MC(h,g)}. Here, we remove the restrictive hypothesis of
$M$ concentrated in non negative degrees.

The main references for this chapter  are \cite{bib Manetti IMNR},
\cite[Sections~5.7 and~5.8]{bib manRENDICONTi} and
\cite[Sections~6 and~7]{bib ManettiPREPRINT}.

\subsection{Notation}

We denote by:

$\Com$ the category of   nilpotent (associative and commutative)
differential graded algebras which are finite dimensional as $
\mathbb{K}$-vector spaces.

$\Com_0$ the full subcategory of $\Com$ of algebras with trivial
multiplication.

\begin{exe}\label{exe defh omega[n] tngent estesi}
Define the complex $(\Omega=\Omega_0 \oplus \Omega_1,d)$, where
$\Omega_0=\K$, $\Omega_1=\K[-1]$ and $d:\Omega_0 \lrg \Omega_1$ is
the canonical linear isomorphism $d(1[0])=1[-1]$. $\Omega \in
\Com_0$ and the projection $p: \Omega \lrg \Omega_0=\K$ and the
inclusion $\Omega_1 \lrg \Omega$ are morphisms in $\Com_0$.

Moreover, $\Omega[n]=\K[n]\otimes \Omega$ is an acyclic complex in
$\Com_0$, for each $n \in \N$.

\end{exe}

\subsection{Definition of extended functors}

Let $A \in \Com$ and $J \subset A$ a differential ideal; then $J
\in \Com$ and the inclusion $J \lrg A$ is a morphism of
differential graded algebras.

\begin{dhef}\label{defin small extension}
A \emph{small extension} in $\Com$ is a short exact sequence
$$
0 \lrg J \lrg B \stackrel{\alpha}{\lrg} A \lrg 0
$$
such that $\alpha$ is a morphism in $\Com$ and $J$ is an ideal of
$B$ such that $BJ=0$; in addition, it is called \emph{acyclic} if
$J$ is an acyclic complex, or equivalently $\alpha$ is a
quasi-isomorphism.
\end{dhef}

\begin{dhef}\label{definizione funct predeformazio}
A covariant functor $F: \Com \lrg \Set$ is called a
\emph{predeformation functor} if the following conditions are
satisfied:

\ref{definizione funct predeformazio}.1) $F(0)=\{*\}$ is the one
point set.

\ref{definizione funct predeformazio}.2) For every $A,B \in \Com$,
the natural map
$$
F(A \times B) \lrg F(A) \times F(B)
$$
is bijective.

\ref{definizione funct predeformazio}.3) For every surjective
morphism $\alpha: A \lrg C$ in $\Com$, with $C $ an acyclic
complex in $ \Com_0$, the natural morphism
$$
F(\ker (\alpha))=F(A\times_C 0) \lrg F(A)\times_{F(C)} F(0)=F(A)
$$
is bijective.

\ref{definizione funct predeformazio}.4) For any pair of morphisms
$\beta: B \lrg A$ and $\gamma: C \lrg A$  in $\Com$, with $\beta$
surjective, the natural map
$$
F(B \times_A C) \lrg F(B) \times_{F(A)} F(C)
$$
is surjective.

\ref{definizione funct predeformazio}.5) For every small acyclic
extension
$$
0 \lrg J \lrg B \stackrel{\alpha}{\lrg}A \lrg 0
$$
the induced map $F(B) \lrg F(A)$ is surjective.
\end{dhef}

\begin{dhef}\label{dhef extended defo functor}
A covariant functor $F: \Com \lrg \Set$ is called a
\emph{deformation functor} if it is a predeformation functor and
$F(J)=0$ for every acyclic complex $J \in \Com_0$.
\end{dhef}

\subsubsection{Examples}

Let $L$ be a differential graded Lie algebra and $A \in \Com$.
Then $L \otimes A$ has a natural structure of (nilpotent) DGLA
given by:
$$
(L \otimes A)^i=\bigoplus_{j \in \Z} L^j \otimes A^{i-j};
$$
$$
d(x \otimes a)=dx\otimes a +(-1)^{\deg(x)}x \otimes da;
$$
$$
[x\otimes a, y \otimes b]=(-1)^{\deg(a)\deg(y)}[x,y]\otimes ab.
$$

\begin{dhef}
The \emph{extended Maurer-Cartan} functor associated with a DGLA
$L$ is
$$
\widetilde{\MC}_L (A)=\{x \in  (L \otimes  A)^1 | \, dx
+\frac{1}{2}[x,x]=0 \}.
$$
\end{dhef}

\begin{lem}\label{lem EXT MC_L is predefo}
$\widetilde{\MC}_L$ is a predeformation functor.
\end{lem}
\begin{proof}
See \cite[Lemma~2.15]{bib Manetti IMNR}. We will also give a proof
in Section~\ref{sezio definizione EXTENDED MC_(f,g)},   since it
is a particular case of Theorem~\ref{teo mc funtor predefor}.
\end{proof}

\begin{oss}\label{oss EXTE MC(C)-Z^1(LxC)}
We note that, for each $C \in \Com_0$, we have:
$$
\widetilde{\MC}_L(C)=\{ x \in  (L \otimes  C)^1 | \, dx=0 \} =
Z^1(L \otimes C).
$$

\end{oss}

\begin{dhef}
The \emph{extended deformation} functor associated with a DGLA $L$
is
$$
\widetilde{\Def}_L (A)=\frac{\{ x \in  (L \otimes A)^1 | \, dx +
\frac{1}{2}[x,x]=0\}} {\mbox{ gauge action of } \exp\, (L\otimes
A)^0},
$$
where the gauge action of $ a \in (L\otimes A)^0$ is the analogous
of the not extended case, i.e.,
$$
e^a * x:= x+\sum_{n=0}^{\infty} \frac{ [a,-]^n}{(n+1)!}([a,x]-da).
$$
\end{dhef}

\begin{lem}\label{lemma EXT DEF_L is deformation funt}
$\widetilde{\Def}_L$ is a deformation functor.
\end{lem}
\begin{proof}
See \cite[Theorem~2.16]{bib Manetti IMNR}. We will also give a
proof in Section~\ref{sezio definizione EXTENDED MC_(f,g)}, since
it is a particular case of Theorem~\ref{teo EXT DEF_(h,g) is
deformation}.

\end{proof}

\begin{oss}
For each $C \in \Com_0$ we note that:
$$
\widetilde{\Def}_L(C)=\frac{ \{x \in  (L \otimes  C)^1 | \, dx=0
\}}{\{da| \, a \in  (L \otimes  C)^0 \} } = H^1(L \otimes C).
$$

\end{oss}

\subsection{Properties}

As in the not extended case, for every predeformation functor $F$
and every $C \in \Com_0$, there exists a natural structure of
vector space on $F(A)$ with:
\begin{itemize}
  \item[-] addition  given by the map $F(C \times C)=F(C)\times
  F(C)\stackrel{+}{\lrg}F(C)$ induced by
  $C \times C \stackrel{+}{\lrg}C$;
  \item[-] scalar multiplication by $ s$ given by the map
  $F(C)\stackrel{\cdot s}{\lrg}F(C)$ induced by
  $C\stackrel{\cdot s}{\lrg}C$.
\end{itemize}

\begin{oss}
$1)$ For any morphism $B \lrg A$ in $\Com_0$, the induced map
$F(B)\lrg F(A)$ is $\C$-linear.

$2)$ For any natural transformation of predeformation functors $F
\lrg G$,   the induced  map $F(C)\lrg G(C)$ is $\C$-linear for
each $C \in \Com_0$.
\end{oss}

\begin{dhef}\label{def tangente funtor predefor}
Let $F$ be a predeformation functor, the \emph{tangent space} of
$F$ is the graded vector space $TF[1]$, where
$$
TF=\bigoplus_{n \in \Z}T^nF, \ \ \ T^{n+1}F=TF[1]^n=\mbox{
coker}(F(\Omega[n]) \stackrel{p}{\lrg} F(\mathbb{K}[n])), n \in \Z
$$
and $p$ is the linear map induced by the projection $\Omega[n]\lrg
\K[n]$ (see  Example~\ref{exe defh omega[n] tngent estesi}).
\end{dhef}

In particular, if $F$ is a deformation functor then
$F(\Omega[n])=0$ for every $n$. Therefore, $T^{n+1}F=TF[1]^n=F(
\kepsi)$, where $\epsi$ is an indeterminate of degree $-n \in \Z$,
such that $\epsi^2=0$.

\begin{dhef}
A natural transformation of predeformation functors $F \lrg G$  is
called a \emph{quasi-isomorphism} if it induces isomorphisms on
tangent spaces, i.e., $T^nF \cong T^n G$.
\end{dhef}

\begin{teo}\label{TEO funzione inversa}
(Inverse function theorem) A natural transformation of deformation
functors is an isomorphism if and only if it is a
quasi-isomorphism.
\end{teo}
\begin{proof}
See \cite[Corollary~3.2]{bib Manetti IMNR} or
\cite[Corollary~5.72]{bib manRENDICONTi}.
\end{proof}

\begin{teo}(Manetti)\label{teo esistenza funtor defor F+}
Let $F$ be a predeformation functor, then there exist a
deformation functor $F^+$ and a natural transformation $\nu:F \lrg
F^+$, that is a quasi-isomorphism, such that for every deformation
functor $G$ and every natural transformation $\phi: F \lrg G$
there exists a unique natural transformation $\psi: F^+ \lrg G$
such that $\phi= \psi \nu$.
\end{teo}

\begin{proof}
See \cite[Theorem~2.8]{bib Manetti IMNR}.
\end{proof}

\begin{oss}
Given a natural transformation of predeformation functors
$\alpha:F \lrg G$, there is a natural transformation of associated
deformation functors $\alpha^+:F^+ \lrg G^+$. Actually, let
$\eta:G \lrg G^+$. Then by composition we get $\beta=\eta \circ
\alpha:F \lrg G^+$ and so, by   Theorem~\ref{teo esistenza funtor
defor F+}, there exists $\alpha^+:F^+ \lrg G^+$.

\end{oss}

\begin{teo}(Manetti)\label{lemma EXTEND spazio tangen(FUNTOR)=
complesso} Let $S$ be a complex of vector spaces and assume that
the functor
$$
\Com_0 \lrg \Set, \ \ \  \ \
$$
$$
\ \ \  C \longmapsto Z^1(S \otimes C)
$$
is the restriction of a predeformation functor $F$. Then for every
complex $C \in \Com_0$, the equality $F^+(C)=H^1(S \otimes C)$
holds; in particular, $T^iF^+=H^i(S)$.
\end{teo}

\begin{proof}
See \cite[Lemma~2.10]{bib Manetti IMNR}.
\end{proof}

\begin{cor}(Manetti)
For every differential graded Lie algebra $L$ there exists a
natural isomorphism $\widetilde{\MC}^+_L \cong
\widetilde{\Def}_L$.
\end{cor}
\begin{proof}
By Remark~\ref{oss EXTE MC(C)-Z^1(LxC)} and Theorem~\ref{lemma
EXTEND spazio tangen(FUNTOR)= complesso}, we have that
$T^i\widetilde{\MC}^+_L=H^i(L)$. Therefore, the natural projection
$\widetilde{\MC}_L \lrg \widetilde{\Def}_L$ induces (by
Theorem~\ref{teo esistenza funtor defor F+}) a natural
transformation ${\widetilde{\MC}^+}_L \lrg \widetilde{\Def}_L$
which is an isomorphism on tangent spaces.
\end{proof}

\subsection{Extended deformation functor of a pair of morphisms}
\label{sezio definizione EXTENDED MC_(f,g)}

\bigskip

Let $L,M,N$ be DGLAs, and $h :L \lrg M$ and $g: N \lrg M$ be
morphisms of DGLAs:

\begin{center}
$\xymatrix{ & L \ar[d]^h \\
           N \ar[r]_g & M.   \\ }$
\end{center}

\begin{dhef}
The \emph{extended Maurer-Cartan functor associated with the pair
$(h,g)$} is
$$
\widetilde{\MC}_{(h,g)}:\Art \lrg \Set,
$$
$$
\widetilde{\MC}_{(h,g)}(A)=\{(x,y,e^p) \in (L \otimes A)^1\times
(N \otimes A)^1\times \operatorname{exp}(M\otimes A)^0  |
$$
$$
dx +\frac{1}{2}[x,x]=0,\   dy +\frac{1}{2}[y,y]=0, \ g(y)=e^p*h(x)
\}.
$$
\end{dhef}
\begin{teo}\label{teo mc funtor predefor}
$ \widetilde{\MC}_{(h,g)}$ is a predeformation functor.
\end{teo}

\begin{proof}
$\widetilde{\MC}_{(h,g)}(0)=0$ and so (\ref{definizione funct
predeformazio}.1) is satisfied.

\noindent For any pair of morphisms $\beta: B \lrg A$ and $\gamma:
C \lrg A$  in $\Com$,  we have
$$
\widetilde{\MC}_{(h,g)} (B \times_A C)=
\widetilde{\MC}_{(h,g)}(B)\times_{\widetilde{\MC}_{(h,g)}(A)}
\widetilde{\MC}_{(h,g)}(C)
$$
and so $\widetilde{\MC}_{(h,g)}$  satisfies properties
(\ref{definizione funct predeformazio}.2),(\ref{definizione funct
predeformazio}.3) and (\ref{definizione funct predeformazio}.4).

\noindent Let $ 0 \lrg J \lrg B \stackrel{\alpha}{\lrg} A \lrg 0 $
be an acyclic small extension and $(x,y,e^p) \in
\widetilde{\MC}_{(h,g)}(A)$. We want to prove the existence of a
lifting $(\overline{x},\overline{y}, e^{ \overline{p}})\in
\widetilde{\MC}_{(h,g)} (B)$ such that the induced map
$\widetilde{\MC}_{(h,g)}(B) \lrg \widetilde{\MC}_{(h,g)}(A)$ is
surjective.

\noindent Since $\alpha$ is surjective, there exists $(r,s,e^t)
\in (L \otimes B)^1\times  (N \otimes B)^1\times
\operatorname{exp} (M\otimes B)^0  $ such that
$\alpha(r)=x,\alpha(s)=y$ and $\alpha(t)=p$.

Let $l \in (L \otimes J)^2$ and $k \in (N \otimes J)^2$ be defined
as follows
$$
l=dr +\frac{1}{2}[r,r] \ \ \mbox{ and } \ \ \  k=ds +
\displaystyle\frac{1}{2}[s,s].
$$
Then
$$
dl= \displaystyle\frac{1}{2} d[r,r]=[dr,r]=[l,r]-
\displaystyle\frac{1}{2}[[r,r],r]
$$
and the same holds for $k$, i.e.,
$$
dk=[k,s]- \displaystyle\frac{1}{2}[[s,s],s].
$$
Since $B\cdot J=0$, we have $[l,r]=[k,s]=0$; moreover, using
Jacobi's identity (see Remark~\ref{oss [a,a]=0 pari [[a,a],a]=0
dispari}), $[[r,r],r]=[[s,s],s]=0$. This implies that $dl=dk=0$.

By hypothesis, $J$ is acyclic and so, by the K\"{u}nneth formula,
the complexes $L\otimes J$ and $N \otimes J$ are acyclic. Thus,
there exist $w \in (L\otimes J)^1$ and  $z \in (N\otimes J)^1$,
such that $dw=l$ and $dz=k$.

Let
$$
  \overline{x}= r-w \in (L\otimes B)^1 \qquad \mbox{ and }
\qquad \overline{y}=s-z \in (N\otimes B)^1;
$$
we have
$$
\alpha(\overline{x})= \alpha (r)-\alpha(w)=\alpha(r)=x,\quad
 \qquad \alpha(\overline{y})=\alpha
(s)-\alpha(z)=\alpha(s)=y,
$$
$$
d\overline{x}+
\displaystyle\frac{1}{2}[\overline{x},\overline{x}]=dr-dw +
\displaystyle\frac{1}{2}[r,r]=l-l=0 \qquad \mbox{ and } \qquad
d\overline{y}+
\displaystyle\frac{1}{2}[\overline{y},\overline{y}]=0. \qquad
\qquad
$$

Therefore, $\overline{x}$ and $\overline{y}$ lift $x$ and $y$ ,
respectively, and they satisfy the Maurer-Cartan equation.

Let $z=e^t* h(\overline{x})-g(\overline{y}) \in (M \otimes B)^1$.
Since $\alpha(z)=e^{\alpha(t)}*h(\alpha(\overline{x}))
-g(\alpha(\overline{y}))=e^p*h(x)-g(y)=0$,   $z \in (M \otimes
J)^1$. Moreover, $dz=0$; indeed
$$
2dz=2d(e^t*h(\overline{x}))-2d(g(\overline{y}))=-[e^t*
h(\overline{x}),e^t*
h(\overline{x})]+[g(\overline{y}),g(\overline{y})]=
$$
$$
-[e^t* h(\overline{x}),e^t* h(\overline{x})]+
[g(\overline{y}),e^t* h(\overline{x})]- [g(\overline{y}),e^t*
h(\overline{x})] +[g(\overline{y}),g(\overline{y})]=
$$
$$
-[e^t* h(\overline{x})-g(\overline{y}),e^t* h(\overline{x})]-
[g(\overline{y}),e^t* h(\overline{x})-g(\overline{y})]=0,
$$
($z=e^t* h(\overline{x})-g(\overline{y}) \in (M \otimes J)^1$ and
so $[e^t* h(\overline{x})-g(\overline{y}), -]=0$).

Since $M\otimes J$ is acyclic, there exists $v \in(M\otimes J)^0$
such that $z=dv$. Therefore, $$ e^t* h(\overline{x})=z
+g(\overline{y})=dv +g(\overline{y})=e^{-v}*g(\overline{y})
$$
and so
$$
e^v e^t* h(\overline{x})=g(\overline{y}).
$$
This implies that $e^{\overline{p}}=e^v e^t$ lifts $e^p$.

Then the triple $(\overline{x},\overline{y}, e^{ \overline{p}})\in
\widetilde{\MC}_{(h,g)}  (B)$ lifts $(x,y,e^p) \in
\widetilde{\MC}_{(h,g)}(A)$ and so (\ref{definizione funct
predeformazio}.5) holds.

\end{proof}

\begin{proof}[Proof of Lemma~\ref{lem EXT MC_L is predefo}]
It  suffices to apply the previous theorem with $M=N=0$.
\end{proof}

\bigskip

\begin{oss}
If the DGLA $M$ is concentrated in non negative degree, then for
every $(A,m_A) \in \Art$,
$\widetilde{\MC}_{(h,g)}(m_A)={\MC}_{(h,g)}(A)$.
\end{oss}

Applying Theorem~\ref{teo esistenza funtor defor F+}, we can
conclude the existence of a deformation functor
$\widetilde{\MC}_{(h,g)}^+$ associated with
$\widetilde{\MC}_{(h,g)}$.

\begin{prop}\label{tangente di MC+ = h^i triplo cilindro}
$T^i\widetilde{\MC}_{(h,g)}^+\cong H^i(\cil^\cdot)$.
\end{prop}

\begin{proof}
For any $C \in\Com_0$ we have
$$
\widetilde{\MC}_{(h,g)}(C)=\{ (l,n,e^m)\in (L\otimes C)^1 \times
(N\otimes C)^1 \times \operatorname{exp}(M\otimes C)^0|
$$
$$
 d l=d n=0, \
g(n)=e^m*h(l)=h(l)-d m \},
$$
and
$$
Z^1(\cil^\cdot \otimes C)=\{(l,n,m)\in (L\otimes C)^1 \times
(N\otimes C)^1 \times (M\otimes C)^0|
$$
$$
\  dl= dn=  -dm-g(n)+h(l)=0 \}.
$$
Therefore, $\widetilde{\MC}_{(h,g)}(-)_{| \Com_0}=Z^1(\cil^\cdot
\otimes -)_{|\Com_0}$. Then we can apply Theorem~\ref{lemma EXTEND
spazio tangen(FUNTOR)= complesso} to complete the proof.
\end{proof}

\bigskip

Next, we consider on $\widetilde{\MC}_{(h,g)}(A)$ the following
equivalence relation $ \approx$:
$$
(x_1,y_1,e^{p_1})\, \approx\,(x_2,y_2,e^{p_2})
$$
if and only if there exist $a \in (L \otimes A)^0$, $b \in (N
\otimes A)^0$ and $c \in (M\otimes A)^{-1} $ such that
$$
x_2=e^a *x_1, \qquad y_2=e^b*y_1
$$
and
$$
e^{p_2}=e^{g(b)}e^T e^{p_1} e^{-h(a)}, \quad \mbox{with }\quad
T=dc+[g(y_1),c]\in Stab_A(g(y_1)).
$$

\begin{lem}
The relation  $ \approx$ is an equivalence relation.
\end{lem}
\begin{proof}
The reflexivity is obvious. As to the symmetry and the
transitivity, we use the following property of the irrelevant
stabilizers $Stab_A(-)$(see Section~\ref{sezio gauge action}): for
each $x\in \widetilde{MC}_M(A)$, $a \in (M\otimes A)^0 $ and $T
=dc+[x,c]\in Stab_A(x)$ ($c \in (M\otimes A)^{-1}$), there exists
$f \in (M\otimes A)^{-1}$ such that
$$
e^a e^{T }=e^{T'}e^a
$$
where $T'=df+[y,f] \in Stab_A(y)$ and $y=e^a*x $.

$Symmetry.$ Let $(x_1,y_1,e^{p_1})\, \approx\,(x_2,y_2,e^{p_2})$.
Therefore, there exist $a \in (L \otimes A)^0$, $b \in (N \otimes
A)^0$ and $c \in (M\otimes A)^{-1} $ such that
$$
x_2=e^a *x_1, \qquad y_2=e^b*y_1
$$
and
$$
e^{p_2}=e^{g(b)}e^T e^{p_1} e^{-h(a)}, \quad \mbox{with }\quad
T=dc+[g(y_1),c].
$$
This implies the existence of $f \in  (M\otimes A)^{-1}$ such that
$$
e^{p_2}=e^{T'}e^{g(b)} e^{p_1} e^{-h(a)} \quad \mbox{with }\quad
T'=df+[g(y_2),f].
$$
Thus, by choosing $\alpha=-a \in (L \otimes A)^0$, $\beta=-b \in
(N \otimes A)^0$ and $\gamma =-f \in  (M\otimes A)^{-1}$ yields
$$
x_1=e^{\alpha} *x_2, \qquad y_1=e^{\beta}*y_2
$$
and
$$
e^{p_1}=e^{-g(b)}e^{-T'} e^{p_2}e^{h(a)}=   e^{g(\beta)}e^{T''}
e^{p_2} e^{-h(\alpha)} \quad \mbox{with }\quad
T''=d\gamma+[g(y_1),\gamma].
$$
Then $(x_2,y_2,e^{p_2})\, \approx\, (x_1,y_1,e^{p_1}).$

$Transitivity.$ Suppose
$$
(x_1,y_1,e^{p_1})\, \approx\,(x_2,y_2,e^{p_2}) \quad \mbox{and
}\quad (x_2,y_2,e^{p_2})\, \approx\, (x_3,y_3,e^{p_3}).
$$
Therefore, there exist $a_1,a_2 \in (L \otimes A)^0$, $b_1,b_2 \in
(N \otimes A)^0$ and $c_1,c_2 \in (M\otimes A)^{-1} $ such that
$$
x_2=e^{a_1} *x_1, \qquad y_2=e^{b_1}*y_1\qquad
e^{p_2}=e^{g(b_1)}e^{T_1} e^{p_1} e^{-h(a_1)}
$$
and
$$
x_3=e^{a_2} *x_2, \qquad y_3=e^{b_2}*y_2\qquad
e^{p_3}=e^{g(b_2)}e^{T_2} e^{p_2} e^{-h(a_2)}
$$
with $ T_1=dc_1+[g(y_1),c_1]$ and $ T_2=dc_2+[g(y_2),c_2]$. Then
$$
e^{p_3}=e^{g(b_2)}e^{T_2} e^{g(b_1)}e^{T_1} e^{p_1} e^{-h(a_1)}
e^{-h(a_2)}=
$$
$$
e^{g(b_2)}  e^{g(b_1)}e^{{T_2}'}e^{T_1} e^{p_1} e^{-h(a_1)}
e^{-h(a_2)},
$$
for some $c'\in (M\otimes A)^{-1}$ and ${T_2}'=dc'+[g(y_1),c']$.

Since $Stab_A(g(y_1))$ is a subgroup, there exists $c   \in
(M\otimes A)^{-1}$ such that $e^{{T_2}'}e^{T_1}=e^T$ with
$T=dc+[g(y_1),c]$.

Let $a=a_2 \bullet a_1   \in (L \otimes A)^0$, $b=b_2\bullet b_1
\in (N \otimes A)^0$ and $c  \in (M\otimes A)^{-1}$ as above. Then
$$
x_3=e^{a } *x_1, \qquad y_3=e^{b }*y_3,\qquad e^{p_3}=e^{g(b
)}e^{T } e^{p_1} e^{-h(a )}
$$
and so $(x_1,y_1,e^{p_1})\,   \approx\, (x_3,y_3,e^{p_3}).$

\end{proof}

\begin{oss}
We note that this equivalence relation generalizes the equivalence
relation induced by the gauge action given in Definition~\ref{def
gauge action MC_(h,g) NO EXT}, when $M$ is concentrated in non
negative degree.
\end{oss}

\begin{dhef}
Define the functor
$$
\widetilde{\Def}_{(h,g)}:\Com \lrg \Set,
$$
$$
\widetilde{\Def}_{(h,g)}(A)=\widetilde{\MC}_{(h,g)}(A)/\approx.
$$
\end{dhef}
\begin{teo}\label{teo EXT DEF_(h,g) is deformation}
$\widetilde{\Def}_{(h,g)}:\Com \lrg \Set$ is a deformation functor
with $T^i \widetilde{\Def}_{(h,g)} =H^i(C^\cdot _{(h,g)})$.
\end{teo}

\begin{proof}
We first prove that $\widetilde{\Def}_{(h,g)}$ is a predeformation
functor.

Since $\widetilde{\Def}_{(h,g)}$ is the quotient of the
predeformation functor $\widetilde{\MC}_{(h,g)}$,  conditions
(\ref{definizione funct predeformazio}.1) and (\ref{definizione
funct predeformazio}.5) are verified.

An easy calculation shows that $\widetilde{\Def}_{(h,g)}$
satisfies (\ref{definizione funct predeformazio}.2).

Next, we verify  condition (\ref{definizione funct
predeformazio}.4). Let $\beta: B \lrg A$ and $\gamma: C \lrg A$ be
morphisms  in $\Com$, with $\beta$ surjective. We prove that the
natural map $ \widetilde{\Def}_{(h,g)}(B \times_A C) \lrg
\widetilde{\Def}_{(h,g)}(B)  \times_{\widetilde{\Def}_{(h,g)}(A)}
\widetilde{\Def}_{(h,g)}(C) $ is surjective.

Let $(x_1,y_1,e^{p_1})\in \widetilde{\MC}_{(h,g)}(B)$ and $
(x_2,y_2,e^{p_2}) \in \widetilde{\MC}_{(h,g)}(C)$, such that
$\beta(x_1,y_1,e^{p_1})$ and $\gamma (x_2,y_2,e^{p_2})$ are the
same element in $\widetilde{\Def}_{(h,g)}(A)$.

Then, there exist $a \in (L \otimes A)^0$, $b \in (N \otimes A)^0$
and  $c \in (M\otimes A)^{-1} $ such that
$$
\gamma(x_2)=e^a*\beta(x_1), \qquad \gamma(y_2)=e^b*\beta(y_1)
$$
and
$$
e^{\gamma(p_2)}=e^{g(b)}e^T e^{\beta(p_1)} e^{-h(a)}
 \quad \mbox{with } T=dc+[g(\beta(y_1)),c].
$$

Let $\tilde{a}\in (L \otimes B)^0$ be a lifting of $a$,
$\tilde{b}\in (N \otimes B)^0$ a lifting of $b$ and $\tilde{c} \in
(M \otimes B)^{-1}$ a lifting of $c$, so that
$\tilde{T}=d\tilde{c} +[g( y_1 ),\tilde{c}] $ lifts $T$.

Up to substituting $(x_1,y_1,e^{p_1})$ with its  gauge equivalent
$(e^{\tilde{a}}*x_1,e^{\tilde{b}}*y_1,
e^{g(\tilde{b})}e^{\tilde{T}}e^{p_1} e^{h(\tilde{a})})$, we can
suppose that\footnote{$ \beta(e^{\tilde{a}}*x_1,e^{\tilde{b}}*y_1,
e^{g(b)}e^{\tilde{T}}e^{p_1} e^{h(\tilde{a})})=
(e^{\beta(\tilde{a})}*\beta(x_1),e^{\beta(\tilde{b})}*\beta(y_1),
e^{g(\beta(\tilde{b}))} e^{\beta(\tilde{T})}e^{\beta(p_1)}
e^{h(\beta(\tilde{a}))} )=(e^a*\beta(x_1),e^b
*\beta(y_1),e^{g(b)}e^T e^{\beta(p_1)}e^{-h(a)})=(
\gamma(x_2),\gamma(y_2), e^{\gamma(p_2)})= \gamma
(x_2,y_2,e^{p_2}) $} $\beta(x_1,y_1,e^{p_1})=\gamma
(x_2,y_2,e^{p_2}) \in \widetilde{\MC}_{(h,g)}(A)$.

Then $((x_1,y_1,e^{p_1}), (x_2,y_2,e^{p_2}))\in
\widetilde{\MC}_{(h,g)}(B) \times_{\widetilde{\MC}_{(h,g)}(A)}
\widetilde{\MC}_{(h,g)}(C) $ and so, since
$\widetilde{\MC}_{(h,g)}$ is a predeformation functor, there
exists a lifting   in $ \widetilde{\MC}_{(h,g)}(B\times_A C)$.
Next,  it is sufficient to take its equivalence class in
$\widetilde{\Def}_{(h,g)}$.

Finally, we prove that condition (\ref{definizione funct
predeformazio}.3) is satisfied.

Let $\alpha:A \lrg C$ be a surjection with $C \in \Com_0$ an
acyclic complex. Let $(x_1,y_1,e^{p_1})$ and $(x_2,y_2,e^{p_2})
\in \widetilde{\Def}_{(h,g)}(\ker \alpha)$ (in particular,
$g(y_1)= e^{p_1}*h(x_1)$ and $g(y_2)=e^{p_2}*h(x_2)$) be gauge
equivalent as elements in $\widetilde{\MC}_{(h,g)}(A)$, i.e.,
there exist $a \in (L\otimes A)^0$, $b \in (N\otimes A)^0$ and $c
\in (M\otimes A)^{-1} $ such that
$$
 x_2=e^a* x_1,   \qquad y_2=e^b* y_1,
$$
$$
e^{p_2}=e^{g(b)}e^T e^{p_1} e^{-h(a)}, \quad \mbox{with } \qquad
T=dc+[g(y_1),c].
$$
We are looking for $\overline{a} \in  (L\otimes \ker \alpha)^0$,
$\overline{b} \in (N\otimes \ker \alpha)^0$ and $\overline{c} \in
(M\otimes \ker \alpha)^{-1}$ such that
$\overline{T}=d\overline{c}+[g(y_1),\overline{c}]$ and
$$
 x_2 =e^{\overline{a}}* x_1 m  \qquad y_2=e^{\overline{b}}* y_1,
\qquad e^{p_2}= e^{g(\overline{b})}e^{\overline{T} }e^{ p_1}
e^{-h(\overline{a})}.
$$
Since $L\otimes C$ and $N\otimes C$ are abelian DGLAs and
$\alpha(x_i)=\alpha(y_i)=0$, for $i=1,2$, we have
$$
0=\alpha(x_2)=e^{\alpha(a)}*\alpha(x_1)=-d\alpha(a)
$$
and
$$
0=\alpha(y_2)=e^{\alpha(b)}*\alpha(y_1)=-d\alpha(y).
$$
Moreover,  $L\otimes C$ and $N\otimes C$ are acyclic; therefore,
there exist $l \in (L\otimes A)^{-1}$ and  $k \in (N\otimes
A)^{-1}$ such that $d\alpha(l)=-\alpha(a)$ and
$d\alpha(k)=-\alpha(b)$. This implies $dl +a \in (L\otimes \ker
\alpha)^0$ and $dk +b \in (N\otimes \ker \alpha)^0$.

Set $w_1=dl +[x_1,l]\in Stab_A(x_1)$,  $w_2=dk +[y_1,k]\in
Stab_A(y_1)$ and define $\overline{a}=a \bullet w_1$ and
$\overline{b}=b \bullet w_2$.

We claim that $\overline{a}  \in (L\otimes \ker \alpha )^{0}$ and
$\overline{b}  \in (N\otimes \ker \alpha)^{0}$. Actually,
$$
\overline{a}=a \bullet w_1\equiv a+w_1 \equiv a+dl\ (mod \,
[L\otimes A,L\otimes A]);
$$
since $A\cdot A \subset \ker \alpha$, we conclude
$\overline{a}\equiv a+dl \equiv 0\ (mod \, L\otimes \ker \alpha)$.
An analogous calculation implies that $\overline{b}  \in (N\otimes
\ker \alpha)^{0}$.

Moreover, we note that
$$
e^{\overline{a}}* x_1=e^a e^{w_1}*x_1=x_2 \qquad \mbox{ and }
\qquad e^{\overline{b}}* y_1= e^b e^{w_2}*y_1=y_2.
$$
As to $e^{p_1}$, since $e^a Stab_A(x)e^{-a}=Stab_A(y)$ with $
y=e^a*x$, we have
$$
e^{-g(w_2)}e^{T }e^{p_1} e^{h(w_1)} =e^Se^{p_1}
$$
for some $S=df+[g(y_1),f]$ with $f \in (M \otimes A)^{-1}$.

Thus,
$$ e^{p_2}=  e^{g(b)}e^{T }e^{ p_1}   e^{-h(a)}=
$$
$$
e^{g(b)}e^{ g(w_2)}e^{-g(w_2)}e^{T }e^{ p_1} e^{h(w_1)}e^{-h(w_1)}
e^{-h(a)} = e^{g(\overline{b})} e^{S}e^{ p_1}
e^{-h(\overline{a})}.
$$
This implies that $e^S=e^{-g(\overline{b})}
e^{p_2}e^{h(\overline{a})}e^{-p_1}$ lies in the subgroup $\exp((M
\otimes \ker\alpha)^0)$ or equivalently $S=df+[g(y_1),f]\in (M
\otimes \ker\alpha)^0$.

Since $C$ is acyclic, the inclusion $M \otimes \ker \alpha \lrg M
\otimes A$ is a quasi-isomorphism and it remains a
quasi-isomorphism if we consider the deformed differentials
$d(-)+[g(y_1),-]$ on both $M \otimes \ker\alpha$ and $M \otimes A$
($g(y_1) $ satisfies the Maurer-Cartan equation  and so by
Remark~\ref{oss d(-)+[x,-]^2=0 se x in MC} $d(-)+[g(y_1),-]$ is a
differential).

Therefore, since the class of $S$ is trivial in $M \otimes A$, it
is also trivial in $M \otimes \ker \alpha$, i.e., there exists
$\overline{c} \in (M \otimes \ker\alpha)^{-1}$ such that
$S=\overline{T}=d\overline{c}+[g(y_1),\overline{c}]$.

In conclusion, $e^{p_2}=e^{g(\overline{b})} e^{\overline{T}}e^{
p_1} e^{-h(\overline{a})}$ and so $\widetilde{\Def}_{(h,g)}$ is a
predeformation functor.

\medskip

Next, we prove that  $\widetilde{\Def}_{(h,g)}$ is a  deformation
functor.

Let $C \in \Com_0$, then
$$
\widetilde{\MC}_{(h,g)}(C)=\{ (l,n,e^m)\in (L\otimes C)^1 \times
(N\otimes C)^1 \times \operatorname{exp}(M\otimes C)^0|
$$
$$
 d l=d n=0, \
g(n)=e^m*h(l)=h(l)-d m \}=Z^1(\cil \otimes C).
$$
and
\begin{equation}\label{equa EXTE Def(h,g)(C)=H^1 C_(h,g)}
\widetilde{\Def}_{(h,g)}(C)=
\end{equation}
$$
\widetilde{\MC}_{(h,g)}(C)  \ /
$$
$$
\{(-da,-db,dc+g(b)-h(a)) |\, a \in (L\otimes C)^0, \ b \in(N
\otimes C)^0, \ c  \in (M\otimes C)^{-1}\}.
$$
Therefore, $\widetilde{\Def}_{(h,g)}(C)$ is isomorphic to the
first cohomology group of the suspended cone of the pair of
morphisms $h:L \otimes C \lrg M \otimes C$ and $g:N \otimes C \lrg
M \otimes C$.

If $C$ is also acyclic, then $\widetilde{\Def}_{(h,g)}(C)=0$. This
implies that  $\widetilde{\Def}_{(h,g)}$ satisfies the condition
of Definition~\ref{dhef extended defo functor} and so it is a
deformation functor.

Finally, equation (\ref{equa EXTE Def(h,g)(C)=H^1 C_(h,g)}) also
implies that $T^i \widetilde{\Def}_{(h,g)}\cong
H^i(C^\cdot_{(h,g)})$.

\end{proof}
\begin{proof}[Proof of Lemma~\ref{lemma EXT DEF_L is deformation funt}]
It suffices to apply the previous theorem with $M=N=0$.
\end{proof}

\begin{teo}
$\widetilde{\Def}_{(h,g)}\cong \widetilde{\MC}_{(h,g)}^+$.
\end{teo}
\begin{proof}
The projection onto the quotient $\widetilde{\MC}_{(h,g)} \lrg
\widetilde{\Def}_{(h,g)}$ induces,  by Theorem~\ref{teo esistenza
funtor defor F+},  a map $\widetilde{\MC}_{(h,g)}^+ \lrg
\widetilde{\Def}_{(h,g)}$ that is a quasi-isomorphism by
Proposition~\ref{tangente di MC+ = h^i triplo cilindro} and
Theorem~\ref{teo EXT DEF_(h,g) is deformation}.
\end{proof}

\begin{cor}\label{cor Mpositiva exteDEF(h,g)restrict toDEf(h,g)}
Let $M$ be concentrated in non negative degree. Then for every
$(A,m_A) \in \Art$ we have
$\widetilde{\Def}_{(h,g)}(m_A)={\Def}_{(h,g)}(A)$.
\end{cor}
\begin{proof}
Obvious.
\end{proof}

\begin{oss}
Every commutative diagram of  differential graded Lie algebras
\begin{center}
$\xymatrix{ & L\ar[d]_h
\ar[r]^{\alpha'} & P\ar[d]^\eta \\
& M\ar[r]^\alpha & Q \\
N\ar[ur]^g \ar[r]^{\alpha''} & R \ar[ur]_\mu & \\}$
\end{center}
induces a natural transformation   $\widetilde{F}$ of the
associated deformation functors:
$$
\widetilde{F}: \widetilde{\Def}_{(h,g)} \lrg
\widetilde{\Def}_{(\eta,\mu)}.
$$
Moreover, the inverse function Theorem~\ref{TEO funzione inversa}
implies that $\widetilde{F}$ is an isomorphism if and only if the
map $(\alpha',\alpha,\alpha'')$ induces a quasi-isomorphism of
complexes  $\varphi^\cdot:\cil^\cdot \lrg
\operatorname{C}^\cdot_{(\eta,\mu)}$.
\end{oss}

\begin{proof}[Proof  of
Theorem~\ref{teo NO exte quasi iso C(h,g)-C(n,m)then DEf iso DEF}]
It is sufficient  to apply the previous remark and
Corollary~\ref{cor Mpositiva exteDEF(h,g)restrict toDEf(h,g)}.

\end{proof}

\begin{proof}[Proof
of Theorem~\ref{teo NO EXT H prodotto fibrato Q.I then ISO fun}]
It  suffices to apply  the inverse function Theorem~\ref{TEO
funzione inversa} and Corollary~\ref{cor Mpositiva
exteDEF(h,g)restrict toDEf(h,g)}.

\end{proof}

\subsection{Fibred product}
In Example~\ref{exe definizio M[t,dt]}, we have defined a DGLA
structure on $M[t,dt]=M \otimes \C[t,dt]$  and for each $a \in \C$
an evaluation morphism $e_a$ which is a surjective
quasi-isomorphism:
$$
e_a:M[t,dt] \lrg M, \ \  \ \ \ e_a(\sum m_it^i +n_it^i dt)=\sum
m_i a^i.
$$
Define $K \subset L \times N \times M[t,dt] \times M[s,ds]$ as
follows
$$
K=\{(l,n,m_1(t,dt),m_2(s,ds))|  \ h(l)=e_1(m_2(s,ds)),
g(n)=e_0(m_1(t,dt))  \}.
$$
$K$ is a DGLA with bracket and differential $\delta$ defined as
the natural ones on  each component.

Define the following morphisms of DGLAs:
$$
e_0: K \lrg M, \qquad (l,n,m_1(t,dt),m_2(s,ds))\longmapsto
e_0(m_1(t,dt))
$$
and
$$
e_1: K \lrg M, \qquad (l,n,m_1(t,dt),m_2(s,ds))\longmapsto
e_1(m_2(s,ds)).
$$

Then we can construct the following simplicial diagram of DGLAs
\begin{equation}\label{equa diagramma simplicial DGLA}
\xymatrix{L \times N \ar[rr]^{G \qquad \qquad } \ar@<-1ex>[dd]_h
\ar@<1ex>[dd]^g & & {\ \ \ \  K \subset L } \times N \times
M[t,dt]\times M[s,ds] \ar@<-1ex>[dd]_{e_1}
\ar@<1ex>[dd]^{e_0}  \\
& & \\
M  \ar[rr]^{id \qquad \qquad } &  & M   \\ }
\end{equation}
with:
\begin{center}
$\xymatrix{(l,n) \ar[rr]^G \ar@<-1.6ex>[ddd]_h \ar@<1ex>[dd]^g & &
{
\ \ \  (l,n,g(n),h(l))} \ar@<-1.6ex>[ddd]_{e_1} \ar@<1ex>[dd]^{e_0}  \\
& & \\
\ \ \  g(n) \ar[rr]^{id} & & \ \ \  g(n) \\
          h(l) \ar[rr]^{id} &  & h(l).   \\ }$
\end{center}
The  diagram is commutative in a simplicial meaning and $G$ is a
quasi-isomorphism.

\begin{lem}\label{lemma C_(f,g)quasiisomorfo C_(e_0,e_1)}
The complexes $C^\cdot_{ e_1-e_0 }$ and $C^\cdot_{ h-g }$ are
quasi-isomorphic.
\end{lem}
\begin{proof}
Consider the following commutative diagram of complexes

\begin{center}
$ \xymatrix{0 \ar[r] & (M^{\cdot-1},-d ) \ar[r] \ar[d]^{id} &
C^\cdot_{ h-g } \ar[r] \ar[d]^\alpha & (L^\cdot
\oplus N^\cdot,d) \ar[d]^G \ar[r] & 0  \\
0 \ar[r] &(M^{\cdot-1},-d )  \ar[r]  & C^\cdot_{ e_1-e_0 } \ar[r]
& (K^\cdot,\delta) \ar[r] & 0,  \\ }$
\end{center}
with $\alpha(l,n,m)= (l,n,g(n),h(l),m)$, for $(l,n,m) \in
C^\cdot_{ h-g }$. Since $id$ and $G$ are quasi-isomorphism,
$\alpha$ is a quasi-isomorphisms.
\end{proof}

\begin{prop}\label{prop mc_(f,g)^+ isomorfo mc_(e_0,e_1)+}
$\widetilde{\Def}_{(h,g)} \cong \widetilde{\Def}_{(e_1,e_0)} $.
\end{prop}
\begin{proof}
Since $G$ is a morphism of DGLAs, using diagram (\ref{equa
diagramma simplicial DGLA}), we can define a morphism of
Maurer-Cartan functors:
$$
G': \widetilde{\MC}_{(h,g)} \lrg \widetilde{\MC}_{(e_1,e_0)},
$$
$$
  (l,n,m) \longmapsto
(G(l,n),m), \ \ \ \mbox{ with } G(l,n)=(l,n,g(n),h(l)).
$$
It is well-defined since $G(l,n)\in K$ satisfies the Maurer-Cartan
equation and
$$
e^m*e_1(G(l,n))=e^m*e_1(h(l))=e^m*h(l)=g(n)=e_0(g(n))=e_0(G(l,n)).
$$
Therefore, $G'$ induces a morphism between the associated
deformation functors $ \widetilde{\Def}_{(h,g)} $ and
$\widetilde{\Def}_{(e_1,e_0)} $ that is a quasi-isomorphism by
Lemma~\ref{lemma C_(f,g)quasiisomorfo C_(e_0,e_1)} (and so an
isomorphism by Theorem~\ref{TEO funzione inversa}).
\end{proof}

\subsubsection{Definition of $\hil$, properties and barycentric
subdivision}\label{sezione def H}

Let $H \subset K$  be defined as follow
$$
H=\{(l,n,m_1(t,dt),m_2(s,ds))\in K| \ e_1(m_1(t,dt))=
e_0(m_2(s,ds)) \},
$$
or written in  detail
$$
H=\{(l,n,m_1(t,dt),m_2(s,ds))\in L \times  N \times M[t,dt]\times
M[s,ds]\ |
$$
$$
h(l)=e_1(m_2(s,ds)), \ g(n)=e_0(m_1(t,dt)),\ e_1(m_1(t,dt))=
e_0(m_2(s,ds)) \}.
$$

\medskip
Let $k=(l,n,m_1(t,dt),m_2(s,ds))\in K$. Then the pair $m_1(t,dt)$
and $m_2(s,ds)$ have fixed values at one of the extremes of the
unit interval. More precisely, the value of $ m_1(t,dt)$ is fixed
at the origin and $ m_2(s,ds)$ is fixed at 1, i.e.,
\begin{center}
\begin{xy}
,(-5,5);(15,5)**\dir{-} ,(-5,5)*{\cdot} ,(15,5)*{\cdot}
,(-5,3)*{0} ,(15,3)*{1} ,(5,3)*{t} ,(5,8)*{m_1(t,dt)}
,(35,5);(55,5)**\dir{-}  ,(35,5)*{\cdot} ,(55,5)*{\cdot}
,(35,3)*{0} ,(55,3)*{1} ,(-5,14)*{e_0(m_1)=g(n)\ \ \ \ }
,(55,14)*{e_1(m_2)=h(l)\ \ \ \ } ,(45,3)*{s} ,(45,8)*{m_2(s,ds)}
  ,(-5,12);,(-5,6)**\dir{.} ,(55,12);,(55,6)**\dir{.} ,(57,5)*{.}
\end{xy}
\end{center}
If $k$ also lies in $H$, then there are conditions on the other
extremes: the value of $m_1(t,dt)$ at 1 has to coincide with the
value of $m_2(s,ds)$ at 0.

Let
\begin{equation}\label{equa def H breve f=e_0 g=e_1}
\hil=
\end{equation}
$$
\{ (l,n,m(t,dt)) \in L\times N \times M[t,dt]\ | \,
h(l)=e_1(m(t,dt)), \ g(n)=e_0(m(t,dt)) \}.
$$
Since $e_i$ are morphisms of DGLAs, it is clear that $\hil$ is a
DGLA.

Moreover, considering the barycentric subdivision we get an
injective quasi-isomorphism
$$
\hil \hookrightarrow H,
$$
$$
(l,n,m(t,dt))\longmapsto (l,n,m\,(\, \displaystyle\frac{1}{2}\,
t,dt),m\,(\,\frac{s+1}{2},ds)).
$$


\begin{dhef}\label{dhef DGLA H(h,g)}
$\hil$ is the \emph{differential graded Lie algebra associated
with the pair $(h,g)$}.
\end{dhef}

\begin{prop}
$\hil$ is  quasi-isomorphic to the complex $C^\cdot_{e_1-e_0}$.
\end{prop}

\begin{proof}
It is sufficient to consider the following commutative diagram of
complexes
\begin{center}
$\xymatrix{\hil\ar[d]\ar@{^{(}->}[r] & K \ar[d]^{e_1-e_0} \\
           0\ar[r] & M,   \\ }$
\end{center}
where $e_1-e_0$ is surjective.
\end{proof}

\begin{prop}\label{prop if g-h surj then
EXTE DEF_H qiso DEf_(h,g)} Let $h: L \lrg M$ and $g:N \lrg M$ be
morphisms of DGLAs. If the morphism $g-h:N\times L \lrg M$ is
surjective, then $\widetilde{\Def}_{L\times _N M}$ is isomorphic
to $\widetilde{\Def}_{(h,g)}$.
\end{prop}

\begin{proof}
We recall that, by definition,
$$
\widetilde{\MC}_{L\times _M N}(A)=
$$
$$
\{(l,n) \in (L\otimes A)^1 \times (N\otimes A)^1| dl+
\displaystyle\frac{1}{2}[l,l]=0, \ dn+
\displaystyle\frac{1}{2}[n,n]=0, \ h(l)=g(n) \},
$$
and
$$
\widetilde{\MC}_{(h,g)}(A)=\{(l,n,m)\in (L\otimes A)^1\times
(N\otimes A)^1 \times (M\otimes A)^0 |
$$
$$
dl+ \displaystyle\frac{1}{2}[l,l]=0, \ dn+
\displaystyle\frac{1}{2}[n,n]=0,\  g(n)=e^m *h(l) \}.
$$
Moreover, $ T^i \widetilde{\Def}_{L\times _M N}\cong H^i(L\times
_M N)$ and $T^i \widetilde{\Def}_{(h,g)}  \cong H^i(
C^\cdot_{(h,g)})$.

Let
$$
\psi: \widetilde{\MC}_{L\times _M N}(A) \lrg
\widetilde{\MC}_{(h,g)}(A),
$$
$$
(l,n) \longmapsto (l,n,0),
$$
and $   \varphi:\widetilde{\Def}_{L\times _N
M}\lrg\widetilde{\Def}_{(h,g)}$ the induced map between the
associated extended deformation functors. Then $\varphi$ is a
quasi-isomorphism.

For completeness we state all details, denoting by the same
$\varphi$ the map induced on cohomology.

$\varphi $ is injective. Suppose that $
\varphi([(l,n)])=[(l,n,0)]=0$ in $H^i( C^\cdot_{(h,g)})$. Then,
$[(l,n,0)]=(dr,ds, -dt -g(s)+h(r) ) $ with $r \in (L  )^{i-1}$, $s
\in (N )^{i-1}$ and $t \in (M )^{i-2}$.

Since $g-h$ is surjective, there exist  $p \in L^{i-2} $ and   $q
\in N ^{i-2} $ such that $g(q)-h(p)=-t$ and so $g(dq)-h(dp)=-dt$.
Let $(l',n')=(r-dp,s-dq)$; then
$h(l')=h(r)-h(dp)=g(s)+dt-dt-g(dq)=g(n')$.

Therefore, $(l',n')\in L\times_M N$ and $d(l',n')=(dr,ds)=(l,n)$
and so $[(l,n)]=0 \in H^i(L\times_M N)$.

\smallskip

$\varphi $ is surjective. Let $[(l,n,m)] \in H^i(
C^\cdot_{(h,g)})$, i.e., $dl=dn=0$ and $-dm-g(n)+h(l)=0$. Since
the class $[(l,n,m)]$ coincides with the class $[l+dr,n+ds,
m-dt-g(s)+h(r)] \in H^i( C^\cdot_{(h,g)}) $, we are looking for $r
\in L^{i-1}  $, $s \in N^{i-1} $ and $t  \in M^{i-2} $ such that
$(l+dr,n+ds) \in H^i(L\times_M N)$ and $m-dt-g(s)+h(r)=0$ (thus,
$[(l,n,m)] \in \image (\varphi)$).

Since $g-h$ is surjective,  there exist  $r \in L^{i-1} $ and $s
\in N^{i-1} $, such that $g(s)-h(r)=m-dt$ and so $g(ds)-h(dr)=dm$.
Therefore, $h(l+dr)=h(l)+g(ds)-dm =h(l)+g(ds)-h(l)+g(n)=g(n+ds)$,
that is, $(l+dr,n+ds) \in L\times_M N$ and $\varphi
([(l+dr,n+ds)])=[(l,n,m)]$.
\end{proof}

\begin{oss}
Let $\alpha_0, \alpha_1:L \lrg M$ be morphisms of DGLAs with
$\alpha_0- \alpha_1$ surjective and $T=\{t \in L | \, \alpha_0(t)
= \alpha_1(t)\}$ ($T$ is called the equalizer of  $\alpha_0$ and $
\alpha_1$). In this particular case, the previous proposition
implies $\widetilde{\Def}_T\cong \widetilde{\Def}_{(\alpha_1,
\alpha_0)} $.
\end{oss}

In conclusion, we have the following theorem.

\begin{teo}\label{teo EXT DEF_H=DEF(h,g)}
$\widetilde{\Def}_{\hil}\cong \widetilde{\Def}_{(h,g)} $.
\end{teo}

\begin{proof}
It  is sufficient to apply Proposition~\ref{prop if g-h surj then
EXTE DEF_H qiso DEf_(h,g)} to $e_0,e_1:K\lrg M$ (with $e_0-e_1$
surjective) to conclude that $\widetilde{\Def}_{\hil}\cong
\widetilde{\Def}_{(e_1,e_0)}$;  then Proposition~\ref{prop
mc_(f,g)^+ isomorfo mc_(e_0,e_1)+} is used to complete the proof.
\end{proof}

\chapter{Deformations  of holomorphic maps}

This chapter  is devoted to the main topic of this thesis:
infinitesimal deformations of holomorphic maps of complex compact
manifolds.

These deformations were first studied during the 70s    by
E.~Horikawa in  his works (\cite{bib HorikawaI} and \cite{bib
Horikawa III}) and then by M.~Namba \cite{bib Namba} and Z.~Ran
\cite{bib RAn mappe}.

Our purpose is to study these deformations using a technique based
on differential graded Lie algebras.



\medskip
\noindent{\bf Note.} Unless otherwise specified, \emph{$X$ and $Y$
are compact complex connected smooth varieties}.

\section{Deformation functor $\Def(f)$ of holomorphic
maps}\label{sezio definizio funtore DEF(f)}

The main references for this section are \cite{bib HorikawaI} and
\cite[Section~3.6]{bib Namba}.

\begin{dhef}\label{def infin deform f FIXED target domain}
Let $f:X \lrg Y$ be a  holomorphic map and $A  \in \Art$. An
\emph{infinitesimal deformation of f with fixed domain and
codomain over $\Spec(A)$} is a commutative diagram
\begin{center}
$\xymatrix{X\times S\ar[rr]^{\mathcal{F}} \ar[dr] &  & Y \times S \ar[ld] \\
            & S, &  \\ }$
\end{center}
where $S=\Spec(A)$, the morphisms onto $S$ are the projections,
$\mathcal{F}$ is a  holomorphic map and $f$ coincides with the
restriction of $\mathcal{F}$ to the fibers over the closed point
of $S$.

If $A=\kepsi$ we have a { \it first order deformation of $f$  with
fixed domain and codomain}.
\end{dhef}

Two infinitesimal deformations of $f$ with fixed domain and
codomain
\begin{center}
$\xymatrix{X\times S\ar[rr]^{\mathcal{F}} \ar[dr]  &  & Y\times S
\ar[ld] &   \mbox{ and }&   X\times S\ar[rr]^{\mathcal{F}'}
\ar[dr]  & &
Y\times S \ar[ld] \\
            & S & &   &    &S \\ }$
\end{center}
are equivalent if there exist  automorphisms $\phi:X \lrg X $ and
$\psi:Y \lrg Y$ such that the following diagram is commutative:
\begin{center}
$\xymatrix{X\times S\ar[rr]^{\mathcal{F}} \ar[d]_{(\phi,id)} &  &
Y\times S
\ar[d]^{(\psi,id)} \\
X\times S\ar[rr]^{\mathcal{F}'} &  & Y\times S. \\ }$
\end{center}

\begin{dhef}
The \emph{infinitesimal deformation  functor } \\ $\Def
(X\stackrel{f}{\lrg} Y)$ of  infinitesimal deformations of a
holomorphic map $f$ with fixed domain and codomain  is defined as
follows:
$$
\Def (X\stackrel{f}{\lrg} Y): \Art \lrg \Set,
$$
$$
  A\longmapsto \Def(X\stackrel{f}{\lrg} Y)
(A)=\left\{\begin{array}{c} \mbox{ isomorphism  classes  of }\\
 \mbox{  infinitesimal   deformations   of   $f$ } \\
 \mbox{  with   fixed   domain   and   codomain}\\
   \mbox{ over $\Spec(A)$} \\
\end{array} \right\}.
$$
\end{dhef}

\begin{oss}\label{oss def X e Y fixed =Gamma cXxY XxY fixed}
When the domain and codomain are fixed, an infinitesimal
deformation of a  holomorphic map $f$ can be viewed as an
infinitesimal deformation  of the graph of the map $f$ in the
product $ X\times Y$, with $X \times Y$ fixed.
\end{oss}

\bigskip

In this case, we are just deforming the map $f$. In general, we
can also deform both the domain and the codomain.

\begin{dhef}
Let $f:X \lrg Y$ be a  holomorphic map and $A \in \Art$. An { \it
infinitesimal deformation of f  over $\Spec(A)$} is a commutative
diagram of complex spaces
\begin{center}
$\xymatrix{X_A\ar[rr]^{\mathcal{F}} \ar[dr]_\pi &  & Y_A
\ar[ld]^\mu \\
            & S, &  \\ }$
\end{center}
where $S=\Spec(A)$, $(X_A,\pi,S)$ and $(Y_A,\mu,S)$ are
infinitesimal deformations of $X$ and $Y$, respectively
(Definition~\ref{dhef infintesimal defor of X}), $\mathcal{F}$ is
a  holomorphic map that restricted to the fibers over the closed
point of $S$ coincides with $f$.

If $A=\kepsi$ we have a \emph{first order deformation} of $f$.
\end{dhef}

\begin{dhef}

Let
\begin{center}
$\xymatrix{X_A\ar[rr]^{\mathcal{F}} \ar[dr]_\pi &  & Y_A
\ar[ld]^\mu & & \mbox{ and }&   &X'_A\ar[rr]^{\mathcal{F}'}
\ar[dr]_{\pi'} & &
Y'_A \ar[ld]^{\mu'}\\
            & S & & & &    &  &S \\ }$
\end{center}
be two infinitesimal deformations of $f$. They are
\emph{equivalent} if there exist bi-holomorphic maps $\phi : X_A
\lrg X'_A$ and $\psi: Y_A\lrg Y'_A$ (that are equivalence of
infinitesimal deformations of $X$ and $Y$, respectively) such that
the following diagram is commutative:
\begin{center}
$\xymatrix{X_A\ar[rr]^{\mathcal{F}} \ar[d]_\phi &  & Y_A
\ar[d]^\psi \\
X'_A\ar[rr]^{\mathcal{F}'} &  & Y'_A. \\ }$
\end{center}
\end{dhef}

\begin{dhef}\label{def funore DEF(f)}
The \emph{functor of infinitesimal deformations} of a  holomorphic
map $f:X \lrg Y$ is
$$
\Def(f): \Art \lrg \Set,
$$
$$
\ \ \ \ \ \  \ \ \ \  \ \ A\longmapsto \Def(f)
(A)=\left\{\begin{array}{c} \mbox{ isomorphism  classes  of}\\
 \mbox{ infinitesimal   deformations   of }\\
 f  \mbox { over }  \Spec(A) \\
\end{array} \right\}.
$$
\end{dhef}

\begin{prop}
$\Def(f)$  is a deformation functor, since it satisfies the
conditions of Definition~\ref{defin deformation funtore}.
\end{prop}
\begin{proof}
It follows from the fact that the functors $\Def_X$  and $\Def_Y$
of infinitesimal deformations of $X$ and $Y$ are deformation
functors.
\end{proof}

\begin{oss}\label{oss defo mapp=Def grafo prod defX.defY}
In this general case, the infinitesimal deformations of $f$ can be
interpreted as infinitesimal deformations  $\widetilde{\Gamma} $
of the graph $\Gamma$ of the map $f$ in the product $X \times Y$,
such that the induced  deformations $\widetilde{X \times Y}$ of $X
\times Y$ are products of infinitesimal deformations of $X$ and of
$Y$. Since not all the deformations of a product are product of
deformations (see Remark~\ref{oss DEF(XxY)> DEF(X) x DEF(Y)}), we
are not just considering the deformations of the graph in the
product.

Moreover, with this interpretation, two infinitesimal deformations
$\widetilde{\Gamma} \subset \widetilde{X \times Y} $ and
$\widetilde{\Gamma}' \subset \widetilde{X \times Y}'$ are
equivalent if there exists an isomorphism $\phi:\widetilde{X
\times Y}\lrg \widetilde{X \times Y}'$ of infinitesimal
deformations of $X\times Y$ such that
$\phi(\widetilde{\Gamma})=\widetilde{\Gamma}'$.

\end{oss}

\subsection{Tangent and obstruction  spaces of  $ \Def(f)$}

Let $\mathcal{U}=\{U_i\}$ and  $\mathcal{W}=\{W_j\}$ be Stein open
covers of $X$ and $Y$, respectively, such that $f(U_i)\subset V_i$
for each $i$. For any integer $p\geq 0$, let
$$
\check{C}^p(\mathcal{U},\mathcal{W},\Theta_X,\Theta_Y,f^*\Theta_Y)=
\check{C}^p(\mathcal{U},\Theta_X) \oplus
\check{C}^p(\mathcal{W},\Theta_Y ) \oplus
\check{C}^{p-1}(\mathcal{U},f^*\Theta_Y),
$$
(with $\check{C}^{-1}(\mathcal{U},f^*\Theta_Y)=0$). Define a
linear map
$$
\check{D}: \check{C}^p(\mathcal{U},\mathcal{W},\Theta_X,
\Theta_Y,f^*\Theta_Y) \lrg
\check{C}^{p+1}(\mathcal{U},\mathcal{W},\Theta_X,
\Theta_Y,f^*\Theta_Y)
$$
$$
(x,y,z) \longmapsto (\check{\delta} x,\check{\delta}
y,\check{\delta} z + (-1)^p(f_*x-f^*y)).
$$
Using the equalities  $f_* \check{\delta}=\check{\delta} f_*$ and
$f^* \check{\delta}=\check{\delta} f^*$, we conclude that
$\check{D}$ is a differential ($\check{D} \circ \check{D}=0 $).
Therefore, the cohomology $\C$-vector spaces
$\check{H}^p(\mathcal{U},\mathcal{W}, \Theta_X,\Theta_Y,f^*
\Theta_Y)$ are well defined.

\begin{lem}
$\check{H}^\cdot (\mathcal{U},\mathcal{W},\Theta_X,\Theta_Y,f^*
\Theta_Y)$ does not  depend on the choice of the covers and so we
denote it by $\check{H}^\cdot (\Theta_X,\Theta_Y,f^*\Theta_Y)$.
\end{lem}
\begin{proof}
The following linear maps are well defined:
$$
\check{H}^p(\mathcal{U},\Theta_X)\oplus
\check{H}^p(\mathcal{V},\Theta_Y) \lrg
\check{H}^p(\mathcal{U},f^*\Theta_Y)
$$
$$
(\{n_1\},\{ n_2\})\longmapsto \{f_*n_1 -f^*n_2\};
$$
$$
\check{H}^{p }(\mathcal{U},f^*\Theta_Y) \lrg \check{H}^{p+1}
(\mathcal{U},\mathcal{W},\Theta_X,\Theta_Y,f^*\Theta_Y)
$$
$$
 \{a\} \longmapsto \{(0,0,a)\};
$$
$$
H^p (\mathcal{U},\mathcal{W},\Theta_X,\Theta_Y,f^*\Theta_Y)  \lrg
\check{H}^p(\mathcal{U},\Theta_X)\oplus
\check{H}^p(\mathcal{V},\Theta_Y)
$$
$$
\{(n_1,n_2,a)\} \longmapsto (\{n_1\},\{ n_2\}).
$$
Then the sequence below is exact:
$$
\cdots \lrg \check{H}^{p }(\mathcal{U},\Theta_X)\oplus
\check{H}^{p }(\mathcal{V},\Theta_Y) \lrg \check{H}^{p }
(\mathcal{U},f^*\Theta_Y) \lrg
$$
$$
\lrg \check{H}^{p+1}
(\mathcal{U},\mathcal{W},\Theta_X,\Theta_Y,f^*\Theta_Y)  \lrg
$$
$$
\lrg   \check{H}^{p+1 }(\mathcal{U},\Theta_X)\oplus \check{H}^{p+1
}(\mathcal{V},\Theta_Y) \lrg \check{H}^{p+1}
(\mathcal{U},f^*\Theta_Y) \lrg \cdots.
$$
Moreover, $\check{H}^\cdot (\mathcal{U},\Theta_X)$,
$\check{H}^\cdot (\mathcal{U}, \Theta_X)$ and $\check{H}^\cdot
(\mathcal{U}, f^*\Theta_Y)$ does not  depend on the choice of the
Stein open covers $\mathcal{U}$ and $\mathcal{V}$. Hence applying
the so called  \lq\lq five lemma", $\check{H}^\cdot
(\mathcal{U},\mathcal{W}, \Theta_X,\Theta_Y,f^* \Theta_Y)$ does
not  depend on the choice of the covers.
\end{proof}

In \cite{bib HorikawaI}, E.~Horikawa used the vector space
$\check{H}^\cdot (\Theta_X,\Theta_Y,f^*\Theta_Y)$ to describe the
tangent and obstruction spaces of the deformation functor
$\Def(f)$.

\begin{teo}\label{teo TANG e OSTR DEF(f) horikawa}
$\check{H}^1 (\Theta_X,\Theta_Y,f^*\Theta_Y)$ is in 1-1
correspondence with the first order deformations of $f:X \lrg Y$.

The obstruction space of the functor $\Def(f)$ is naturally
contained in $\check{H}^2(\Theta_X,\Theta_Y,f^*\Theta_Y)$.
\end{teo}

\begin{proof}
See \cite[Section~3.6]{bib Namba}.

\end{proof}

\begin{oss}\label{oss class H^2(TX,TY,f*TY) def by  Xe,Ye}
Consider a first order deformation $f_\epsi$ of $f$: in
particular, we are considering   first order deformations
$X_\epsi$ and $Y_\epsi$  of $X$ and  $Y$.  Using Theorem~\ref{teo
sernesi first ordine=H^1(X,T_X)}, we associate with $X_\epsi$ a
class $x \in \check{H}^1(X,\Theta_X)$ and to $Y_\epsi$ a class $y
\in \check{H}^1(Y,\Theta_Y)$.

Therefore, the class in $\check{H}^1 (\Theta_X,\Theta_Y,
f^*\Theta_Y)$ associated with $f_\epsi$ is $[(x,y,z)]$ with $z \in
\check{C}^{0}(\mathcal{U},f^*\Theta_Y)$ such that $\check{\delta}
z = f_*x-f^*y$.

Analogously, let  $0\lrg J \lrg B \lrg A\lrg 0$ be a small
extension and $f_A$ an infinitesimal deformation of $f$ over
$\Spec (A)$. If $h \in \check{H}^2(X,\Theta_X)$ and $k \in
\check{H}^2(Y,\Theta_Y)$ are the obstruction classes associated
with $X_A$ and $Y_A$, respectively, then the obstruction class in
$\check{H}^2 (\Theta_X,\Theta_Y, f^*\Theta_Y)$ associated with
$f_A$ is $[(h,k,r)]$, with $r \in
\check{C}^{1}(\mathcal{U},f^*\Theta_Y)$ such that $\check{\delta}
r = -(f_*h-f^*k)$.
\end{oss}

\begin{oss}
Using \v{C}ech's cohomology, we have defined the $\C$-vector space
$\check{H}^* ( \Theta_X, \Theta_Y,f^* \Theta_Y)$. For convenience
we reinterpret it, using  Dolbeault's cohomology.

Let $(B^\cdot,D_{\debar})$  be the complex below
$$
B^p=A_X^{(0,p)}( \Theta_X) \oplus A_Y^{(0,p)}(\Theta_Y) \oplus
  A_X^{(0,p-1)}( f^*\Theta_Y)
$$
and
$$
D_{\debar}: B^p \lrg B^{p+1}, \qquad (x,y,z)\longmapsto(\debar
x,\debar y, \debar z+ (-1)^p(f_*x-f^*y)).
$$

\begin{lem}\label{lemm C horikawa quasi iso to B}
The complexes $(\check{C}^\cdot(\mathcal{U},\mathcal{W}, \Theta_X,
\Theta_Y,f^*\Theta_Y), \check{D})$  and  $(B^\cdot,D_{\debar})$
  are quasi-isomorphic.
\end{lem}
\begin{proof}
Let $\mathcal{U}=\{U_i\}$ and  $\mathcal{W}=\{W_j\}$ as above and
denote by $\phi_1:\check{C}^* ( \mathcal{U} ,\Theta_X)\lrg
A_X^{0,*}(\Theta_X)$, $\phi_2:\check{C}^* ( \mathcal{V}
,\Theta_Y)\lrg A_Y^{0,*}(\Theta_Y)$ and $\phi_3:\check{C}^* (
\mathcal{U} ,f^*\Theta_X)\lrg A_X^{0,*}(f^*\Theta_X)$  the
quasi-isomorphisms of complexes of Leray's theorem, defined in
Section~\ref{sottosezione cech e leray theorem}. We recall that
$\phi_i \check{\delta}=\debar \phi_i$, $f_*\phi_i =\phi_i f_*$ and
$f^*\phi_i =\phi_i f^*$, for each $i=1,2,3$ (see
Section~\ref{sezione f_* and f^*}).

Next, define the following morphism
$$
\gamma: (\check{C}^\cdot(\mathcal{U},\mathcal{W},\Theta_X,
\Theta_Y,f^*\Theta_Y), \check{D}) \lrg (B^\cdot,D_{\debar}),
$$
$$
\gamma(x,y,z)=(\phi_1(x),\phi_2(y),\phi_3(z)).
$$
Then $\gamma$ is a morphism of complexes:; actually, for  any
$(x,y,z) \in \check{C}^p(\mathcal{U},\mathcal{W},\Theta_X,
\Theta_Y,f^*\Theta_Y)$
\begin{center}
$\xymatrix{(x,y,z) \ar[r]^{\gamma\qquad \ \ } \ar[d]_{\check{D} }
&
 (\phi_1(x),\phi_2(y),\phi_3(z)) \ar[d]^{D_{\debar}}\\
(\check{\delta} x,\check{\delta} y,\check{\delta} z
+(-1)^p(f_*x-f^*y)) \ar[r]  & b \\ }$
\end{center}
where $b=(\debar\phi_1( x),\debar\phi_2(y),\debar\phi_3(z)
+(-1)^p(f_*\phi_1(x)-f^*\phi_2(y)))$.

Moreover, we have the following commutative diagram

\begin{center}
$\xymatrix{0\ar[r] &  \check{C}^{\cdot-1}(\mathcal{U},f^*\Theta_Y)
\ar[r] \ar[d]^{\phi_3} & \mathcal{\check{C}}^\cdot
\ar[r]\ar[d]^\gamma & \check{C}^\cdot(\mathcal{U},\Theta_X) \oplus
\check{C}^\cdot(\mathcal{W},\Theta_Y)\ar[r]
\ar[d]^{(\phi_1,\phi_2)}  &0 \\
0\ar[r] &   A_X^{(0,\cdot-1)}( f^*\Theta_Y) \ar[r]  & B^\cdot
\ar[r] & A_X^{(0,\cdot)} ( \Theta_X) \oplus A_Y^{(0,\cdot)}
(\Theta_Y) \ar[r]&0  \\ }$
\end{center}
where $\mathcal{\check{C}}^\cdot$ stands for
$\check{C}^\cdot(\mathcal{U},\mathcal{W},
\Theta_X,\Theta_Y,f^*\Theta_Y)$.

Since $\phi_3$ and $(\phi_1,\phi_2)$ are quasi-isomorphisms,
$\gamma$ is a quasi-isomorphism.

\end{proof}
\end{oss}

\section{Infinitesimal deformations of holomorphic maps}
\label{sez teo Def_(h,g)Def (f)}

Let $f:X \lrg Y$ be a  holomorphic map and $\Gamma$ its graph in
$Z:=X\times Y$.

Let
$$
F:X \lrg \Gamma \subseteq Z:=X \times Y,
$$
$$
\ \ \ x \longmapsto (x,f(x)),
$$
and $p:Z\lrg X$ and $q:Z \lrg Y$  the natural projections.

Then we have the following commutative diagram:
\begin{center}
$\xymatrix{X \ar[rrrr]^{F} \ar[ddrrr]_{id} \ar[ddrrrrr]^f & &  & &
Z\ar[ddl]_{ p} \ar[ddr]_{ \ \ q} & \\
& & &  & &  \\
& & & X & & Y. \\ }$
\end{center}
In particular, $F^*\circ p^*=id$ and $F^* \circ q^*=f^*$.

Since $\Theta_Z=p^*\Theta_X \oplus q^*\Theta_Y$, it follows that
$F^*(\Theta_{Z})=\Theta_X \oplus f^*\Theta_Y$.

Define the morphism $\gamma:\Theta_Z \lrg f^*\Theta_Y$ as the
product
$$
\gamma: \Theta_Z \stackrel{F^*}{\lrg}\Theta_X \oplus f^* \Theta_Y
\stackrel{(f^*,-id)}{\lrg} f^*\Theta_Y
$$
and let $\pi$ be the following surjective morphism:
$$
\shA^{0,j}_{Z}(\Theta_{Z}) \stackrel{\pi}{\lrg}
\shA_X^{0,j}(f^*\Theta_Y) \lrg 0,
$$
$$
\pi(\omega\, u)= F^*(\omega)\gamma(u), \qquad \forall \ \omega \in
\shA^{0,j}_{Z}, \ u \in \Theta_{Z}.
$$
Since each $u \in \Theta_{Z}$ can be written as $u=p^*v_1+q^*v_2$,
for some $v \in\Theta_X$ and $w \in \Theta_Y$, we also have
$$
\pi (\omega u)=F^*(\omega)(f_*(v_1)-f^*(v_2)).
$$
Since $F^*\debar=\debar F^*$, $\pi$ is a morphism of complexes.

For the reader's convenience, we give an explicit description of
the map $\pi$.

Let $\mathcal{U}=\{U_i\}_{i \in I}$ and $\mathcal{V}=\{V_i\}_{i
\in I}$ be finite Stein open  covers of $X$ and $Y$, respectively,
such that $f(U_i) \subset V_i$ ($U_i$ is allowed to be empty).
Moreover, let  $z=(z_1,z_2,\ldots ,z_n)$ on $U_i$ and $w=(w_1,w_2,
\ldots w_m)$ on $V_i$ be local holomorphic coordinate systems for
each $i \in I$. As in Section~\ref{sezione f_* and f^*}, if $
\displaystyle v_1= \sum_{i=1}^n \varphi_i(z) \frac{\de}{\de z_i}$
and $ \displaystyle v_2= \sum_{h=1}^m \psi_h(w) \frac{\de}{\de
w_h}$ then $ \displaystyle f_* v_1 =\sum_{i=1}^n \varphi_i(z)
\sum_{h=1}^m \frac{\de f_{h}}{\de z_i} \frac{\de }{\de w_h }  $
and  $ \displaystyle  f^*v_2=\sum_{h=1}^m \psi_h(f(z))   \frac{\de
}{\de w_h }$. Let $K$ and $J$ be   multi-indexes of length
  $k$ and $j$, respectively, and fix $\omega =\Phi(z,w)
d\overline{z}_K \wedge d\overline{w}_J \in \shA^{0,k+j}_{Z}$. Then
$$
\pi (\omega u)=F^*(\omega)(f_*(v_1)-f^*(v_2))=
$$
$$
\Phi(z,f(z)) d\overline{z}_K \wedge \debar f_J \sum_{h=1}^m \left(
\sum_{i=1}^n \varphi_i(z)  \frac{\de f_{h}}{\de z_i} -
\psi_h(f(z)) \right) \frac{\de }{\de w_h }.
$$

Let $\mathcal{L}$ be the kernel of $\pi$:
\begin{equation}\label{equa def L con f^*T_X}
0 \lrg \mathcal{L}\stackrel{h}{ \lrg} \shA^{0,*}_{Z}(\Theta_{Z})
\stackrel{\pi}{\lrg} \shA_X^{0,*}(f^*\Theta_Y) \lrg 0
\end{equation}
and $h:\mathcal{L}\lrg \shA^{0,*}_{Z}(\Theta_{Z}) $ the inclusion.

\begin{lem}\label{lemma Lfascio is DGLA}
$\mathcal{L}$ is a sheaf of differential graded subalgebras of
$\shA^{0,*}_{Z}(\Theta_{Z})$ and $h$ is a morphism of differential
graded Lie algebras.
\end{lem}
\begin{proof}
There is a canonical  isomorphism between  the normal bundle
$N_{\Gamma|Z}$ of $\Gamma$ in $Z$  and the pull-back $f^*T_Y$.

Therefore, (\ref{equa def L con f^*T_X}) reduces to
$$
0 \lrg \mathcal{L}\stackrel{h}{ \lrg} \shA^{0,*}_{Z}(\Theta_{Z})
\stackrel{\pi}{\lrg} \shA_\Gamma^{0,*}(N_{\Gamma|Z}) \lrg 0.
$$

\bigskip
Then,  by Lemma~\ref{lemm L manetti is DGLA}, $ \mathcal{L}$ is a
sheaf of differential graded subalgebras of $
A^{0,*}_{Z}(\Theta_{Z})$.

\end{proof}

Let $L$ be the differential graded  Lie algebra of global sections
of $\mathcal{L}$.

\smallskip

Let $M$ be the Kodaira-Spencer algebra of the product $X\times Y$:
$M=KS_{X\times Y}=KS_Z$ and $h:L \lrg M$ be the inclusion.
\smallskip

Let $N=KS_X \times KS_Y$ be the product of the Kodaira-Spencer
algebras of $X$ and of $Y$  and  $g:KS_X \times KS_Y\lrg
KS_{X\times Y}$ be given by $g=p^*+q^*$, i.e., $g(n_1,n_2)=p^*n_1
+q^*n_2$ (for $n=(n_1,n_2)$ we  also use the notation $g(n)$).

Therefore, we get a diagram
\begin{center}
\begin{equation}\label{equa diagrL-> KS(XxY) <-KS(X)xKS(Y)}
\xymatrix{ & & L \ar@{^{(}->}[d]^h   \\
N=KS_X \times KS_Y \ar[rr]^{g=(p^*,q^*)} &  & M= KS_{X\times Y}. \\
}
\end{equation}
\end{center}

\begin{oss}

Given the morphisms of DGLAs $h:L \lrg {KS} _{X \times Y}$ and
$g:KS_X\times  KS_Y \lrg {KX} _{X \times Y}$ we can consider   the
complex $(\cil^\cdot,D)$, with $C^i_{(h,g)}=L^i \oplus
KS_X^i\oplus KS_Y^i \oplus {KS}^{i-1}_{X \times Y}$ and
differential $D(l,n_1,n_2,m)=(\debar l, \debar n_1, \debar n_2,
-\debar m -p^*n_1 -q^* n_2 + h(l))$.

Using the morphism $\pi : KS_{X\times Y} \lrg
A_X^{0,*}(f^*\Theta_Y)$ we can define a morphism
$$
\beta :(\cil^\cdot, D) \lrg (B^\cdot, D_{\debar}),
$$
$$
\beta (l,n_1,n_2,m)=(n_1, n_2,(-1)^i\pi (m)) \qquad  \forall \
(l,n_1,n_2,m)\in \cil^i.
$$

\begin{prop}\label{prop isomor complesso Namba con Cono}
$\beta:(\cil^\cdot, D) \lrg (B^\cdot, D_{\debar}) $ is a morphism
of complexes which is a quasi-isomorphism.

\end{prop}
\begin{proof}
$\beta$ commutes with the differentials, i.e., $\beta \circ D=
D_{\debar} \circ  \beta$; in fact,  for each $(l,n_1,n_2,m) \in
C^i_{(h,g)}$ we have the following commutative diagram
\begin{center}
$\xymatrix{(l,n_1,n_2,m) \ar[r]^{\beta} \ar[d]_D  &  (n_1,
n_2,(-1)^i\pi (m)) \ar[d]^{D_{\debar}} \\
 c \ar[r]^\beta & b, \\ }$
\end{center}
where $c=(\debar l, \debar n_1,  \debar n_2, -\debar m -p^*n_1
-q^* n_2 + h(l))$ and $b=\beta(c)= (\debar n_1,\debar n_2, (-1)^{i
} ( \debar \pi(m) +f_*n_1 -f^* n_2))= D_{\debar} (n_1,
n_2,(-1)^i\pi (m)) )$.

Therefore, $\beta$ induces a map in cohomology that we again  call
$\beta$  that is a quasi-isomorphism. The proof is standard but we
write it.

\emph{$\beta$ is injective}. Let $[(l,n_1,n_2,m)]\in H^i(C^
\cdot_{(h,g)})$\footnote{In particular, this implies that $-\debar
m -p^*n_1-q^*n_2+ l =0$ } be a class such that $\beta
([(l,n_1,n_2,m)])= [(n_1, n_2,(-1)^i\pi (m))]$ is zero in
$H^i(B^\cdot)$.

Then there exists $(a_1,a_2,b) \in B^{i-1}$ such that $ n_1=\debar
a_1$, $  n_2=\debar a_2$, and $(-1)^i \pi (m)=\debar b
+(-1)^{i-1}(f_*a_1-f^*a_2)$, so that $\pi(m)=(-1)^i \debar
b-f_*a_1+f^*a_2$.

Let $n_1'= a_1 \in KS_X^{i-1}$ and $n_2 '= a_2 \in KS_Y^{i-1}$;
then $ \debar n_1'=n_1$ and  $ \debar n_2'=n_2$.

Let $z \in KS_{X\times Y}^{i-2}$ be a lifting of $-b$ (i.e.,
$\pi(z)=-b$) and $l'=  m +(-1)^i\debar z+ p^*n_1'+q^*n'_2 \in
KS_{X\times Y}{i-1}$. Then $l' \in (L )^{i-1}$; indeed $\pi(l')=
\pi( m) -(-1)^i\debar b+f_*a_1-f^*a_2=0$; moreover,  $\debar
l'=\debar m+ p^* n_1 +q^*n_2=l$.  Therefore, $
[(l,n_1,n_2,m)]=[\debar (l',n_1',n_2',(-1)^i z)]$ is zero in $
H^i(C^ \cdot_{(h,g)})$.

\emph{$\beta$ is surjective}. Let $[(a_1,a_2,b)] \in
H^i(B^\cdot)$. Then $\debar a_1=\debar a_2=0$ and $\debar b
+(-1)^i(f_*a_1-f^*a_2)=0$. Let $m \in KS_{X \times Y}^{i-1}$ be a
lifting of $(-1)^ib$, i.e., $\pi(m)=(-1)^i b$. Let $n_1= a_1\in
KS_X^i$ and $n_2=  a_2\in KS_Y^i$. Finally, let $l=\debar m +
p^*a_1+ q^*a_2 \in KS_{X \times Y}^{i}$. Then $l \in L ^i$; in
fact $\pi(l)= \pi (\debar m) + f_*a_1-f^*a_2 =(-1)^i  \debar b +
f_*a_1-f^*a_2 =0$.

Since $ \debar l=\debar n_1=\debar n_2=0$ and $ -\debar m
-p^*n_1-q^*n_2+l=0 $, $[(l,n_1, n_2,m)] \in H^i(C^ \cdot_{(h,g)})$
and $\beta[(l,n_1, n_2,m)]=[( n_1, n_2,(-1)^i
\pi(m))]=[(a_1,a_2,b)] \in H^i(B)$.
\end{proof}
\end{oss}

\subsection{$\Def_{(h,g)}$ is isomorphic to $\Def(f)$}
Using the notation above and diagram (\ref{equa diagrL-> KS(XxY)
<-KS(X)xKS(Y)}), consider the functor $\Def_{(h,g)}$. Since $h$ is
injective, by Remark~\ref{oss DEF_(h,g) con h iniettivo}, for each
$(A,m_A) \in \Art$ we have

$$
\Def_{(h,g)}(A)=\{(n,e^m) \in  (N^1 \otimes m_A)\times
\operatorname{exp}(M^0 \otimes m_A)  |
$$
$$
dn +\frac{1}{2}[n,n]=0,\  e^{-m}*g(n) \in L^1 \otimes m_A
\}/gauge.
$$

\begin{oss}\label{oss (n,e^m) DEF(h,g)-> def T in def Xx defY}
Let $(n,e^m) \in \Def_{(h,g)}$. In particular, $n=(n_1,n_2)$
satisfies the Maurer-Cartan equation and so $n_1 \in \MC_{KS_X}$
and $n_2 \in \MC_{KS_Y}$. Therefore, there are associated with $n$
infinitesimals deformations $X_A$ of $X$ (induced by $n_1$) and
$Y_A$ of $Y$ (induced by $n_2$). Moreover, since $g(n)$ satisfies
Maurer-Cartan equation in $M=KS_{X \times Y}$, it defines an
infinitesimal deformation $(X \times Y)_A$ of $ X \times Y$. By
construction, the deformation $(X \times Y)_A$ is the product of
the deformations $X_A$ and  $Y _A$.

\end{oss}

Let $i^*:\shA_{Z}^{0,*} \lrg \shA^{0,*}_\Gamma$ be the restriction
morphism and let $I=\ker i^* \cap \shO_{Z}$ be the holomorphic
ideal sheaf of the graph  $\Gamma $of $f$ in $Z$.

\medskip

Consider an infinitesimal deformation of the holomorphic map $f$
over $\Spec(A)$ as an infinitesimal deformation
$\widetilde{\Gamma}$ of $\Gamma$ over $\Spec(A)$ and
$\widetilde{Z}$ of $Z$ over  $\Spec(A)$, with $\widetilde{Z}$
product of deformations of $X$ and  $Y$ over $\Spec(A)$.

By applying  Remark~\ref{oss (n,e^m) DEF(h,g)-> def T in def Xx
defY} and Theorem~\ref{teo def_k =Def_X}, the condition on the
deformation $\widetilde{Z}$ is equivalent to requiring
$\shO_{\widetilde{Z}}=\shO_A(g(n))$, for some Maurer-Cartan
element $n\in  KS_X \times KS_Y$.

The deformation $\widetilde{\Gamma}$ of the graph corresponds to
an infinitesimal deformation $I_A\subset \shO_{\widetilde{Z}} $ of
the holomorphic ideal sheaf $I$ over $\Spec(A)$, that is, $I_A$ is
a sheaf flat over $A$ such that $I_A\otimes _A \C\cong I$.

In conclusion, to give an infinitesimal deformation of $f$ over
$\Spec(A)$ (an element in $\Def(f)(A)$) it is sufficient to give
an ideal sheaf $I_A\subset \shO_A(g(n))$ (for some $n\in \MC_{KS_X
\times KS_Y}$)  with $I_A$ $A$-flat and $I_A\otimes _A \C\cong I$.

\begin{teo}\label{teo Def_(h,g)  Def (f)}
Let $h,g$ and $i^*$ be   as above. Then there exists an
isomorphism of functors
$$
\gamma: \Def_{(h,g)}\lrg \Def (f).
$$
Given a local Artinian $\mathbb{C}$-algebra $A$ and  an element
$(n,e^m) \in \MC_{(h,g)}(A)$, we define a deformation of $f$ over
$\Spec(A)$ as a deformation $I_A(n,e^m)$ of the holomorphic ideal
sheaf of the graph  of $f$ in the following way
$$
\gamma(n,e^m)=I_A(n,e^m):=(\ker(\shA_{Z}^{0,0}\otimes A
\stackrel{\debar+\bl_{g(n)}}{\lrg} \shA^{0,1}_{Z}\otimes A))\cap
e^m (\ker i^* \otimes A)=
$$
$$
=\shO_{A}(g(n))\cap e^m(\ker i^* \otimes A),
$$
where $\shO_{A}(g(n))$ is the infinitesimal deformation of $Z$
that corresponds to $g(n)\in MC_{KS_{X \times Y}}$ (see
Theorem~\ref{teo def_k =Def_X}).
\end{teo}

\begin{proof}
For each $(n,e^m) \in \MC_{(h,g)}(A)$ we have defined
$$
I_A(n,e^m) =\shO_{A}(g(n))\cap e^m(\ker i^* \otimes A).
$$
First of all, we verify that this sheaf $I_A(n,e^m)\subset
\shO_{A}(g(n))$ defines an infinitesimal  deformation of $f$;
therefore, we need to prove that $I_A(n,e^m)$ is flat over $A$ and
$I_A\otimes _A \C\cong I$. It is equivalent to verify these
properties for $e^{-m}I_A(n,e^m)$.

Applying Lemma~\ref{lemma e^a(d+l_x)e^-a = d+e^a *x} yields
$$
e^{-m}(\shO_A(g(n)))=\ker(\debar + e^{-m} *g(n):
\shA_Z^{0,0}\otimes A \lrg \shA_Z^{0,1}\otimes A)
$$
and also
$$
e^{-m}I_A(n,e^m)=e^{-m}(\shO_A(g(n)))\cap (\ker i^* \otimes A)=
$$
$$
=\ker(\debar +e^{-m}*g(n))\cap(\ker i^* \otimes A).
$$
Since flatness is a local property, we can assume that $Z$ is a
Stein manifold. Then $H^1(Z,\Theta_Z)=0$ and $H^0(Z,\Theta_Z) \lrg
H^0(Z,N_{\Gamma|Z})$ is surjective. Since the following sequence
is exact
$$
\cdots \lrg H^0(Z,\Theta_Z) \lrg H^0(Z,N_{\Gamma|Z}) \lrg H^1(Z,L)
\lrg H^1(Z,\Theta_Z)\lrg \cdots,
$$
we conclude that  $H^1(L)=0$ or, equivalently, that the tangent
space of the functor $\Def_{L }$ is trivial. Therefore, by
Corollary~\ref{cor funtoreF banale sse t_F=0}, $\Def_{L }$ is the
trivial functor.

This implies the existence of $\nu\in L^0 \otimes m_A$ such that
$e^{-m}*g(n)=e^\nu*0$ (by hypothesis, $e^{-m}*g(n)$ is a solution
of  the Maurer-Cartan equation  in $L$). Moreover, we recall that
if $a \in L^0\otimes m_A$ then $e^a(\ker  i^* \otimes A)=\ker i^*
\otimes A$ (see Section~\ref{sezio submanifold defz L'}).

Therefore,
$$
e^{-m}I_A(n,e^m)=\ker(\debar +e^{\nu}*0)\cap(\ker i^* \otimes
A)=\shO_A(e^\nu *0)\cap (\ker i^* \otimes A)
$$
$$
=e^\nu(\shO_A(0))\cap e^\nu(\ker i^* \otimes A)=e^\nu(I\otimes A).
$$
Thus, $I_A(n,e^m)$  defines a deformation of $f$ and the morphism
$$
\gamma: \MC_{(h,g)}\lrg \Def (f)
$$
is well defined, such that
$$
\gamma(A): \MC_{(h,g)}(A)\lrg \Def (f)(A)\qquad
$$
$$
\qquad  \qquad (n,e^m) \longmapsto \gamma(n,e^m)=I_A(n,e^m).
$$

Moreover, $\gamma$ is  well defined on
$\Def_{(h,g)}(A)=\MC_{(h,g)}(A)/gauge$. Actually, for each $a \in
L^0 \otimes m_A$ and $b \in N^0\otimes m_A$,  we have
$$
\gamma(e^{b}*n,e^{g(b)}e^m e^{a})= \shO_A(e^{g(b)}*g(n)) \cap
e^{g(b)}e^m e^{a}(\ker i^* \otimes A)=
$$
$$
e^{g(b)}\shO_A(g(n)) \cap e^{g(b)}e^m (\ker i^* \otimes
A)=e^{g(b)}\gamma(n,e^m).
$$
This implies that the deformations $\gamma(n,e^m)$ and $\gamma(
e^b*n,e^{g(b)}e^m e^{a})$ are isomorphic (see Remark~\ref{oss defo
mapp=Def grafo prod defX.defY}).

In conclusion, $\gamma: \Def_{(f,g)} \lrg \Def(f)$ is a well
defined natural transformation of functors.

In order to prove that $\gamma$ is an isomorphism   it is
sufficient to prove that
\begin{itemize}
  \item[i)] $\gamma$ is injective;
  \item[ii)] $\gamma$ induces a bijective map on the tangent
  spaces;
  \item[iii)] $\gamma$ induces an injective map on the
  obstruction spaces.
\end{itemize}
Actually, by Corollary~\ref{cor etale se bietti su tg e inj su
ostru}, conditions $ii)$ and $iii)$ imply that $\gamma$ is
\'{e}tale and so surjective.

\medskip

i)  $\gamma\  is\  injective$. Suppose that $\gamma(n,e^m)=
\gamma(r,e^s)$, i.e., the deformations induced  by $(n,e^m)$ and
$(r,e^s)$, respectively, are isomorphic. We want to conclude that
$(n,e^m)$ is gauge equivalent to $ (r,e^s)$, i.e., there exist $a
\in L^0 \otimes m_A$ and $b \in N^0 \otimes m_A$ such that
$e^b*r=n$ and $e^{g(b)}e^se^a=e^m$.

By hypothesis, $\gamma(n,e^m)$ and $ \gamma(r,e^s)$ are isomorphic
deformations; then, in particular, the deformations induced on $Z$
are isomorphic. This implies that there exists $b \in N^0 \otimes
m_A$ such that $e^b*r=n$ and so $e^{g(b)}(\shO_A(g(r)) )
=\shO_A(g(n))$. Up to substituting $(r,e^s)$ with its equivalent
$(e^b*r,e^{g(b)}e^s)$, we can assume   to be  in the following
situation
$$
\shO_A(g(n))\cap e^m(\ker i^* \otimes A)=\shO_A(g(n))\cap e^{m'}
(\ker i^* \otimes A).
$$
Let $e^{a}=e^{-m'}e^{m}$, then
$$
e^{a}(e^{-m}(\shO_A(g(n)))\cap (\ker i^* \otimes
A))=e^{-m'}(\shO_A(g(n)))\cap (\ker i^* \otimes A).
$$
In particular, $e^a(e^{-m} (\shO_A(g(n)))\cap (\ker i^* \otimes
A)) \subseteq \ker i^* \otimes A$.

Next, we prove, by induction, that $a \in L^0 \otimes m_A$ (thus
$e^m= e^{m'}e^a=e^{g(b)}e^se^a$).

Let $z_1,\ldots,z_n$ be holomorphic coordinates on $Z$ such that
$Z\supset \Gamma=\{z_{t+1}=\cdots=z_n=0\}$. Consider the
projection on the residue field
$$
e^{-m}(\shO_A(g(n)))\cap (\ker i^* \otimes A )\lrg \shO_Z \cap
\ker i^*.
$$
Then $z_i\in \ker i^*\cap \shO_Z$, for $i >t$. Since
$e^{-m}(\shO_A(g(n)))\cap (\ker i^* \otimes A )$ is flat over $A$,
we can lift $z_i$ to $\tilde{z}_i=z_i + \varphi_i\in
e^{-m}(\shO_A(g(n)))\cap (\ker i^* \otimes A ) $, with $\varphi_i
\in \ker i^* \otimes m_A$. By hypothesis,
\begin{equation}\label{equa e^a(ztilde)=e^a(z)+e^a(phi)}
e^a(\tilde{z}_i)=e^a(z_i)+e^a(\varphi_i) \in \ker i^*\otimes A.
\end{equation}
By Lemma~\ref{lemm L manetti is DGLA}, in order to prove that
$a\in L^0 \otimes m_A$ it is sufficient to verify that
$e^a(z_i)\in \ker i^*\otimes A $ and so, by (\ref{equa
e^a(ztilde)=e^a(z)+e^a(phi)}), that $e^a(\varphi_i) \in \ker i^*
\otimes A $.

If $A=\C[\epsi]$, then $\varphi_i \in \ker i^* \otimes \C\epsi$, $
a \in \shA^{0,0}_Z \otimes \C\epsi$, this implies
$e^a(\varphi_i)=\varphi_i \in \ker i^* \otimes \C\epsi$.

Next, let $0 \lrg J \lrg B \stackrel{\alpha}{\lrg} A\lrg 0$ be a
small extension. By hypothesis, $\alpha(a) \in L^0 \otimes m_A$,
that is, $\displaystyle \alpha(a)=\sum_{j=1}^n \overline{a}_j
\frac{\de}{\de z_j}$ with $\overline{a}_j \in \ker i^* \otimes
m_A$ for $j
>t$.

Let $a'_j$ be liftings of $\overline{a}_j$. Then  $a'_j \in \ker
i^* \otimes m_B$ for $j >t$ , $\displaystyle a'= \sum_{j=1}^n a'_j
\frac{\de}{\de z_j} \in L^0 \otimes m_B $ and
$e^{a'}(\varphi_i)\in \ker i^* \otimes m_B$. Since
$\alpha(a)=\alpha(a')$, then $a=a'+j$ with $j \in M^0 \otimes J$.
This implies that $e^a(\varphi_i) =e^{a'+j} (\varphi_i)=
e^{a'}(\varphi_i) \in \ker i^* \otimes m_B$.

\bigskip

As to tangent and obstruction spaces, by Theorem~\ref{teo TANG e
OSTR DEF(f) horikawa} and Lemma~\ref{lemm C horikawa quasi iso to
B}, the tangent space of $\Def(f)$ is isomorphic to $H^1(B^\cdot)$
and the obstruction space  is naturally contained in
$H^2(B^\cdot)$, where $(B^\cdot,D_{\debar})$ is the complex with $
B^p=A_X^{(0,p)}( \Theta_X) \oplus A_Y^{(0,p)}(\Theta_Y) \oplus
A_X^{(0,p-1)}( f^*\Theta_Y) $ and $ D_{\debar} (x,y,z)=(\debar
x,\debar y, \debar z+ (-1)^p(f_*x-f^*y)) $, for each $(x,y,z) \in
B^p$. As to the functor $\Def_{(h,g)}$, in Section~\ref{sez
tangent e ostru MC(h,g) DEF(h,g)}, we have proved that the tangent
space is $H^1(\cil ^\cdot)$ and the obstruction space is naturally
contained in  $H^2(\cil^\cdot)$, where $\cil^\cdot$ is the
suspended cone associated with the pair $(h,g)$
(Section~\ref{sezione triplo cono}). Moreover,
Proposition~\ref{prop isomor complesso Namba con Cono} shows the
existence of a quasi-isomorphism $\beta$ between the previous
complexes $B^\cdot$ and $\cil ^\cdot$.

Next, we want to prove that the maps induced by $\gamma$ on
tangent and obstruction spaces coincide with the isomorphism
induced by $\beta$.

\bigskip

$ii)$ $\gamma $ induces a  bijection on the tangent spaces  (we
prove that $\gamma$ coincides with the isomorphism $\beta$).

Let $(n_1,n_2,m) \in \Def_{(h,g)}(\C[\epsi])$. Then $\debar
n_1=\debar n_2=0 $, $\debar m +g(n)=\debar m + p^* n_1 +q^* n_2
\in H^1(L )$ and so  $(n_1,n_2,m)$ determines the class
$[(g(n)+\debar m,n_1,n_2,m)] \in H^1(C_{(h,g)}^\cdot)$. Moreover,
we note that $\debar \pi (m) +f_*n_1 -f^*n_2=0$ and
$\beta[(g(n)+\debar m,n_1,n_2,m)]=[( n_1, n_2,-\pi(m))]$.

\bigskip

Next, consider the first order deformation of $f$ induced by
$(n_1,n_2,m)$, i.e., $\gamma(n_1,n_2,m)=\shO_{A}(g(n))\cap
e^m(\ker i^* \otimes A)$. Using Theorem~\ref{teo TANG e OSTR
DEF(f) horikawa} and Remark~\ref{oss class H^2(TX,TY,f*TY) def by
Xe,Ye}, we associate with $\gamma(n_1,n_2,m)$ the class
\footnote{By Theorem~\ref{teo def_k =Def_X}, the class associated
with the first order deformation of $X$ induced by $n_1$ is $n_1$
itself  and the class associated with the first order deformation
of $Y$ induced by $n_2$ is $n_2$ itself.} $[(n_1,n_2,a)] \in H^1
(B ^\cdot)$, such that $\debar a=f_* n_1-f^*n_2$. Then $\debar a=-
\debar \pi(m)$ and so $ [(n_1,n_2,a)]=[(n_1,n_2,-\pi(m))]= \beta
([(g(n)+\debar m,n_1,n_2,m)]) \in H^1 (B ^\cdot)$.

This implies that the map induced by $\gamma$ on the tangent
spaces coincides with the  isomorphism $\beta$.

\bigskip

$iii)$ $\gamma'$ induces an injective map on the obstruction
spaces. In this case too, we prove that  the map induced by
$\gamma$ on the obstruction space coincides with  $\beta$ and so
it is injective.

Actually, let
$$
0 \lrg J \lrg B \stackrel{\alpha}{\lrg} A\lrg 0
$$
be a small extension. The obstruction class associated with
$(n_1,n_2,m)\in \Def_{(h,g)}(A)$ is the class $[(k,h_1,h_2,r)] \in
H^2(\cil ^\cdot)\otimes J $ defined in \\ Lemma~\ref{lem calcolo
ostruzione MC_{(h,g)}=H^2(cil)}. We note that $\debar r +p^*h_1
+q^* h_2 \in H^2(L)$) and so $\pi(\debar r)=-f_*h_1 +f^* h_2 $.

Next, again by Theorem~\ref{teo TANG e OSTR DEF(f) horikawa} and
Remark~\ref{oss class H^2(TX,TY,f*TY) def by Xe,Ye}, the
obstruction class associated with the deformation
$\gamma(n_1,n_2,m)$ of $f$  induced by $(n_1,n_2,m)$  is
$[(h_1,h_2,a)] \in H^2(B^\cdot)\otimes J $ with $\debar a=-( f_*
h_1-f^*h_2)$ and as above $[(h_1,h_2,a)]= \beta ( [(k,h_1,h_2,r)])
\in H^2(B^\cdot)\otimes J$.

\end{proof}

In conclusion, by suitably choosing  $L,M,$ and $h,g$,
Theorem~\ref{teo Def_(h,g) Def (f)} shows that the infinitesimal
deformation functor $\Def(f)$ associated with a holomorphic map
$f$ is isomorphic to the functor ${\Def}_{(h,g)} $.

\begin{teo}\label{teo esiste hil governa DEF(f)}
Let $f:X \lrg Y$ be a  holomorphic map. Then the DGLA $\hil$
associated with the above morphisms $h:L \hookrightarrow
KS_{X\times Y}$ and $g=(p^*,q^*):KS_X \times KS_Y \lrg KS_{X\times
Y}$ (see Definition~\ref{dhef DGLA H(h,g)}) controls infinitesimal
deformations of $f$:
$$
\Def_{\hil} \cong \Def(f).
$$
\end{teo}
\begin{proof}
It is sufficient to apply  Theorem~\ref{teo Def_(h,g) Def (f)},
Corollary~\ref{cor Mpositiva exteDEF(h,g)restrict toDEf(h,g)} and
Theorem~\ref{teo EXT DEF_H=DEF(h,g)}.
\end{proof}

\begin{oss}
This theorem  is very interesting from the point of view of \lq
\lq guiding principle", since it shows the existence of a DGLA
which controls the geometric problem  of infinitesimal
deformations of holomorphic maps.

Anyway, in  Chapter~\ref{capitolo application and examples} we
will see that in the applications it is more convenient  to use
the functor ${\Def}_{(h,g)} $  than $\Def_{\hil}$.
\end{oss}

\begin{oss}
Consider the diagram
\begin{center}
$\xymatrix{  & & L \ar[d]^h   & \\
KS_{X} \times  KS_Y \ar[rr]^{g}  \ar@/^/[rrrd]_{\pi \circ g}& &
KS_{X \times Y}\ar[dr]^\pi &    \\
& & & A_X^{0,*}(f^*T_Y). &
 \\ }$
\end{center}

Since $h$ is injective, Lemma~\ref{lem h inj then C_(h,g)=C_pi
dopo g} implies the existence of    a quasi-isomorphism of
complexes $(\cil^\cdot,D)$ and $(C^\cdot_{\pi \circ
g},\check{\delta})$.

Therefore, we get the following exact sequence
\begin{equation}\label{succ esatta mia con ro uguale ziv}
\cdots \lrg H^1(C^\cdot_{\pi \circ g} )\stackrel{\varrho^1}{\lrg}
H^1(X,\Theta_X)\oplus H^1(Y,\Theta_Y)\lrg H^1(X,f^*\Theta_Y)\lrg
\end{equation}
$$
\qquad \lrg H^2(C^\cdot_{\pi \circ g} )\stackrel{\varrho^2}{\lrg}
H^2(X,\Theta_X)\oplus H^2(Y,\Theta_Y)\lrg H^2(X,f^*\Theta_Y)\lrg
\cdots
$$
where $\varrho^1$ and $\varrho^2$ are the projections on the first
factors and they are induced by the projection morphism
$\varrho:\Def(f) \lrg \Def_{KS_{X} \times  KS_Y }$ (see
Remark~\ref{oss proj varrho induce no EXT Def(h,g)->Def_N}).

In particular, $\varrho:\Def(f) \lrg \Def_{KS_{X} \times  KS_Y }$
associates with an infinitesimal deformation of $f$ the induced
infinitesimal deformations of $X$ and  $Y$.

Then $\varrho^1$ associates with a first order deformation  of $f$
the induced first order deformations of $X$ and $Y$ and
$\varrho^2$ is a morphism of obstruction theory: the obstruction
to deform $f$ is mapped to the induced obstructions to deform $X$
and $Y$ (see also Remark~\ref{oss (n,e^m) DEF(h,g)-> def T in def
Xx defY}).

\end{oss}

In \cite{bib RAn mappe}, Z. Ran studied the infinitesimal
deformations of a holomorphic map $f:X \lrg Y$ of singular compact
complex spaces. He introduced some algebraic objects $T^i_f$,
$i=1,2$, that classify the deformations of a map $f$ and obtained
the following exact sequence
\begin{equation}\label{successione esatta ziv}
T^1_f \lrg \Ext^1_{\shO_X}(\Omega_X,\shO_X) \oplus
\Ext^1_{\shO_Y}(\Omega_Y,\shO_Y)  \lrg
\Ext^1_{\shO_X}(f^*\Omega_Y, \mathcal{O}_X)\lrg
\end{equation}
$$
 \lrg T^2_f \lrg \Ext^2_{\shO_X}(\Omega_X,\shO_X) \oplus
\Ext^2_{\shO_Y}(\Omega_Y,\shO_Y)   \lrg
\Ext^2_{\shO_X}(f^*\Omega_Y, \mathcal{O}_X).
$$

\begin{lem}
If $X$ and $Y$ are compact complex manifold, then the exact
sequence (\ref{successione esatta ziv}) reduces to (\ref{succ
esatta mia con ro uguale ziv}).
\end{lem}
\begin{proof}
If $X$ and $Y$ are smooth, then $\Omega_X$ and $\Omega_Y$ are
locally free. Then applying the spectral sequence associated with
$\Ext$ (see \cite[Lemme~7.4.1]{bib Godement}) we get
$\Ext^i_{\shO_X}(\Omega_X,\shO_X) \cong H^i(X,\Theta_X)$,
$\Ext^i_{\shO_Y}(\Omega_Y,\shO_Y) \cong H^i(Y,\Theta_Y)$ and
$\Ext^i_{\shO_X}(f^*\Omega_Y,\shO_X) \cong H^i(X,f^*\Theta_Y)$.

\end{proof}

\chapter{Semiregularity maps}\label{capitolo application and examples}

In the previous  chapter  we have studied the infinitesimal
deformations of holomorphic maps.

More precisely, let $f:X \lrg Y$ be a  holomorphic map of compact
complex manifolds and $\Gamma\subset X \times Y$ its graph.

Consider the DGLAs $M=KS_{X\times Y}$, $N=KS_X \times KS_Y$ and
the morphism  $g=(p^*,q^*): KS_X \times KS_Y\lrg KS_{X\times Y}$,
where $p$ and $q$ are  the projections  of $X \times Y$ onto $X$
and $Y$, respectively.

Moreover, let $L=A_{X \times Y}^{0,*}(-log\, \Gamma)$ be the DGLA
defined by the following exact sequence (see Section~\ref{sez teo
Def_(h,g)Def (f)}):
$$
0 \lrg L\stackrel{h}{ \lrg} KS_{X\times Y} \stackrel{\pi}{\lrg}
 A_X^{0,*}(f^*\Theta_Y) \lrg 0.
$$
Then we get the following diagram
\begin{equation}\label{diagram L M N in chapter semiregolar}
\xymatrix{& & L \ar[d]^h   \\
KS_X \times KS_Y\ar[rr]^{(p^*,q^*)} &  & KS_{X\times Y}. \\ }
\end{equation}

Theorem~\ref{teo Def_(h,g)  Def (f)} shows the existence of an
isomorphism between the functor $\Def(f)$ of infinitesimal
deformations of $f$ and the functor $\Def_{(h,g)}$ associated with
the pair of morphisms $(h,g)$:
$$
\Def(f) \cong \Def_{(h,g)} \qquad \mbox{(Theorem~\ref{teo
Def_(h,g) Def (f)}).}
$$

Furthermore, Theorem~\ref{teo esiste hil governa DEF(f)} gives an
explicit description of the DGLA $H_{(h,g)}$ (Definition~\ref{dhef
DGLA H(h,g)}) that controls the infinitesimal deformations of $f$:
$$
\Def(f) \cong \Def_{H(h,g)} \qquad \mbox{(Theorem~\ref{teo esiste
hil governa DEF(f)}).}
$$

\medskip

In particular, if $i:X \hookrightarrow Y$ is an inclusion, we can
find an easy description of the DGLA associated with $\Def(i)$
(see Section~\ref{sezio inclusion}). Actually, let $L'$ be the the
DGLA $L'$ introduced in Section~\ref{sezio submanifold defz L'}
(see also \cite[Section~5]{bib ManettiPREPRINT})
$$
0 \lrg L' \lrg A_Y^{0,*}(\Theta_Y) \stackrel{\pi'} \lrg
A_X^{0,*}(N_{X|Y})\lrg 0.
$$
Thus, $L'$ controls the deformations of  $i$ (Corollary~\ref{cor L
governa def coppia (XcY)}).

\medskip

In general, without an easy description of $H_{(h,g)}$ it is
convenient to use the deformation  functor associated with the
previous  diagram (\ref{diagram L M N in chapter semiregolar}).

For example, if we want to study the infinitesimal deformations of
$f$ with fixed domain or fixed codomain, it is sufficient to
consider diagram (\ref{diagram L M N in chapter semiregolar}) with
$ N=KS_Y$ (Section~\ref{sezione fixed target})  or $N=KS_X$
(Section~\ref{sezione fixed domain}), respectively.

\bigskip

Anyway, the main application of the techniques developed in the
previous chapters concerns  the study of obstructions to deforming
a holomorphic map $f$ and the  \lq \lq semiregularity" maps.

In general, we can find a vector space $V$ (most of the times a
cohomology vector space) that contains the obstruction space, but
we do not  know an explicit description of the elements which are
obstructions.

Then the idea is to restrict the vector space $V$, for example by
defining a map on $V$, the so called \lq \lq semiregularity" map,
that contains the obstructions in the kernel.

In particular, let $f:X \lrg Y$ be a  holomorphic map and consider
the   infinitesimal deformations of $f$ with fixed codomain $Y$.
In  \cite{bib HorikawaI} E.~Horikawa   proved the following
theorem.

\begin{teo} (Horikawa)\label{teo horikawa}
Let $f:X \lrg Y$ be a  holomorphic map and consider the functor of
infinitesimal deformations of $f$ with fixed codomain $Y$.   Then
the tangent space is the hypercohomology vector space
$\mathbb{H}^1\left(X,\shO(\Theta_X)
\stackrel{f_*}{\lrg}\shO(f^*\Theta_Y)\right)$ and the obstruction
space is contained in $\mathbb{H}^2\left(X,\shO(\Theta_X)
\stackrel{f_*}{\lrg}\shO(f^*\Theta_Y)\right)$.
\end{teo}

Using the techniques introduced in the previous section,  we can
improve this result in the case of K\"{a}hler manifolds, defining
a map that contains the obstruction space in the kernel.

\begin{teos}
Let $f:X \lrg Y$ be a  holomorphic map of  compact K\"{a}hler
manifolds. Let $p=dim Y - dim X$. Then the obstruction space of
deformations of $f$ (with fixed $Y$) is contained in the kernel of
a map
$$
\sigma : \mathbb{H}^2\left(X,\shO(\Theta_X)
\stackrel{f_*}{\lrg}\shO(f^*\Theta_Y)\right) \lrg
H^{p+1}(Y,\Omega_Y^{p-1}).
$$
\end{teos}
The previous map is the generalization of the semiregularity map
defined by Bloch (see \cite{bib bloch} or \cite[Section~9]{bib
ManettiPREPRINT}) obtained when  $f$ is the inclusion map
$X\hookrightarrow Y$, i.e.,
$$
\sigma: H^1(X,N_{X|Y}) \lrg H^{p+1}(Y,\Omega_Y^{p-1}).
$$

The proof of this theorem is postponed to
Section~\ref{sottosection semiregolarity} (Corollary~\ref{coroll
semiregolari-ostruzion}), where we also give an explicit
description of the map $\sigma$.

\section{Deformations  with fixed codomain}
\label{sezione fixed target} This section is devoted to   studying
infinitesimal deformations of a holomorphic map $f:X \lrg Y$, with
fixed codomain $Y$.

In this case the DGLA $N$ reduces to $KS_X$ and so diagram
(\ref{diagram L M N in chapter semiregolar}) reduces to

\begin{center}
$\xymatrix{ & L \ar[d]^h   & \\
A_X^{0,*}(\Theta_X) \ar[r]^{p^*}  \ar@/^/[rrd]_{f_*} &
A_{X\times Y}^{0,*}(\Theta_{X\times Y})\ar[dr]^\pi &    \\
& & A_X^{0,*}(f^*\Theta_Y). &
 \\ }$
\end{center}
where $f_*$ is the product $\pi \circ p^*$ (see Section~\ref{sez
teo Def_(h,g)Def (f)}).


Using this diagram and Theorem~\ref{teo Def_(h,g)  Def (f)}, we
can easily prove  Theorem~\ref{teo horikawa} due to E.~Horikawa.

\begin{prop}\label{prop mia uguale horikawa}
The tangent space of the infinitesimal deformation functor of a
holomorphic map $f:X \lrg Y$, with $Y$ fixed, is
$H^1(C^\cdot_{f_*})$ and the obstruction space is naturally
contained in $H^2(C^\cdot_{f_*})$.
\end{prop}
\begin{proof}
Theorem~\ref{teo Def_(h,g)  Def (f)} implies that the
infinitesimal deformation  functor of $f$,  with $Y$ fixed, is
isomorphic to $\Def_{(h,p^*)}$. Therefore, the tangent space is
$H^1(C^\cdot_{(h,p^*)})$ and the obstruction space is naturally
contained in $H^2(C^\cdot_{(h,p^*)})$. Since $h$ is injective,
Lemma~\ref{lem h inj then C_(h,g)=C_pi dopo g} implies that, for
each $i$, $H^i (C^\cdot_{(h,p^*)})\cong H^i(C^\cdot_{\pi \circ
p^*})=H^i(C^\cdot_{f_*})$.
\end{proof}

\medskip

As we already announced  we improve this theorem in the case of
K\"{a}hler manifolds (see   Section~\ref{sottosection
semiregolarity}).

The next section is devoted to some preliminary lemmas.

\subsection{Preliminaries } \label{sezio HTP}

Let $Z$ be a complex manifold.

Then $KS_Z=A_Z^{0,*}(\Theta_Z)$ is the Kodaira-Spencer algebra of
$Z$ and in Section~\ref{sec cotrction map e Lie deriv}  we have
defined the contraction map $\bi$:
$$
\bi : KS_Z \lrg \operatorname \Hom^*(A_Z,A_Z),
$$
$$
\bi_a(\omega)=a \contr \omega
$$
for any $a \in KS_Z $ and $ \omega \in  A_Z^{*,*}$.

\noindent Therefore, $\bi(A_Z^{0,j}(\Theta_Z))\subset \oplus_{h,l}
\Hom^0(A_Z^{h,l},A_Z^{h-1,l+j}) \subset \Hom^{j-1}(A_Z,A_Z)$.

\medskip

In order to interpret $\bi$ as a morphism of DGLAs, the key idea,
due to M.~Manetti \cite[Section~8]{bib ManettiPREPRINT}, is to
substitute $\Hom^*(A_Z,A_Z)$ with the graded  vector space $\Htp
\left (\ker(\de), \dfrac{A_Z} {\de A_Z}\right)=\bigoplus_i
\Hom^{i-1} \left (\ker(\de), \dfrac{A_Z} {\de A_Z}\right)$ (see
Example~\ref{exe definizio Htp(V,W)}). Consider on $\Htp \left
(\ker(\de), \dfrac{A_Z} {\de A_Z}\right)$  the following
differential $\delta$ and bracket $\{\ ,\ \}$:
$$
\delta(f)=-\debar f - (-1)^{\deg(f)}f \debar,
$$
$$
\{f,g\}=f \de g -(-1)^{\deg(f)\deg(g)}g \de f.
$$

\begin{prop}
$\Htp \left (\ker(\de), \dfrac{A_Z} {\de A_Z}\right)$ is a DGLA
and the linear map
$$
\bi: A_Z^{0,*}(T_Z) \lrg \Htp\left(\ker(\de),\dfrac{A_Z}{\de
A_Z}\right)
$$
is a morphism of DGLAs.
\end{prop}
\begin{proof}
See \cite[Proposition~8.1]{bib ManettiPREPRINT}. An easy
calculation shows that the graded vector space $\Htp \left
(\ker(\de), \dfrac{A_Z} {\de A_Z}\right)$ is a DGLA.

Moreover, $\bi$ is a linear map that preserves degree and commutes
with differential and bracket. Actually, using Cartan
fomulas\footnote{$ \bi_{\tilde{d}a}=-[\debar, \bi_a],  \
\bi_{[a,b]}= [\bi_a,[\de,\bi_b]]=[[\bi_a,\de],\bi_b],\
[\bi_a,\bi_b]=0. $ } (see Lemma~\ref{lem formule cartan}), we get
$$
\bi_{\tilde{d}a}=-[\debar, \bi_a]=-\debar \bi_a
+(-1)^{\overline{a}-1}\bi_a \debar
$$
and, by definition  of $\delta$,
$$
\delta(\bi_a)=-\debar\bi_a- (-1)^{\overline{a}}\bi_a \debar.
$$
Therefore, $ \bi_{\tilde{d}a}=\delta(\bi_a)$. As  to the bracket,
again by  Cartan fomulas, we have
$$
\bi_{[a,b]}= [\bi_a,[\de,\bi_b]]=[\bi_a,\de \bi_b
-(-1)^{\overline{b}-1}\bi_b \de]=
$$
$$
 \bi_a\de \bi_b
-(-1)^{\overline{b}-1} \bi_a\bi_b \de
-(-1)^{(\overline{a}-1)\overline{b}}(\de \bi_b \bi_a
-(-1)^{\overline{b}-1}\bi_b \de \bi_a )=\footnote{$\bi_a\bi_b \de
$ and  $ \de \bi_b \bi_a$ are zero in $\Htp \left (\ker(\de),
\dfrac{A_Z} {\de A_Z}\right)$ }
$$
$$
\bi_a\de \bi_b +(-1)^{\overline{a}\overline{b} -1} \bi_b \de \bi_a
$$
and, by definition  of $\{\ ,\ \}$,
$$
\{\bi_a,\bi_b\}=\bi_a \de \bi_b
-(-1)^{\overline{a}\,\overline{b}}\bi_b \de \bi_a.
$$
\end{proof}

Next, let $f:X \lrg Y$ be a  holomorphic map, fix $Z=X\times Y$
and $\Gamma$ the graph of $f$ in $Z$. Let $I_\Gamma \subset A_Z$
be the space of the differential forms vanishing on $\Gamma$ and
$L \subset KS_Z$ be the DGLA defined as in Lemma~\ref{lemma
Lfascio is DGLA}:
$$
0 \lrg L \lrg KS_Z \lrg A_X^{0,*}(f^*\Theta_Y) \lrg 0.
$$
We recall that
$$
L \subset \{a \in  A^{0,*}_Z(\Theta_Z)| \ \bi_a(I_\Gamma)\subset
I_\Gamma\};
$$
moreover,
$$
p^*A_X^{0,*}(\Theta_X)\subset \{a \in  A^{0,*}_Z(\Theta_Z)| \
\bi_a(q^*A_Y)=0\},
$$
where $p$ and $q$ are the projection of $Z$ onto $X$ and $Y$,
respectively.

In conclusion, we can define the following commutative diagram of
morphisms of DGLAs:
\begin{center}
\begin{equation}\label{equa diagramma tre piani}
 \xymatrix{ L\ar@{^{(}->}[d]^h \ar[r] & K= \left\{f\in
\Htp\left(\ker(\de),\dfrac{A_Z}{\de A_Z}\right)\mid
f(I_\Gamma\cap\ker(\de))\subset \dfrac{I_\Gamma}{I_\Gamma\cap\de
A_Z}\right\}
\ar@{^{(}->}[d]^\eta \\
 A^{0,*}_Z(\Theta_Z)\ar[r]& \Htp\left(\ker(\de),\dfrac{A_Z}{\de
A_Z}\right) \\
A^{0,*}_X(\Theta_X)\ar[u]^{p^*} \ar[r] & J=\left\{f\in
\Htp\left(\ker(\de),\dfrac{A_Z}{\de A_Z}\right)\mid
f(\ker(\de)\cap q^* A_Y)=0 \right\}, \ar\ar@{^{(}->}[u]_\mu \\ }
\end{equation}
\end{center}
where the horizontal morphisms are all given by $\bi$.

Therefore, diagram (\ref{equa diagramma tre piani}) induces a
morphism of deformation functors:
$$
\mathcal{I}:\Def_{(h,p^*)}\lrg \Def_{(\eta,\mu)}.
$$

\begin{lem}\label{lemma ostruzione nel ker tra H^2}
If the differential graded vector spaces $(\de A_Z,\debar)$, $(\de
A_\Gamma,\debar)$ and $(\de A_Z \cap q^*A_Y, \debar)$ are acyclic,
then the functor $\Def_{(\eta,\mu)}$ is unobstructed. In
particular, the obstruction space of $\Def_{(h,p^*)}$ is naturally
contained in the kernel of the map
$$
H^2(C^\cdot_{(h,p^*)}) \stackrel{\mathcal{I}}{\lrg}
H^2(C^\cdot_{(\eta,\mu)}).
$$
\end{lem}
\begin{proof}This proof is an extension of the proof of
\cite[Lemma~8.2]{bib ManettiPREPRINT}.

The projection $\ker(\de)\to \ker(\de)/\de A_Z$ induces a
commutative diagram
\begin{center}
\begin{equation}\label{eqa diagr 2 colon DEF_(n,m) abeliano}
\xymatrix{ K\ar[d]^\eta & \{f\in K | f(\de A_Z)=0\}\ar[d]_{\eta'}
\ar[l]^\alpha  \\
\Htp\left(\ker(\de),\dfrac{A_Z}{\de A_Z}\right)  &
\Htp\left(\dfrac{\ker(\de)}{\de
A_Z},\dfrac{A_Z}{\de A_Z}\right)\ar[l]^\beta\\
J \ar[u]_\mu & \{f\in J| f(\de A_Z)=0\}.\ar[u]^{\mu'}
\ar[l]^\gamma
\\ }
\end{equation}
\end{center}
Since $\de A_Z$ is acyclic, $\beta$ is a quasi-isomorphism of
DGLAs. Since
$$
\coker(\alpha)=\{ f \in \Htp\left(\de A_Z,\dfrac{A_Z}{\de
A_Z}\right) | f(I_\Gamma \cap \de A_Z)\subset
\dfrac{I_\Gamma}{I_\Gamma\cap\de A_Z}  \},
$$
there exists an exact sequence
$$
0\to \Htp\left(\dfrac{\de A_Z}{I_\Gamma\cap \de
A_Z},\dfrac{A_Z}{\de A_Z}\right) \to \coker(\alpha)\to
\Htp\left(I_\Gamma\cap \de A_Z,\dfrac{I_\Gamma}{I_\Gamma\cap \de
A_Z}\right)\to 0.
$$
Moreover, the exact sequence
$$
0 \lrg I_\Gamma \cap A_Z \lrg \de A_Z \lrg \de A_\Gamma\lrg 0
$$
implies that $I_\Gamma \cap A_Z$ and $\dfrac{\de A_Z}{I_\Gamma\cap
\de A_Z}=\de A_\Gamma$ are acyclic. By Example~\ref{exe dhef
Hom^*(V,W) grade vect spac}, this implies that the complexes
$\Htp\left(\dfrac{\de A_Z}{I_\Gamma\cap \de A_Z},\dfrac{A_Z}{\de
A_Z}\right)$ and $\Htp\left(I_\Gamma\cap \de
A_Z,\dfrac{I_\Gamma}{I_\Gamma\cap \de A_Z}\right)$ are acyclic and
so the same holds for $\coker(\alpha)$; i.e., $\alpha$ is a
quasi-isomorphism.

As to  $\gamma$, we have
$$
\coker(\gamma)=
$$
$$
\{ f \in \Htp\left(\de A_Z,\dfrac{A_Z}{\de A_Z}\right) |\  f(\de
A_Z \cap q^*A_Y)=0 \}=
$$
$$
 \Htp\left(\dfrac{\de A_Z}{\de A_Z \cap
q^*A_Y},\dfrac{A_Z}{\de A_Z}\right).
$$
By hypothesis, $\de A_Z \cap q^*A_Y $ and $\de A_Z$ are acyclic
and so the same holds for $\dfrac{\de A_Z}{\de A_Z \cap q^*A_Y}$.
Then $\coker(\gamma)$ is acyclic, i.e.,   $\gamma$ is also a
quasi-isomorphism.

Therefore,  Theorem~\ref{teo NO exte quasi iso C(h,g)-C(n,m)then
DEf iso DEF} implies the existence of  an isomorphism of functors
$\Def_{(\eta,\mu)}\cong \Def_{(\eta',\mu')}$.

We note that the elements of the three algebras of the second
column of (\ref{eqa diagr 2 colon DEF_(n,m) abeliano}) vanish on
$\de A_Z$. Then, by definition  of the bracket $\{ \ , \ \}$,
these algebras are abelian. Therefore, Lemma~\ref{lemma L M N
abelian ->DEF(h,g)liscio} implies that the functor
$\Def_{(\eta,\mu)} \cong \Def_{(\eta',\mu')}$ is smooth. Finally,
Proposition~\ref{prop DEF_>DEF liscio,ostruz in ker tra H^2 }
guarantees that the obstruction space lies in the kernel of $
H^2(C^\cdot_{(h,p^*)}) \stackrel{\mathcal{I}}{\lrg}
H^2(C^\cdot_{(\eta,\mu)}). $
\end{proof}

\subsection{Semiregularity for deformations with fixed codomain}
\label{sottosection semiregolarity}

With the notation of the previous section we have the following
theorem.
\begin{teo}\label{teo ostruzi in ker H^2}
Let $f:X \lrg Y$ be a  holomorphic map of compact K\"{a}hler
manifolds. Then the obstruction space to the infinitesimal
deformations of $f$ with fixed codomain is contained in the kernel
of the following map
$$
H^2(C^\cdot_{f_*})\stackrel{\mathcal{I\,'}}{\lrg}
H^1\left(\Htp(I_\Gamma \cap \ker(\de)\cap q^*A_Y,A_\Gamma)\right).
$$
\end{teo}
\begin{proof}
Lemma~\ref{lemma kaler aciclico  q^* A_Y} implies that the
complexes $(\de A_Z,\debar)$, $(\de A_\Gamma,\debar)$ and $(\de
A_Z \cap q^*A_Y, \debar)$ are acyclic. Then we can apply
Lemma~\ref{lemma ostruzione nel ker tra H^2} to conclude that  the
obstruction space lies in the kernel of the following map
$$
H^2(C^\cdot_{(h,p^*)}) \stackrel{\mathcal{I }}{\lrg}
H^2(C^\cdot_{(\eta,\mu)}).
$$
Since $h$ is injective, by Lemma~\ref{lem h inj then C_(h,g)=C_pi
dopo g}, $H^2(C^\cdot_{(h,p^*)})\cong H^2(C^\cdot_{\pi \circ p^*})
\cong H^2(C^\cdot_{f_*})$. Then the obstructions lie  in the
kernel of $\mathcal{I} : H^2(C^\cdot_{f_*}) \lrg
H^2(C^\cdot_{(\eta,\mu)})$.

As to $H^2(C^\cdot_{(\eta,\mu)})$, consider
$$
K= \left\{f\in \Htp\left(\ker(\de),\dfrac{A_Z}{\de A_Z}\right)\mid
f(I_\Gamma\cap\ker(\de))\subset \dfrac{I_\Gamma}{I_\Gamma\cap\de
A_Z}\right\}
$$
and the exact sequence
$$
0 \lrg K \stackrel{\eta}{\lrg} \Htp\left(\ker
(\de),\dfrac{A_Z}{\de A_Z}\right)
\stackrel{\pi'}{\lrg}\coker(\eta)\lrg 0,
$$
with $\coker(\eta)= \Htp\left(I_\Gamma\cap \ker(\de)
,\dfrac{A_\Gamma}{ \de A_\Gamma}\right)$.

Applying again   Lemma~\ref{lem h inj then C_(h,g)=C_pi dopo g},
there exists an isomorphism \\ $ H^2(C^\cdot_{(\eta,\mu)}) \cong
H^2(C^\cdot_{\pi' \circ \mu}) $. By Remark~\ref{oss ho mappa da
C_h ->coker (h)} there also exists   a map $\mathcal{I\,
''}:H^2(C^\cdot_{\pi' \circ \mu}) \lrg H^1(\coker(\pi' \circ
\mu))$.

More precisely,
$$
J= \left\{f\in \Htp\left(\ker(\de),\dfrac{A_Z}{\de A_Z}\right)\mid
f(\ker(\de)\cap q^* A_Y)=0 \right\},
$$
$$
\pi' \circ \mu: J \lrg \Htp\left(I_\Gamma\cap \ker(\de)
,\dfrac{A_\Gamma}{ \de A_\Gamma}\right)
$$
and
$$
\coker(\pi' \circ \mu)=\Htp(I_\Gamma \cap \ker (\de)\cap q^* A_Y,
\dfrac{A_\Gamma}{\de A_\Gamma}).
$$

Finally, since the complex $\de A_\Gamma$ is acyclic, the
projection
$$
H^1(\Htp(I_\Gamma \cap \ker (\de)\cap q^* A_Y,\dfrac{A_\Gamma}{\de
A_\Gamma}) \lrg H^1(\Htp(I_\Gamma \cap \ker (\de)\cap q^* A_Y,
A_\Gamma))
$$
is an isomorphism.

Therefore, the obstruction space is contained in the kernel of
$\mathcal{I\, '}: H^2(C^\cdot_{f_*}) \lrg  H^1(\Htp(I_\Gamma \cap
\ker(\de)\cap q^*A_Y,A_\Gamma))$, i.e.,
\begin{center}
$\xymatrix{ H^2(C^\cdot_{f_*}) \ar[rrd]^{\mathcal{I\, '}}
\ar[r]^{\mathcal{I}}  & H^2(C^\cdot_{\pi' \circ
\mu})\ar[r]^{\mathcal{I\, ''}\qquad \qquad\ \  } &
H^1\left(\Htp(I_\Gamma \cap \ker (\de)\cap q^* A_Y,
\dfrac{A_\Gamma}{\de A_\Gamma}) \right)\ar[d]^\cong  \\
& & \oplus_i \Hom(H^i(I_\Gamma \cap \ker (\de)\cap q^* A_Y),
H^i(A_\Gamma)).
 \\ }$
\end{center}

\end{proof}

\begin{cor}\label{coroll semiregolari-ostruzion}
Let $f:X \lrg Y$ be a  holomorphic map of  compact K\"{a}hler
manifolds. Let $p=dim Y - dim X$. Then the obstruction space to
the infinitesimal deformations of $f$  with fixed $Y$  is
contained in the kernel of the map
$$
\sigma : \mathbb{H}^2\left(X,\shO(\Theta_X) \stackrel{f_*}
{\lrg}\shO(f^*\Theta_Y)\right) \lrg H^{p+1}(Y,\Omega_Y^{p-1}).
$$
\end{cor}

\begin{proof}
Let $n=dim X$, $p=dim Y - dim X$  and $\mathcal{H}$ be the space
of harmonic forms on $Y$ of type $(n+1,n-1)$. By Dolbeault's
theorem and Serre's duality we obtain the equalities
$\mathcal{H}^\nu=(H^{n-1} (Y,\Omega_Y^{n+1}))^\nu=
H^{p+1}(Y,\Omega_Y^{p-1})$.

Let $\omega \in \mathcal{H}$ such that $f^*\omega=0$. By
considering the contraction with $\omega$ we define a morphism of
complexes
$$
(A_X^{0,*}(f^*\Theta_Y), \debar) \stackrel{\contr \omega}{\lrg}
(A_X^{n,*+n-1}, \debar ),
$$
$$
\contr \omega( \phi f^* \chi)=\phi f^*(\chi \contr \omega)\in
A_X^{n,p+n-1}  \qquad \forall\  \phi f^* \chi \in
A_X^{0,p}(f^*\Theta_Y).
$$
Actually, since $\debar \omega=0$, then $\debar (\phi f^*(\chi
\contr \omega))=(\debar  \phi)f^*(\chi \contr \omega)= \contr
\omega(\debar  \phi f^* \chi)$.

In particular, since $f^* \omega=0$, using the identity  of
Lemma~\ref{lemma f^*(x-|w)=f_*x-|f^*w}\footnote{$f^*(\chi \contr
\omega)= \eta \contr f^* \omega$, for each $\omega \in A_Y^{*,*}$,
$\chi \in A_Y^{0,*}(\Theta_Y)$ and $\eta \in A_X^{0,*}(\Theta_X)$
such that $f^*\chi=f_*\eta$.}, we have the   following commutative
diagram
\begin{center}
$\xymatrix{A_X^{0,*}(f^*\Theta_Y)  \ar[r]^{\ \ \contr \omega}
& A_X^{n,*+n-1} \\
A_X^{0,*}(\Theta_X) \ar[u]^{f_*} \ar[r]  & 0 \ar[u]_\alpha.   \\
}$
\end{center}
Then we get  a morphism between the second cohomology groups of
the cones associated with the morphisms $f_*$ and $\alpha$:
$$
H^2(C^\cdot_{f_*}) \lrg H^2(C^\cdot_\alpha)\cong H^{n}(X,\Omega^n_
X).
$$
By taking the product of the previous morphism with the
integration on $X$, we get
$$
\sigma : \mathbb{H}^2\left(X,\shO(\Theta_X)
\stackrel{f_*}{\lrg}\shO(f^*\Theta_Y)\right) \lrg
H^{p+1}(Y,\Omega_Y^{p-1}).
$$
Since $ q^* \mathcal{H}$ is contained in $I_\Gamma \cap \ker
\debar \cap \ker \de \cap q^* A_Y$, we conclude the proof applying
Theorem~\ref{teo ostruzi in ker H^2}.

\end{proof}

\begin{oss}
We recall that, as already observed in Remark~\ref{oss de-debar
lemma per X,Y,XxY, Gamma}, the   K\"{a}hler hypothesis is just
used to have the   $\de \debar$-lemma  on $A_X$,$A_Y$, $A_{X
\times Y}$ and $A_\Gamma$.
\end{oss}

\section{Deformations of a map with fixed domain}
\label{sezione fixed domain}

In this section  we study  infinitesimal deformations of a
holomorphic map $f:X \lrg Y$,  with fixed domain  $X$.

In this case the DGLA $N$ reduces to $KS_Y$ and so diagram
(\ref{diagram L M N in chapter semiregolar}) reduces to

\begin{center}
$\xymatrix{ & L \ar[d]^h   & \\
A_Y^{0,*}(\Theta_Y) \ar[r]^{q^*}  \ar@/^/[rrd]_{-f^*} &
A_{X\times Y}^{0,*}(\Theta_{X\times Y})\ar[dr]^\pi &    \\
& & A_X^{0,*}(f^*\Theta_Y). &
 \\ }$
\end{center}
where $-f^*$ is the product $\pi \circ q^*$.


Using this diagram and Theorem~\ref{teo Def_(h,g)  Def (f)},
analogously to the case of fixed codomain we can prove the
following proposition.

\begin{prop}
The tangent space of the infinitesimal deformation functor of a
holomorphic map $f:X \lrg Y$, with $X$ fixed, is
$H^1(C^\cdot_{f^*})$ and the obstruction space is naturally
contained in $H^2(C^\cdot_{f^*})$.
\end{prop}
\begin{proof}
Analogous to the proof of Proposition~\ref{prop mia uguale
horikawa}.
\end{proof}

\bigskip




\section{Deformations with fixed
domain and codomain}\label{sezio fixed target domain}

Let $f:X \lrg Y$ be a  holomorphic map and consider the
infinitesimal deformations of $f$ with fixed codomain and domain
(see Definition~\ref{def infin deform f FIXED target domain}).

As we have already observed in Remark~\ref{oss def X e Y fixed
=Gamma cXxY XxY fixed}, these deformations can be interpreted as
 infinitesimal deformations of the graph $\Gamma$ in $X \times
Y$, with $X \times Y$ fixed.

Therefore, in diagram (\ref{diagram L M N in chapter
semiregolar}), we do not  need to consider the DGLA $N=KS_X \times
KS_Y$ and $g$. Thus, the functor $\Def_{(h,g)} $ reduces to
$\Def_h$ with
$$
h :L =A_{X \times Y}^{0,*}(-log\, \Gamma)\hookrightarrow M=
KS_{X\times Y}.
$$
This implies that Theorem~\ref{teo Def_(h,g) Def (f)} reduces to:

\begin{cor}
Let $f:X \lrg Y$  be a  holomorphic map. Then the functor $\Def (X
\stackrel{f}{\lrg} Y)$ of infinitesimal  deformations of $f$ with
fixed codomain and domain is isomorphic to $\Def_h$:
$$
 Def_h \cong \Def (X \stackrel{f}{\lrg} Y).
$$
\end{cor}

\begin{proof}
Apply Theorem~\ref{teo Def_(h,g)  Def (f)} with $N=g=0$.
\end{proof}

\begin{oss}
This corollary is equivalent to \cite[Theorem~5.2]{bib
ManettiPREPRINT}.
\end{oss}

\subsection{Semiregularity map}

As  to the obstructions to deform  a map $f$ fixing both $X$ and
$Y$, Lemma~\ref{lemma ostruzione nel ker tra H^2} has the form
below.

\begin{lem}\label{lem ostr kern X e Y fixed}
If the differential graded vector spaces $(\de A_Z,\debar)$ and
$(\de A_\Gamma,\debar)$  are acyclic, then the functor
$\Def_{\eta}$ is unobstructed. In particular, the obstruction
space of $\Def_{h}$ is naturally contained in the kernel of the
map
$$
H^2(C^\cdot_{h}) \stackrel{\mathcal{I}}{\lrg}
H^2(C^\cdot_{\eta})\cong \oplus_i \Hom\left( H^i(I_\Gamma \cap
\ker(\de)),H^i(\frac {A_Z}{\de A_Z}) \right).
$$

\end{lem}

\begin{proof}
See \cite[Lemma~8.2]{bib ManettiPREPRINT}. It follows from
Lemma~\ref{lemma ostruzione nel ker tra H^2}, with $N=J=0$.

\end{proof}

We note that $H^2(C^\cdot_{h})\cong H ^1 (X,f^*\Theta_Y)$ (see
Remark~\ref{oss ho mappa da C_h ->coker (h)}). Moreover, we have
an analogue of Theorem~\ref{teo ostruzi in ker H^2} and then
Corollary~\ref{coroll semiregolari-ostruzion} becomes:

\begin{cor}\label{cor OSTRUZIO embd uguale a MANE}
Let $f:X \lrg Y$ be a  holomorphic map of  compact K\"{a}hler
manifolds. Let $p=dim Y - dim X$. Then the obstruction space to
the infinitesimal deformations of $f$,  with fixed $X$ and $Y$, is
contained in the kernel of the map
$$
\sigma :  H ^1 (X,f^*\Theta_Y) \lrg H^{p+1}(Y,\Omega_Y^{p-1}).
$$
\end{cor}

\begin{proof}
See \cite[Corollary~9.2]{bib ManettiPREPRINT}. It follows from
 Lemma~\ref{lem ostr kern X e Y fixed} and
Corollary~\ref{coroll semiregolari-ostruzion}.

\end{proof}

\section{Deformations of the inclusion map}
\label{sezio inclusion}

In this section, we focus our attention on the infinitesimal
deformations of an inclusion $i:X \hookrightarrow Y$ of compact
complex  manifolds. The DGLAs $L,M,N$ and the morphisms $g,h$ and
$\pi$ are as before.

Consider the DGLA $L'$ introduced in Section~\ref{sezio
submanifold defz L'} (see also \cite[Section~5]{bib
ManettiPREPRINT})
$$
0 \lrg L' \lrg A_Y^{0,*}(\Theta_Y) \stackrel{\pi'} \lrg
A_X^{0,*}(N_{X|Y})\lrg 0.
$$

\begin{cor}\label{cor L governa def coppia (XcY)}
$L'$ controls the infinitesimal deformations of the inclusion $i:X
\hookrightarrow Y$.
\end{cor}

The proof is postponed to the end of this section after two
preliminary lemmas.

\begin{lem}\label{lem i:X-Y inclu allora h-g suriet}
If $i:X \hookrightarrow Y$ is the inclusion, then the morphism
$g-h:L  \times N \lrg M$ is surjective.
\end{lem}
\begin{proof}
We want to prove that for each $\phi \in M=A_{X\times
Y}^{0,*}(\Theta_{X\times Y})$ there exist $n_1 \in
A_X^{0,*}(\Theta_X)$ and $n_2 \in  A_Y^{0,*}(\Theta_Y)$ such that
$g(n_1,n_2)-\phi =p^*n_1 +q^*n_2 - \phi \in L=\ker \pi$, that is,
$\pi (\phi)=i_*n_1 -i^*n_2\in A_X^{ 0,* }(i^*\Theta_Y)=A_X^{ 0,*
}({\Theta_Y}_{|X})$. Then the proof immediately follows  from the
fact that the restriction morphism $i^* :A_Y^{0,*}(\Theta_Y) \lrg
A_X^{ 0,* }({\Theta_Y}_{|X})$ is surjective.
\end{proof}

\begin{lem}\label{lem se inclusione  Lx_M N=L'}
With the notation above, $L\times_M N \cong L'$.
\end{lem}
\begin{proof}
By definition, $L \times_M N=\{(l,n_1,n_2)\in L \times KS_X \times
KS_Y \ | \ h(l)=p^*n_1 +q^*n_2  \}$ and so $0= \pi h(l)=i_* n_1 -
i^* n_2$.

Define $\gamma:L \times_M N \lrg KS_Y$ as the projection onto
$KS_Y$.

Then $\gamma$ is an injective morphism of DGLAs with $L'$ as
image.

Actually, suppose that $\gamma(l,n_1,n_2)=n_2=0$; then $i_* n_1=0$
and so  $n_1=n_2=l=0$.

As to the image, consider the following exact sequences
\begin{center}
$\xymatrix{ 0\ar[r] & L' \ar[r] \ar[d]&  A_Y^{0,*}(\Theta_Y)
\ar[r]^{\pi'} \ar[d]^{i^*} & A_X^{0,*}(N_{X|Y})\ar[d]_{\cong}
\ar[r]&0  \\
0 \ar[r]& A_X^{0,*}(\Theta_X)\ar[r]^{i_*} &
A_Y^{0,*}({\Theta_Y}_{|X}) \ar[r]^\beta &A_X^{0,*}(N_{X|Y}) \ar[r]
& 0.
\\ }$
\end{center}
Let $ \gamma(l,n_1,n_2)=n_2 \in KS_Y$; then $\beta(i^* n_2)=\beta
(i_* n_1)=0$ and so  $\pi'(n_2)=0$. This implies that
$\gamma(l,n_1,n_2) \in L'$.
\end{proof}

\begin{proof}[Proof of Corollary~\ref{cor L governa def coppia
(XcY)}] By Theorem~\ref{teo Def_(h,g)  Def (f)}, the infinitesimal
deformation functor $\Def(i)$ is isomorphic to the functor
$\Def_{(h,g)}$. By Lemma~\ref{lem i:X-Y inclu allora h-g suriet}
and Proposition~\ref{prop if g-h surj then EXTE DEF_H qiso
DEf_(h,g)}, $\Def_{L \times_M N} \cong \Def_{(h,g)}$. Finally,
  Lemma~\ref{lem se inclusione  Lx_M N=L'} implies that
$\Def_{L'}\cong\Def_{L \times_M N}$. Then $\Def(i) \cong
\Def_{L'}$.
\end{proof}

\subsection{Example}

We can generalize the case of one inclusion $i: X \hookrightarrow
Y$ considering more subvarieties.

For example, let $X$ be a manifold of dimension $n$ and $D_1,
\ldots , D_m$ smooth hypersurfaces, with $0 < m \leq n-2$.
Moreover, assume that $D_1, \ldots , D_m$ intersect transversally
in a smooth subvariety $S$.

Define
$$
\Def_{X;D_1, \ldots , D_m} : \Art \lrg \Set
$$
as the functor of infinitesimal deformations of the holomorphic
map
$$
f:\bigcup^{\circ}D_i \lrg X,
$$
where $f_{|D_i}$ is the inclusion. Equivalently, for any $A \in
\Art$, $\Def_{X;D_1, \ldots , D_m}(A)$ provides   an infinitesimal
deformation $X_A$ of $X$ over $\Spec(A)$ and an infinitesimal
deformations $\mathcal{D}_i \subset X_A$ of the $D_i$'s.

\medskip

Let $\Theta_X(-log\,D)\subset \Theta_X$ be the subsheaf of vector
fields that are tangent to $D_i$, for every $i$, and $N_{D_i|X}$
be the normal bundle of $D_i$ in $X$.

Define $L':=A_X^{0,*}(\Theta_X(-log\, D))$ as in
Section~\ref{sezio submanifold defz L'} and $D=\cup _i D_i \subset
X$. Then $L'$ is a DGLA and we have the following exact sequence
$$
0 \lrg L' \lrg A_X^{0,*}(\Theta_X)\stackrel{\pi'}{\lrg}
A_D^{0,*}(\oplus _i
 N_{D_i|X} )\lrg 0.
$$

Denoting for convenience $\displaystyle D^\circ=
\bigcup^{\circ}D_i$, we define  $ M=KS_{D^\circ\times X}$,
$N=KS_{D^\circ} \times KS_X$ and the morphism
$$
g=(p^*,q^*):N=KS_{D^\circ}  \lrg M=KS_{D^\circ\times X},
$$
where $p$ and $q$ are the projections of $D^\circ \times X$ onto
$D^\circ$ and $X$, respectively.

Finally, let $L$ be the DGLA defined, as in Section~\ref{sez teo
Def_(h,g)Def (f)}, by the following exact sequence:
$$
0 \lrg L \stackrel{h}{\lrg} M \stackrel{\pi}{\lrg}
A_{D^\circ}^{0.*}(f^*\Theta_X)\lrg 0.
$$

\begin{cor}
$L'$ controls the functor $\Def_{X;D_1, \ldots , D_m} $. In
particular, the tangent space of $\Def_{X;D_1, \ldots , D_m}$ is
$H^1(X,\Theta_X(-log\,D))$ and the obstruction space is naturally
contained in $H^2(X,\Theta_X(-log\,D))$.
\end{cor}

\begin{proof}
With the notation above, applying Theorem~\ref{teo Def_(h,g)  Def
(f)}, there exists an isomorphism of functors $\Def_{X;D_1, \ldots
, D_m} \cong \Def_{(h,g)}$. Then proceeding as in the case of
inclusion it   suffices to prove the two steps below.
\begin{itemize}
  \item[$Step\  1$.] The morphism $g-h:L\times N\lrg M$ such that
 $(g-h)(l,n)=g(n)-h(l)$ is surjective (analogous of
 Lemma~\ref{lem i:X-Y inclu allora h-g suriet}).
  \item[$ Step \ 2$ .] $L' \cong L \times_M N$
  (analogous of Lemma~\ref{lem se inclusione  Lx_M N=L'}).
\end{itemize}
Actually, $Step\ 1$ and Proposition~\ref{prop if g-h surj then
EXTE DEF_H qiso DEf_(h,g)} imply that $\Def_{X;D_1, \ldots , D_m}
\cong \Def_{(h,g)} \cong \Def_{L \times_M N}$. Finally, $ Step \
2$ implies $\Def_{X;D_1, \ldots , D_m} \cong \Def_{L'}$.

\smallskip

\emph{Proof of Step 1.} We want to prove that for each $\phi \in
M=A_{D^\circ \times X}^{0,*}(\Theta_{D^\circ \times X})$ there
exist $n_1 \in KS_{D^\circ }=A_{D^\circ }^{0,*}(\Theta_{D^\circ
})$ and $n_2 \in A_X^{0,*}(\Theta_X)$ such that $g(n_1,n_2)-\phi
=p^*n_1 +q^*n_2 - \phi \in L=\ker \pi$, or equivalently  $\pi
(\phi)=f_*n_1 -f^*n_2\in A_{D^\circ }^{ 0,* }(f^*\Theta_X) $. Then
we have to prove that $(f^*,-f_*):A_X^{0,*}(\Theta_X) \times
A_{D^\circ }^{ 0,* }(\Theta_ {D^\circ }) \lrg A_{D^\circ }^{ 0,*
}(f^*\Theta_X)$ is surjective and  this follows from the
hypothesis on $D_i$.

\emph{Proof of Step 2.} By definition, $L \times_M
N=\{(l,n_1,n_2)\in L \times KS_{D^\circ} \times KS_X   \ | \
h(l)=p^*n_1 +q^*n_2 \}$ and so $0= \pi h(l)=f_* n_1 - f^* n_2$.

Define $\gamma:L \times_M N \lrg KS_X$ as the projection onto
$KS_X$.

Then $\gamma$ is an injective morphism of DGLAs and its image is
$L'=A_X^{0,*}(\Theta_X(-log\, D)) $.

Actually, suppose that $\gamma(l,n_1,n_2)=n_2=0$; then $f_* n_1=0$
and so  $n_1=n_2=l=0$.

As to the image, consider the following exact sequences
\begin{center}
$\xymatrix{ 0\ar[r] & L' \ar[r] \ar[d]&  A_X^{0,*}(\Theta_X)
\ar[r]^{\pi'} \ar[d]^{f^*} & A_D(\oplus _i
 N_{D_i|X} ) \ar[d]_{\cong}
\ar[r]&0  \\
0 \ar[r]& A_{D^\circ }^{ 0,* }(\Theta_ {D^\circ })\ar[r]^{f_*} &
A_{D^\circ }^{ 0,* }(f^*\Theta_X) \ar[r]^\beta & A_D(\oplus _i
 N_{D_i|X} ) \ar[r] & 0.
\\ }$
\end{center}
Let $ \gamma(l,n_1,n_2)=n_2 \in KS_X$, then $\beta(f^* n_2)=\beta
(f_* n_1)=0$ and so  $\pi'(n_2)=0$. This implies that
$\gamma(l,n_1,n_2) \in L'$.

\end{proof}

\subsection{Semiregularity map for the inclusion} Let $i:X
\hookrightarrow Y$ be the inclusion of a submanifold $X$ in $Y$
and  $L'$   the DGLA (defined in Section~\ref{sezio submanifold
defz L'}):
$$
0 \lrg L' \lrg A_Y^{0,*}(\Theta_Y) \stackrel{\pi'} \lrg
A_X^{0,*}(N_{X|Y})\lrg 0.
$$
In Corollary~\ref{cor L governa def coppia (XcY)}, we proved that
$L'$ controls the infinitesimal deformations  of the inclusion
$i$. In particular, this implies that the obstructions are
naturally contained in $H^2(L')$.

Moreover, as in Section~\ref{sezio HTP},  we can consider the
following morphism of DGLAs
$$
\bi: L' \lrg  K'= \left\{f\in \Htp\left(\ker(\de),\dfrac{A_Y}{\de
A_Y}\right)\mid f(I_X\cap\ker(\de))\subset \dfrac{I_X}{I_X\cap\de
A_Y}\right\}.
$$

Then we get the following corollary, whose proof is essentially
contained in Manetti \cite[Corollary~9.2]{bib ManettiPREPRINT}.

\begin{cor}\label{cor semiregularity inclusion}
Let $i:X \lrg Y$ be the inclusion of a submanifold $X$ in a
compact K\"{a}hler manifold $Y$. Let $p=dim Y - dim X$. Then the
obstruction space to the infinitesimal deformations of $i$  is
contained in the kernel of the map
$$
\sigma :  H ^2 (L') \lrg H^{p-1,p+2}(Y,I_X),
$$
where $I_X=\ker i^* \subset A_Y^{*,*}$ is the subcomplex of
differential forms vanishing on $X$.
\end{cor}

\begin{proof}
It follows from   Corollary~\ref{cor OSTRUZIO embd uguale a MANE},
by recalling that $ L' \subset \{a \in A^{0,*}_Y(\Theta_Y)| \
\bi_a(I_X)\subset I_X\} $.

\end{proof}


\clearpage

\begin{thebibliography}{99}
\bibitem{bib Artin} M. Artin,
{\sl Deformations of Singularities,} Tata Institute of
Foundamental Research, Bombay, (1976).


\bibitem{bib atyha Mc donald} M.F. Atiyah, I.G. MacDonald,
{\sl Introduction to Commutative Algebra,} Addison-Wesley
Publishing Company, (1969).


\bibitem{bib bloch} S. Bloch,
\emph{Semi-regularity and de Rham cohomology}, Invent. Math {\bf
17}, (1972), 51-66.

\bibitem{bib Catanesecime} F. Catanese,
\emph{Moduli of algebraic surfaces,} Theory of moduli (Montecatini
Terme,1985), Lecture Notes in Mathematics, {\bf 1337},
Springer-Verlag, New York/Berlin, (1988), 1-83.

\bibitem{bib clemens} H. Clemens,
\emph{Geometry of formal Kuranishi theory,}  Adv. Math., {\bf 198}
(No. 1), (2005),   311-365.


\bibitem{bib Fantec-Manet}B. Fantechi, M. Manetti,
\emph{Obstruction calculus for functors of Artin rings, I,} J.
Algebra, {\bf 202}, (1998), 541-576.

\bibitem{bib fiorenz-Manet} D. Fiorenza, M. Manetti,
\emph{$L_{\infty}$ structures on mapping cones.} Preprint
\texttt{arXiv:math.QA/0601312}.


\bibitem{bib fukaya} K. Fukaya,
\emph{Deformation theory, homological algebra and mirror
symmetry,} Geometry and physics of branes (Como, 2001), Ser. High
Energy Phys. Cosmol. Gravit., IOP Bristol, (2003), 121-209.
Electronic version available at
\texttt{http://www.math.kyoto-u.ac.jp/\%7Efukaya/como.dvi}



\bibitem{bib Godement} R. Godement,
{\sl Topologie alg\'{e}brique et th\'{e}orie des
faisceaux,}Actualités scientifiques et industrielles, {\bf 1252},
Hermann, Paris, (1958).



\bibitem{bib GoMil III} W.M. Goldman, J.J. Millson,
\emph{The homotopy invariance of the Kuranishi space,} Ill.  J.
Math., {\bf34}, (1990), 337-367.



\bibitem{bib Griffone} Ph. Griffiths, J. Harris,
{\sl Principles of algebraic geometry,} Wiley Classics Library,
Wiley, New York, (1978).



\bibitem{bib Grothendieck} A. Grothendieck,
\emph{Technique de Descente et th\'{e}or\`{e}mes d'existence en
g\'{e}om\'{e}trie alg\'{e}brique. II. Le  th\'{e}or\`{e}mes
d'existence en th\'{e}orie formelle des modules}, S\'{e}minaire
Bourbaki, t. 12, Exp. no. {\bf 195}, (1959-1960).


\bibitem{bib hartshorne} R. Hartshorne,
{\sl Algebraic Geometry,} Graduate texts in mathematics, {\bf 52},
Springer-Verlag, New York/Berlin, (1977).



\bibitem{bib HorikawaI} E. Horikawa,
\emph{On deformations of holomorphic maps I, II, } J. Math. Soc.
Japan, {\bf 25} (No.3), (1973), 372-396; {\bf 26} (No.4), (1974),
647-667.

\bibitem{bib Horikawa III} E. Horikawa,
\emph{On deformations of holomorphic maps III, }   Math. Annalen,
{\bf 222}, (1976), 275-282.


\bibitem{bib kodaira libro} K. Kodaira,
{\sl Complex Manifolds and Deformation of Complex Structures,} Die
Grundlehren der mathematischen Wissenschaften in
Einzeldarstellungen, {\bf 283}, Springer-Verlag, New York/Berlin,
(1986).


\bibitem{bib Kodaira SpencerII} K. Kodaira, D.C. Spencer,
\emph{On Deformations of Complex Analytic Structures, II}  Ann. of
Math., \textbf{67} (2), (1958), 403-466.




\bibitem{bib kontsevich} M. Kontsevich,
\emph{Deformation quantization of Poisson manifolds, I.} Letters
in Mathematical Physics, {\bf 66}, (2003), 157-216;
\texttt{arXiv:q-alg/9709040}.


\bibitem{bib kuranishi} M. Kuranishi,
{\sl Deformations of compact complex manifolds,} S\'eminaire de
Mathematiques Sup\'erieures, No. \textbf{39} (\'Et\'e 1969), Les
Presses de l'Universit\'e de Montreal, Montreal, (1971).



\bibitem{bib manetPISA} M. Manetti,
\emph{Deformation theory via differential graded Lie algebras,}
{\sl Seminari di Geometria Algebrica 1998-1999,} Scuola Normale
Superiore (1999).

\bibitem{bib Manetti IMNR} M. Manetti,
\emph{Extended deformation functors.} Int. Math. Res. Not.,
{\bf14}, (2002), 719-756; \texttt{arXiv:math.AG/9910071}.


\bibitem{bib mane COSTRAINT} M. Manetti,
\emph{Cohomological constraint to deformations of compact K\"ahler
manifolds.}, Adv. Math.,  {\bf 186}, (2004), 125-142;
\texttt{arXiv:math.AG/0105175}.

\bibitem{bib manRENDICONTi} M. Manetti,
\emph{Lectures on deformations of complex manifolds,} Rend. Mat.
Appl. (7), {\bf24}, (2004), 1-183; \texttt{arXiv:math.AG/0507286}.

\bibitem{bib ManettiPREPRINT} M. Manetti,
\emph{Lie description of higher obstructions to deforming
submanifolds,} Preprint \texttt{arXiv:math.AG/0507287} v2
7Oct2005.

\bibitem{bib Manetti PRODOTTI}  M. Manetti,
\emph{Deformations of products and new examples of obstructed
irregular surfaces,} Preprint \texttt{arXiv:math.AG/0601169}.



\bibitem{bib Matsumura} H Matsumura,
{\sl Commutative ring theory,} Cambridge Studies in Advanced
Mathematics, {\bf 8}, Cambridge Univesity Press, Cambridge,
(1986).




\bibitem{bib Namba} M. Namba,
{\sl Families of meromorphic functions on compact Riemann
surfaces,} Lecture Notes in Mathematics, {\bf 767},
Springer-Verlag, New York/Berlin, (1979).


\bibitem{bib RAn mappe}Z. Ran,
\emph{Deformations of maps}, Algebraic curves and projective
geometry (Trento 1989),  Lecture Notes in Mathematics, {\bf 1389},
Springer-Verlag, New York/Berlin, (1989), 246-253.


\bibitem{bib Ran ostruzio}Z. Ran,
\emph{Semiregularity, obstructions and deformations of Hodge
classes}, Ann. Scuola Norm. Pisa Cl. Sci., {\bf 28} (No. 4),
(1999), 809-821.

\bibitem{bib Ran ostruzioII}Z. Ran,
\emph{Universal variations of Hodge structure an
Calabi-Yau-Schottky relations}, Invent. Math., {\bf 138}, (1999),
425-449.





\bibitem{bib Sernesi} E. Sernesi,
{\sl Deformations of Algebraic Schemes}, Grundlehren der
mathematischen Wissenschaften, {\bf 334}, Springer-Verlag, New
York/Berlin, 2006. Electronic version available at
\texttt{http://www.mat.uniroma3.it/users/sernesi/defsch.html}



\bibitem{bib Schlessinger }M. Schlessinger,
\emph{Functors of Artin rings, } Trans. Amer. Math. Soc., {\bf
130}, (1968), 208-222.

\bibitem{bib Voisin} C. Voisin,
{\sl Th\'{e}orie de Hodge et g\'{e}om\'{e}trie alg\'{e}brique
complexe I,} Soci\'{e}t\'{e} Math\'{e}matique de France, Paris,
(2002).


\end{thebibliography}
\end{document}